\documentclass[reqno,11pt,twoside]{amsart}
\usepackage{amlt}
\RequirePackage[none]{hyphenat}

\usepackage[bookmarks]{hyperref}
\hypersetup{%
    naturalnames,%
    bookmarksnumbered,%
    plainpages=false,%
    pdfstartview=FitH,%
    pdfborderstyle={/S/U/W 1},%
    pdfdisplaydoctitle=true,%
    pdftitle={Log Structures on Generalized Semi-Stable Varieties},%
    pdfauthor={Li Ting},%
    pdfsubject={Algebraic Geometry},%
    pdfkeywords={Log structure, normal crossing singularity, semistable variety}%
}

\numberwithin{equation}{section}
\newcommand{\lC}{local chart\xspace}
\newcommand{\lCs}{local charts\xspace}
\newcommand{\rlC}{refined local chart\xspace}
\newcommand{\rlCs}{refined local charts\xspace}
\newcommand{\RlCs}{Refined Local Charts}
\newcommand{\wnC}{weakly normal crossing\xspace}
\newcommand{\cpC}{\mathrm{CP}}

\newcommand{\Mor}[1]{\xymatrix@1@C-.7em{#1}}
\newcommand{\bMor}[2]{\Mor{#1 \ar@<.6ex>[r] \ar@<-.2ex>[r] & \hspace{0.1em}#2}}

\newcommand{\CartS}{\ar@{}[dr]|-*+{\square}}

\newcommand{\rSch}[2]{#1_{[#2]}}

\newcommand{\prM}[1]{\mathrm{pr}_{#1}}

\begin{document}

\clearP
\pdfbookmark[1]{Abstract}{amltb}
\begin{center}
\large \bfseries Log Structures on Generalized Semi-Stable Varieties

\bigskip

{\normalfont\normalsize Ting Li\\
Directed by Professor Zhao Chunlai}
\plainfoot{April 12, 2005}%
\plainfoot{Department of Mathematics, Peking University, Beijing 100871, P. R. China}%
\plainfoot{EMail: moduli@live.cn}%

\bigskip
Abstract
\bigskip
\end{center}

In this paper we study the log structures on generalized semistable varieties, generalize the result by F.~Kato~and~M.~Olsson, and prove the canonicity of log structure when it can be expected.

In out text we first give the definitions of \lC and \wnC morphism. Then we study the invariants of complete noetherian local ring coming from \wnC morphisms. These invariants enable us to further define the \rlCs and prove that all log structures induced by \rlCs are locally isomorphic. Let $f \colon X \to S$ be a surjective, proper and \wnC morphism of locally noetherian schemes which satisfies the conditions (\dag) and (\ddag) in \ref{Se:global:3} and certain local conditions stated at the beginning of \ref{Se:loc}. Then the obstructions for the existence of semistable log structures on $X$ is an invertible sheaf $\mathscr{L}(f)$ on a finite \mDash{X}{scheme} $E = E(f)$. The main result of local case with respect to base schemes is:

\newtheorem*{enT}{Theorem}
\begin{enT}\
\begin{enumerate}
\item There exists a semistable log structure on $X$ if and only if $\mathscr{L}(f) \cong \mathcal{O}_E$.
\item The semistable log structure on $X$ is unique up to (not necessarily canonical) isomorphisms if it exists.
\end{enumerate}
\end{enT}

The main result of global case with respect to base schemes is:

\begin{enT}
Let $X$ and $S$ be locally noetherian schemes, $f \colon X \to S$ a surjective proper \wnC morphism without powers. If $f$ satisfies the condition \rmB{\dag} in \ref{Se:global:3} and every fiber of $f$ is geometrically connected, then
\begin{enumerate}
\item There exists a semistable log structure for $f$ if and only if for every point $y \in S$, $\mathscr{L}_{\bar{y}}$ is trivial on $E_{\bar{y}}$.
\item Let $(\mathscr{M}_1, \mathscr{N}_1, \sigma_1, \tau_1, \phi_1)$ and $(\mathscr{M}_2, \mathscr{N}_2, \sigma_2, \tau_2, \phi_2)$ be two semistable log structures for $f$. Then there exist isomorphisms of log structures $\varphi \colon \mathscr{M}_1 \xrightarrow{\sim} \mathscr{M}_2$ and $\psi \colon \mathscr{N}_1 \xrightarrow{\sim} \mathscr{N}_2$ such that $\varphi \circ \phi_1 = \phi_2 \circ f^{\ast}\psi$, $\sigma_2 \circ \varphi = \sigma_1$ and $\tau_2 \circ \psi = \tau_1$. Moreover such pair $(\varphi,\psi)$ is unique.
\end{enumerate}
\end{enT}

We further prove that the existence of semistable log structures remains under fibred products, base extension, inverse limits, flat descent. Finally we study the semistable curves. The main result is:

\begin{enT}
Any semistable curve over a locally noetherian scheme is a \wnC morphism without powers and has a canonical semistable log structure.
\end{enT}


\clearP
\tableofcontents
\clearP

\section*{Introduction}
A major advantage of logarithmic geometry is that it enables us to treat some kind of singularity as smooth case. To achieve this, we must equip a singular morphism with suitable log structure so that it becomes \emph{smooth} in the sense of logarithmic geometry. 

In this paper, we study the existence and uniqueness of semistable log structures on the morphisms of schemes which locally have the form
\[\Spec \frac{A[T_{11},\ldots,T_{1q_1},\ldots,T_{p1},\ldots,T_{pq_p},T_1,\ldots,T_m]}{\Bigl(\prodT\limits^{q_1}_{j_1=1}T_{1j_1}^{e_{1j_1}} - a_1, \ldots, \prodT\limits^{q_p}_{j_p=1}T_{pj_p}^{e_{pj_p}} - a_p\Bigr)} \,.\]
If $p$ and all $e_{ij_i}$ are equal to $1$, then the singularity is the so called \emph{normal crossing} in the classical sense.

The study of normal crossing singularities began with Deligne and Mumford \cite{PDDM}, where they showed that any curve with normal crossing singularities deforms to a smooth curve. For higher dimensional spaces, Friedman \cite{RFri1} discovered that an obstruction for the existence of smoothenings with regular total space is an invertible sheaf on the singular locus. In \cite[\S11-12]{FKato2} and \cite{FKato1}, F.~Kato introduced log structures for normal crossing varieties over fields. And in \cite{MCOL1}, M.~Olsson generalizes them to morphisms $f \colon X \to S$, where $X$ is locally isomorphic to
\[\Spec \mathcal{O}_S[T_1,\ldots,T_l]/(T_1 \cdots T_r - t) \,,\]
with $t \in \varGamma(S,\mathcal{O}_S)$ a fixed section. Also in \cite{FKato3}, F.~Kato considered the existence of log structures on pointed stable curves.

In this paper we generalize the results in \cite{FKato1} and \cite{MCOL1}, mainly add nontrivial powers and remove the fixed section $t$ in \cite{MCOL1}. Roughly speaking, we construct an obstruction at every morphism $X \times_S \Spec \mathcal{O}_{S,\bar{y}} \to \Spec \mathcal{O}_{S,\bar{y}}$. Then we prove that the semistable log structure exists if and only if all these obstructions vanish (see Theorem \thmref{670P2} and \thmref{671P3}). In the case of no power, we shall see that this kind of semistable log structures is canonical (i.e.~unique up to a unique isomorphism), which was not discussed in \cite{FKato1} and \cite{MCOL1}.

In Section 1 we generalize the concept of normal crossing to the so called ``\wnC{}". In Section 2 we study the invariants of complete local rings, which is of fundamental importance. In Section 3 we define the concept of \rlC. On each \rlC, we may define a log structure, which is the tile for building the global semistable log structures. In Section 4 we list some technique and notations in cohomology theory which are needed in later sections.

In Section 5 we study the local case. In other words, for a \wnC morphism $f \colon X \to S$, we focus on morphisms $X_V \to V$ for every \'{e}tale neighborhood $V$ on $S$ which is \emph{small enough}, especially the case when the base scheme $S$ is the spectrum of a strictly Henselian local ring. For \wnC morphisms $f \colon X \to S$ with nontrivial power, the theory can only be built on local cases, because semistable log structures on $X_V \to V$ may not be unique (up to isomorphism). In Section 6 we prove that for a \wnC morphism without powers, the semistable log structures exist if and only if all local obstructions vanish. If so, then it must be canonical.

In Section 7 we study properties of \wnC morphisms under base change. We shall show that the semistable log structures constructed in \S6 have good functorial properties. In Section 8 we show that on semistable curves, our constructed obstructions are always trivial. So there exists a canonical log structure on any semistable curve which make it log smooth.

\bigskip

{\bfseries Notation and Conventions.} Throughout the paper, rings, algebras and monoids are all assumed to be commutative and have multiplicative identity elements. A homomorphism of rings (resp.~monoids) is assumed to preserve identity element. A subring (resp.~monoid) is assumed to contain the identity of the total ring (resp.~monoid).

If $n$ is a positive integer, we use $S_n$ to denote the symmetric group on $\{1,2,\ldots,n\}$.

For every pair of integers $m$ and $n$, we define a set $[m,n]$ as
\[[m,n] \defeq \begin{cases}
\{m,m+1,\ldots,n\} \,, & \mbox{if $m \leqslant n$,} \\
\qquad \emptyset \,, & \mbox{if $m > n$.}
\end{cases}\]

For a field $k$, we use $\bar{k}$ to denote the algebraic closure of $k$ and $k_s$ the separable closure of $k$.

If $X$ is a scheme, $f \in \varGamma(X,\mathcal{O}_X)$ is a section and $x \in X$ is a point, we use $f(x)$ to denote the image of the stalk $f_x$ in the residue field $\kappa(x)$.

If $X$ is a scheme, a geometric point on $X$ is a morphism of schemes $\Spec K \to X$ where $K$ is a separably closed field; if $x$ is a point on $X$, we use $\bar{x}$ to denote the geometric point $\Spec \kappa(x)_s \to X$.

If $S$ is a scheme, $f \colon X \to S$ and $T \to S$ are two \mDash{S}{schemes}, then we define $X_T \defeq X \times_S T$ and let $f_T \colon X_T \to T$ denote the second projection.

For every morphism $f \colon X \to S$ of schemes, we use $\rSch{X}{f}$ to denote the \mDash{S}{scheme} $X$ via $f$.

If $X$ is a scheme and $G$ a monoid (resp.~abelian group), we use $G_X$ to denote the constant sheaf of monoids (resp.~of abelian groups) on $\etSite{X}$ associated to $G$.

If $\mathscr{M}$ is a log structure on a scheme $X$, we write $\overline{\mathscr{M}} \defeq \mathscr{M}/\mathcal{O}_X^{\ast}$.

\section{Definition}
\begin{noteDe}[667E1]
Let $f \colon X \to S$ be a morphism of finite type of locally noetherian schemes, $x$ a point on $X$ and $y \defeq f(x)$. A \emph{\lC} of $f$ at $x$ consists of the following data:
\begin{enumerate}
\item an \'{e}tale \mDash{S}{scheme} $V = \Spec R$ which is a connected affine scheme;
\item a point $y'$ on $V$ which maps onto $y$;
\item an \'{e}{tale} \mDash{V \times_S X}{scheme} $U = \Spec A$ which is a connected affine scheme;
\item a point $x'$ on $U$ which maps onto $x$ and $y'$;
\item a finitely generated \mDash{R}{algebra} $P$ such that $\varOmega_{P/R}$ is a free \mDash{P}{module}, $\Spec P$ is connected and is smooth over $V$;
\item a point $\mathfrak{p}$ in $\Spec P$ which maps onto $y'$;
\item a subset
\[\defset[1]{T_{ij_i}}{i \in [1,p], \ j_i \in [1,q_i]}\]
of $P$ such that $T_{ij_i} \in \mathfrak{p}$ for all $i$ and $j_i$, and
\[\defset[1]{d_{P/R}(T_{ij_i})}{i \in [1,p], \  j_i \in [1,q_i]}\]
is a part of a basis of $\varOmega_{P/R}$;
\item a closed \mDash{V}{immersion} $U \hookrightarrow \Spec P$ which maps $x'$ onto $\mathfrak{p}$ and is defined by the ideal
\[\Bigl(\prodT^{q_1}_{j_1=1}T_{1j_1}^{e_{1j_1}} - a_1, \prodT^{q_2}_{j_2=1}T_{2j_2}^{e_{2j_2}} - a_2, \ldots, \prodT^{q_p}_{j_p=1}T_{pj_p}^{e_{pj_p}} - a_p\Bigr) \,,\]
where $a_i \in R$ such that $a_i(y') = 0$, and $e_{ij_i} \geqslant 1$ are integers which are invertible in $R$, such that for every $i \in [1,p]$, $\sum\limits^{q_i}_{j_i=1}e_{ij_i} > 1$ and
\begin{equation}\label{E:defn1e}
D_i(U/V) \defeq \Spec \Bigl(P \Big/\sumT^{q_i}_{j_i = 1}P \cdot \bigl(T_{i1}^{e_{i1}} \cdots T_{i,j_i-1}^{e_{i,j_i-1}} T_{ij_i}^{e_{ij_i}-1} T_{i,j_i+1}^{e_{i,j_i+1}} \cdots T_{i,q_i}^{e_{i,q_i}} \bigr)\Bigr)
\end{equation}
is connected.
\end{enumerate}

We use
\begin{equation}\label{E:defn0e}
U \to \Spec \Bigl(P \Bigl/ \Bigl(\ldots,\prodT^{q_i}_{j_i=1}T_{ij_i}^{e_{ij_i}} - a_i, \ldots \Bigr)\Bigr)
\end{equation}
or $U/V$, or simply $U$ to denote the \lC.
\end{noteDe}

\begin{noteRe}
Note that all these \emph{connectedness} can be satisfied by contracting $\Spec P$, $U$ and $V$ suitably, so they are not essential restriction.
\end{noteRe}

The following theorem shows that if a point has a \lC, then all points in some of its open neighborhood have \lCs.

\begin{noteTh}[667E0]
Let $R$ be a noetherian ring, $P$ a finitely generated \mDash{R}{algebra},
\[\defset[1]{T_{ij_i}}{i \in [1,p], \ j_i \in [1,q_i]}\]
a subset of $P$, $\Sequ{p}{a} \in R$,
\[\defset[1]{e_{ij_i}}{i \in [1,p], \ j_i \in [1,q_i]}\]
a set of positive integers which are invertible in $R$. For each $i \in [1,p]$, put
\[b_i \defeq \prodT^{q_i}_{j=1}T_{ij}^{e_{ij}} - a_i \,.\]
Put $A \defeq P/(\Sequ{p}{b})$, $S \defeq \Spec R$ and $X \defeq \Spec A$. Assume that
\begin{enumerate}
\renewcommand{\theenumi}{\alph{enumi}}
\item $P$ is smooth over $R$;
\item $\varOmega_{P/R}$ is a free \mDash{P}{module};
\item $\defset[1]{d_{P/R}(T_{ij_i})}{i \in [1,p], \, j_i \in [1,q_i]}$ is a part of a basis of $\varOmega_{P/R}$;
\item for any $i \in [1,q_i]$, $\sum\limits^{q_i}_{j_i=1}e_{ij_i} > 1$.
\end{enumerate}
Then we have
\begin{enumerate}
\item $\Sequ{p}{b}$ is a \mDash{P}{regular} sequence.
\item $X \to S$ is a flat and local complete intersection morphism.
\item For every point $x$ on $X$, there is a \lC at $x$.
\end{enumerate}
\end{noteTh}

\begin{proof}
(1) and (2). Since $P$ is smooth over $R$ and $\{\cdots,d(T_{ij_i}),\cdots\}$ is a part of a basis of $\varOmega_{P/R}$, $\{\ldots,T_{ij_i},\ldots\}$ are algebraically independent over $R$ and $P$ is smooth over $R[\cdots,T_{ij_i},\cdots]$. So we may assume that
\[P = R[\cdots,T_{ij_i},\cdots]\]
is a polynomial algebra over $R$ with indeterminates $\{\ldots,T_{ij_i},\ldots\}$. Then (1) is by \cite[(20.F), COROLLARY 2]{Mats1} and induction on $p$. So $X \to S$ is a local complete intersection morphism. By \cite[Corollary of Theorem 22.5]{Mats2}, $X$ is flat over $S$.

(3) $x$ defines a prime ideal $\mathfrak{P}$ of $P$. Put $\mathfrak{p} \defeq R \cap \mathfrak{P}$. Assume that $a_i \in \mathfrak{p}$ for $i \in [1,l]$ and $a_i \notin \mathfrak{p}$ for $i \in [l+1,p]$. And for each $i \in [1,l]$, we assume that $T_{ij_i} \in \mathfrak{P}$ for $j_i \in [1,s_i]$ and $T_{ij_i} \notin \mathfrak{P}$ for $j_i \in [s_i+1,q_i]$. Obviously, for all $i \in [1,l]$ we have $s_i \geqslant 1$. Assume that $\sum\limits^{s_i}_{j=1}e_{ij} > 1$ when $i \in [1,r]$, and $\sum\limits^{s_i}_{j=1}e_{ij} = 1$ when $i \in [r+1,l]$. By taking an affine open neighborhood of $\mathfrak{P}$ in $\Spec P$ and an affine open neighborhood of $\mathfrak{p}$ in $\Spec R$, we may assume that $a_i \in R^{\ast}$ whenever $i \in [l+1,p]$, and $T_{ij_i} \in P^{\ast}$ whenever
\[i \in [l+1,p] \vee \bigl(i \in [1,l]  \wedge j \in [s_i+1,q_i] \bigl)\]
is valid. Then
\[P' \defeq P/(b_{r+1},\ldots,b_p)\]
is smooth over $R$. For each $i \in [1,r]$, since $e_{i1}$ is invertible in $\SH{P}_{\mathfrak{P}}$, there exists an element $u_i \in \SH{P}_{\mathfrak{P}}$ such that
\[u_i^{e_{i1}} = \prod^{q_i}_{j_i = s_i+1}T_{ij_i}^{e_{ij_i}} \qquad \mbox{(if $s_i = q_i$, we let $u_i = 1$)} \,.\]
By taking an affine \'{e}tale neighborhood of $\mathfrak{P}$ in $\Spec P$, we may assume that $u_i \in P$ for all $i \in [1,r]$. For each $i \in [1,r]$, let $T'_{i1}$ be the image of $u_iT_{i1}$ in $P'$, and for each $j_i \in [2,s_i]$, $T'_{ij_i}$ be the image of $T_{ij_i}$ in $P'$. Then we have
\[A = P' \Big/ \Bigl(\prodT^{s_1}_{j_1=1}(T'_{1j_1})^{e_{1j_1}} - a_1, \prodT^{s_2}_{j_2=1}(T'_{2j_2})^{e_{2j_2}} - a_2, \ldots, \prodT^{s_r}_{j_r=1}(T'_{rj_r})^{e_{rj_r}} - a_r\Bigr) \,.\]
Moreover, $P'$ is smooth over $R$ and $\{\cdots,d(T'_{ij_i}),\cdots\}$ is a part of basis of $\varOmega_{P'/R}$.
\end{proof}

\begin{noteDe}
Let $f \colon X \to S$ be a morphism of locally noetherian schemes. We say that $f$ is \emph{\wnC} if it is of finite type, and for every point $x \in X$, either $f$ is smooth at $x$ or there exists a \lC at $x$.

A \wnC morphism $f \colon X \to S$ is said to be \emph{without powers} if in every \lC of $f$ as \eqref{E:defn0e}, all the powers $e_{ij_i}$ are equal to $1$.
\end{noteDe}

By Theorem \thmref{667E0}, if $f \colon X \to S$ is \wnC, then $f$ is a flat and local complete intersection morphism.

The following lemma is obvious.

\begin{noteTh}
Let
\[\xymatrix{X' \ar[r]^{f'} \ar[d] \ar@{}[dr]|-*+{\square} & S' \ar[d] \\ X \ar[r]_{f} & S}\]
be a Cartesian square of locally noetherian schemes. If $f$ is \wnC, so is $f'$.
\end{noteTh}

\section{Invariants of Complete Local Rings}
In this section we study the invariants of complete noetherian local ring coming from \wnC morphisms, which ensure that all log structures induced by \lCs are locally isomorphic.

Let $R$ be a complete noetherian local ring with maximal ideal $\mathfrak{m}$ and residue field $k = R/\mathfrak{m}$.

Let $P$ and $Q$ be rings of power series over $R$ in variables
\[\defset[1]{X_{ij}}{i \in [1,p], \  j \in [1,q_i]} \cupT\, \bigl\{\Sequ{m}{X}\bigr\}\]
and
\[\defset[1]{Y_{i'j'}}{i' \in [1,p'], \  j' \in [1,q'_{i'}]} \cupT\, \bigl\{\Sequ{m'}{Y}\bigr\}\]
respectively.

For each $i \in [1,p]$, let $e_{i1},e_{i2},\ldots,e_{iq_i}$ be positive integers which are invertible in $R$ with $\sum\limits^{q_i}_{j=1}e_{ij} > 1$, $a_i$ an element in $\mathfrak{m}$, and
\[F_i \defeq \prod^{q_i}_{j=1}X_{ij}^{e_{ij}} - a_i \in P \,.\]

For each $i' \in [1,p']$, let $e'_{i'1},e'_{i'2},\ldots,e'_{i'q'_{i'}}$ be positive integers which are invertible in $R$ with $\sum\limits^{q'_{i'}}_{j'=1}e'_{i'j'} > 1$, $b_{i'}$ an element in $\mathfrak{m}$, and
\[G_{i'} \defeq \prod^{q'_{i'}}_{j'=1}Y_{i'j'}^{e'_{i'j'}} - b_{i'} \in Q \,.\]

Put
\[A \defeq P/(\Sequ{p}{F}) \qquad \mbox{and} \qquad B \defeq Q/(\Sequ{p'}{G}) \,.\]
Let $x_{ij}$, $x_k$ and $y_{i'j'}$, $y_{k'}$ be the images of $X_{ij}$, $X_k$ and $Y_{i'j'}$, $Y_{k'}$ in $A$ and $B$ respectively. Let $\mathfrak{M}_1$ and $\mathfrak{M}_2$ be the maximal ideals of $A$ and $B$, $\mathfrak{N}_1$ and $\mathfrak{N}_2$ the nilradicals of $A$ and $B$.

The following theorem is the main result of this section.

\begin{noteTh}[667T1]
Let $\varphi \colon A \xrightarrow{\sim} B$ be an isomorphism of \mDash{R}{algebras}. Then $p = p'$, $m = m'$; and there exists a $\sigma \in S_p$ such that for each $i \in [1,p]$, we have
\begin{enumerate}
\item $q_i = q'_{\sigma(i)}$,
\item $a_i = u_ib_{\sigma(i)}$ for some $u_i \in R^{\ast}$,
\item there exists a $\tau_i \in S_{q_i}$ such that for each $j \in [1,q_i]$, we have $e_{ij} = e'_{\sigma(i),\tau_i(j)}$ and $\varphi(x_{ij}) = v_{ij}y_{\sigma(i),\tau_i(j)}$ for some $v_{ij} \in B^{\ast}$.
\end{enumerate}
\end{noteTh}

To prove Theorem \thmref{667T1}, we note the following simple fact.

\begin{noteLe}[667L2]
Every element in $A$ can be uniquely written as a power series:
\begin{equation}\label{E:local0}
\sum c(\ldots,\alpha_{ij},\ldots;\ldots,\beta_k,\ldots) \cdot \left(\prod^p_{i=1}\prod^{q_i}_{j=1}x_{ij}^{\alpha_{ij}}\right) \cdot \left(\prod^m_{k=1}x_k^{\beta_k}\right) \,,
\end{equation}
where $\alpha_{ij}$, $\beta_k$ are in $\mathbb{N}$, $c(\cdots)$ are in $R$ satisfying the following conditions: for every $i \in [1,p]$, there exists a $j \in [1,q_i]$ such that $\alpha_{ij} < e_{ij}$. \rmB{So we may talk \emph{monomials} and \emph{coefficients} etc.} Furthermore, if
\[a_1 = a_2 = \cdots = a_p = 0 \,,\]
then $A = \bigoplus A_n$ is a graded ring, where $A_n$ consists of homogeneous polynomials of degree $n$.
\end{noteLe}

We first prove Theorem \thmref{667T1} in a simple but fundamental case.

\begin{noteLe}[667L3]
If $R = K$ is a field, then Theorem \thmref{667T1} is valid.
\end{noteLe}

\begin{proof}
Without lose of generality, we may assume that
\[q_i \begin{cases}
> 1 \,, & i \in [1,r] \,, \\ =1 \,, & i \in [r+1,p] \,;
\end{cases}
\qquad \mbox{and} \qquad
q'_{i'} \begin{cases}
> 1 \,, & i' \in [1,r'] \,, \\ =1 \,, & i' \in [r'+1,p'] \,,
\end{cases}\]
for some $r \in [0,p]$ and $r' \in [0,p']$. Firstly it easy to see that
\begin{align*}
\mathfrak{N}_1 &= \left(\prodT^{q_1}_{j=1}x_{1j}, \prodT^{q_2}_{j=1}x_{2j}, \ldots, \prodT^{q_p}_{j=1}x_{pj}\right) \,, \\
\mathfrak{N}_2 &= \left(\prodT^{q'_1}_{j'=1}y_{1j'}, \prodT^{q'_2}_{j'=1}y_{2j'}, \ldots, \prodT^{q'_{p'}}_{j'=1}y_{p'j'}\right) \,.
\end{align*}
Note that $\varphi$ induces an isomorphism of vector spaces over $K$:
\[\bar{\varphi} \colon \mathfrak{N}_1/(\mathfrak{N}_1 \cap \mathfrak{M}_1^2) \xrightarrow{\sim} \mathfrak{N}_2/(\mathfrak{N}_2 \cap \mathfrak{M}_2^2) \,.\]
As $\mathfrak{N}_1/(\mathfrak{N}_1 \cap \mathfrak{M}_1^2)$ has a base $\bar{x}_{r+1,1},\ldots,\bar{x}_{p1}$ and $\mathfrak{N}_2/(\mathfrak{N}_2 \cap \mathfrak{M}_2^2)$ has a base $\bar{y}_{r'+1,1},\ldots,\bar{y}_{p'1}$, we have
\begin{equation}\label{E:local9}
p-r = p'-r' = f
\end{equation}
and there is a $D = (d_{i'i}) \in \GL{f}{K}$ such that
\[\bigl(\bar{\varphi}(\bar{x}_{r+1,1}),\ldots,\bar{\varphi}(\bar{x}_{p1})\bigr) = \bigl(\bar{y}_{r'+1,1},\ldots,\bar{y}_{p'1}\bigr) \cdot D \,.\]
For each $i \in [r+1,p]$, put
\begin{equation}\label{E:local1}
g_i \defeq \max \defset[1]{e'_{i'1}}{i' \in [r'+1,p'] \mbox{ such that } d_{i'-r',i-r} \neq 0}
\end{equation}
and let $\sigma(i)$ be the smallest number $i'$ in $[r'+1, p']$ such that $d_{i'-r',i-r} \neq 0$ and $e'_{i'1} = g_i$. Then for every $i \in [r+1,p]$, we may write $\varphi(x_i)$ as
\[\varphi(x_{i1}) = v_iy_{\sigma(i),1} + w_i \,,\]
where if we write $v_i$ and $w_i$ as the form \eqref{E:local0}, then the constant term of $v_i$ is \[d_{\sigma(i)-r',i-r} \in K^{\ast} \qquad \mbox{(so $v_i \in B^{\ast}$)}\]
and $w_i$ does not contain $y_{\sigma(i),1}$. For any $h \in [0,e'_{\sigma(i),1}-1]$, by considering the coefficient of $y_{\sigma(i),1}^h$ in $(v_iy_{\sigma(i),1} + w_i)^h$, we know that $\varphi(x_{i1})^h \neq 0$. Hence $e_{i1} \geqslant e'_{\sigma(i),1}$.

Suppose that $w_i \neq 0$. We write $w_i$ as
\[w_i = c_{i1}L_{i1} + c_{i2}L_{i2} + \cdots + c_{il_i}L_{il_i} + H_i \,,\]
where $L_{i1}, L_{i2}, \ldots, L_{il_i}$ are monic monomials occurred in $w_i$ with lowest degree $n$ ($\geqslant 1$), $c_{i1}, c_{i2}, \ldots, c_{il_i}$ are nonzero elements in $K$, and $H_i$ is the sum of monomials of degree greater than $n$ in $w_i$. The coefficient of $y_{\sigma(i),1}^{e'_{\sigma(i),1}-1}L_{i1}$ in $(v_iy_{\sigma(i),1} + w_i)^{e'_{\sigma(i),1}}$ is equal to
\[e'_{\sigma(i),1} \cdot d_{\sigma(i)-r',i-r} \cdot c_{i1} \neq 0 \,.\]
So $\varphi(x_{i1})^{e'_{\sigma(i),1}} \neq 0$. Hence $e_{i1} > e'_{\sigma(i),1}$.

As $D$ is an invertible matrix over $K$, there exists a $\sigma' \in S_f$ such that
\[d_{\sigma'(1),1},d_{\sigma'(2),2},\ldots,d_{\sigma'(f),r}\]
are all nonzero. By \eqref{E:local1}, for every $i \in [r+1,p]$ we have
\[e'_{\sigma'(i-r)+r',1} \leqslant g_i = e'_{\sigma(i),1} \leqslant e_{i1} \,.\]
Thus
\[\sum^{p'}_{i'=r'+1}e'_{i'1} = \sum^p_{i=r+1}e'_{\sigma'(i-r)+r',1} \leqslant \sum^p_{i=r+1}e_{i1} \,.\]
Applying above analysis to the homomorphism $\varphi^{-1} \colon B \to A$, we have \[\sum^p_{i=r+1}e_{i1} \leqslant \sum^{p'}_{i'=r'+1}e'_{i'1} \,.\]
Hence for each $i \in [r+1,p]$, we have
\[e'_{\sigma'(i-r)+r',1} = e'_{\sigma(i),1} = e_{i1}\]
and $w_i = 0$; and $\sigma(i) = \sigma'(i-r)+r'$ is a bijective from $[r+1,p]$ to $[r'+1,p']$.

In the following we prove that $p = p'$ and extend $\sigma$ to an element in $S_p$. Put
\begin{align*}
J & \defeq [1,q_1] \times [1,q_2] \times \cdots \times [1,q_p] \,, \\
J' & \defeq [1,q'_1] \times [1,q'_2] \times \cdots \times [1,q'_{p'}] \,.
\end{align*}
For each $j_{\cdot} = (\Sequ{p}{j}) \in J$, put
\begin{align*}
\mathfrak{a}_{j_{\cdot}} &= \mathfrak{a}_{j_1,j_2,\ldots,j_p} \defeq \bigl(x_{1j_1}^{e_{1j_1}}, x_{2j_2}^{e_{2j_2}}, \ldots, x_{pj_p}^{e_{pj_p}}\bigr) \,, \\
\mathfrak{p}_{j_{\cdot}} &= \mathfrak{p}_{j_1,j_2,\ldots,j_p} \defeq \bigl(x_{1j_1}, x_{2j_2}, \ldots, x_{pj_p}\bigr) \,;
\end{align*}
and for each $j'_{\cdot} = (j'_1,j'_2,\ldots,j'_{p'}) \in J'$, put
\begin{align*}
\mathfrak{b}_{j'_{\cdot}} &= \mathfrak{b}_{j'_1,j'_2,\ldots,j'_{p'}} \defeq \bigl(y_{1j'_1}^{e'_{1j'_1}}, y_{2j'_2}^{e'_{2j'_2}}, \ldots, y_{p'j'_{p'}}^{e'_{p'j'_{p'}}}\bigr) \,, \\
\mathfrak{q}_{j'_{\cdot}} &= \mathfrak{q}_{j'_1,j'_2,\ldots,j'_{p'}} \defeq \bigl(y_{1j'_1}, y_{2j'_2}, \ldots, y_{p'j'_{p'}}\bigr) \,.
\end{align*}
Then $\mathfrak{p}_{j_{\cdot}} = \sqrt{\mathfrak{a}_{j_{\cdot}}}$ is a prime ideal of $A$ and $\mathfrak{q}_{j'_{\cdot}} = \sqrt{\mathfrak{b}_{j'_{\cdot}}}$ is a prime ideal of $B$. Moreover
\[(0) = \bigcap_{j_{\cdot} \in J}\mathfrak{a}_{j_{\cdot}} \qquad \mbox{and} \qquad (0) = \bigcap_{j'_{\cdot} \in J'}\mathfrak{b}_{j'_{\cdot}}\]
are the primary decompositions of $(0) \subseteq A$ and $(0) \subseteq B$ respectively. Note that for every $j_{\cdot} \in J$,
\begin{equation}\label{E:local2}
\dim(A/\mathfrak{p}_{j_{\cdot}}) = \sum^p_{i=1}q_i - p + m \,,
\end{equation}
which does not depend on $j_{\cdot}$. So all $\mathfrak{p}_{j_{\cdot}}$ are isolated prime ideals belonging to $(0)$. Similarly we have
\begin{equation}\label{E:local3}
\dim(B/\mathfrak{q}_{j'_{\cdot}}) = \sum^{p'}_{i'=1}q'_{i'} - p' + m' \,,
\end{equation}
and all $\mathfrak{q}_{j'_{\cdot}}$ are isolated prime ideals belonging to $(0)$. By the uniqueness of primary decomposition of ideals, there is a bijective $\alpha \colon J \to J'$ such that for every $j_{\cdot} \in J$, $\varphi(\mathfrak{a}_{j_{\cdot}}) = \mathfrak{b}_{\alpha(j_{\cdot})}$ and $\varphi(\mathfrak{p}_{j_{\cdot}}) = \mathfrak{q}_{\alpha(j_{\cdot})}$.
By \eqref{E:local2} and \eqref{E:local3}, we have
\begin{equation}\label{E:local5}
\sum^p_{i=1}q_i - p + m = \dim(A/\mathfrak{p}_{j_{\cdot}}) = \dim(B/\mathfrak{q}_{\alpha(j_{\cdot})}) = \sum^{p'}_{i'=1}q'_{i'} - p' + m' \,.
\end{equation}
Note that $\varphi$ induces an isomorphism of rings:
\[A \Big/\hspace{-.5em} \sumT_{j_{\cdot} \in J} \mathfrak{p}_{j_{\cdot}} \xrightarrow{\sim} B \Big/\hspace{-.5em} \sumT_{j'_{\cdot} \in J'} \mathfrak{q}_{j'_{\cdot}} \,.\]
By comparing the dimensions of both sides, we get $m = m'$. For any $j_{\cdot} = (\Sequ{p}{j})$, $l_{\cdot} = (\Sequ{p}{l}) \in J$, put
\[\mathfrak{d}(j_{\cdot},l_{\cdot}) \defeq \# \defset[1]{i \in [1,p]}{j_i \neq l_i} \,;\]
and for each $j'_{\cdot} = (j'_1,j'_2,\ldots,j'_{p'})$, $l'_{\cdot} = (l'_1,l'_2,\ldots,l'_{p'}) \in J'$, put
\[\mathfrak{d}'(j'_{\cdot},l'_{\cdot}) \defeq \# \defset[1]{i' \in [1,p']}{j'_{i'} \neq l'_{i'}} \,.\]
For each $j_{\cdot},l_{\cdot} \in J$, we have
\begin{align*}
\sum^p_{i=1}q_i - p + m - \mathfrak{d}\bigl(j_{\cdot},l_{\cdot}\bigr)
&= \dim A/(\mathfrak{p}_{j_{\cdot}} + \mathfrak{p}_{l_{\cdot}}) \\
&= \dim B/(\mathfrak{q}_{\alpha(j_{\cdot})} + \mathfrak{q}_{\alpha(l_{\cdot})}) \\
&= \sum^{p'}_{i'=1}q'_{i'} - p' + m -\mathfrak{d}'\bigl(\alpha(j_{\cdot}),\alpha(l_{\cdot})\bigr) \,.
\end{align*}
By \eqref{E:local5}, we have
\begin{equation}\label{E:local6}
\mathfrak{d}\bigl(j_{\cdot},l_{\cdot}\bigr) = \mathfrak{d}'\bigl(\alpha(j_{\cdot}),\alpha(l_{\cdot})\bigr) \,.
\end{equation}
So we get
\begin{align*}
r &= \mathfrak{d}\bigl((1,\ldots,1), (2,\ldots,2,1,\ldots,1)\bigr) \\
&= \mathfrak{d}'\bigl(\alpha(1,\ldots,1), \alpha(2,\ldots,2,1,\ldots,1)\bigr) \\
& \leqslant r' \,.
\end{align*}
Applying above argument to $\alpha^{-1}$, we get $r' \leqslant r$ and hence $r = r'$. By \eqref{E:local9} we obtain $p = p'$.

For each $j_{\cdot} \in J$, we put
\[\alpha(j_{\cdot}) = \bigl(\alpha_1\bigl(j_{\cdot}), \alpha_2(j_{\cdot}), \ldots, \alpha_p(j_{\cdot})\bigr) \,.\]
For each $h \in [1,r]$, let $s$ and $t$ be two different numbers in $[1,q_h]$, and
\[(j_1,\ldots,\Hat{\jmath}_h,\ldots,j_p) \in [1,q_1] \times \cdots \times [1,q_{h-1}] \times [1,q_{h+1}] \times \cdots \times [1,q_p] \,.\]
By \eqref{E:local6}, there is a unique integer
\[\sigma = \sigma(h,s,t; j_1,\ldots,\Hat{\jmath}_h,\ldots,j_p) \in [1,r]\]
such that
\[\alpha_{\sigma}(j_1,\ldots,j_{h-1},s,j_{h+1},\ldots,j_p) \neq \alpha_{\sigma}(j_1,\ldots,j_{h-1},t,j_{h+1},\ldots,j_p) \,,\]
and for all $l \in [1,p] - \{\sigma\}$,
\[\alpha_l(j_1,\ldots,j_{h-1},s,j_{h+1},\ldots,j_p) = \alpha_l(j_1,\ldots,j_{h-1},t,j_{h+1},\ldots,j_p) \,.\]

First we prove that $\sigma(h,s,t; j_1,\ldots,\Hat{\jmath}_h,\ldots,j_p)$ does not depend on $j_1,\ldots,\Hat{\jmath}_h,\ldots,j_p$. For simplicity we assume that $h = 1$ and $\gamma,\delta \in [1,q_2]$ are two different numbers. Put
\begin{align*}
n_1 & \defeq \sigma(1,s,t; \gamma,j_3,\ldots,j_p) \,, \\
n_2 & \defeq \sigma(1,s,t; \delta,j_3,\ldots,j_p) \,.
\end{align*}
Suppose that $n_1 \neq n_2$. Then we have
\begin{align*}
\alpha_{n_2}(s,\gamma,j_3,\ldots,j_p) &= \alpha_{n_2}(t,\gamma,j_3,\ldots,j_p) \,, \\
\alpha_{n_2}(s,\delta,j_3,\ldots,j_p) & \neq \alpha_{n_2}(t,\delta,j_3,\ldots,j_p) \,.
\end{align*}
So either
\begin{align}\label{E:local7}
\alpha_{n_2}(s,\gamma,j_3,\ldots,j_p) & \neq \alpha_{n_2}(s,\delta,j_3,\ldots,j_p) \,, \\
\intertext{or} \label{E:local8}
\alpha_{n_2}(t,\gamma,j_3,\ldots,j_p) & \neq \alpha_{n_2}(t,\delta,j_3,\ldots,j_p) \,.
\end{align}
Assume that \eqref{E:local7} is valid, then
\[\sigma(2,\gamma,\delta;s,j_3,\ldots,j_p) = n_2 \,.\]
Thus for all $l \in [1,p] -\{n_2\}$, we have
\[\alpha_l(s,\gamma,j_3,\ldots,j_p) = \alpha_l(s,\delta,j_3,\ldots,j_p) = \alpha_l(t,\delta,j_3,\ldots,j_p) \,,\]
i.e.,
\[\mathfrak{d}'\bigl(\alpha(s,\gamma,j_3,\ldots,j_p), \alpha(t,\delta,j_3,\ldots,j_p)\bigr) \leqslant 1 \,,\]
which contradicts to \eqref{E:local6}. Similarly the validity of \eqref{E:local8} leads to a contradiction. Hence $n_1 = n_2$. So $\sigma(h,s,t; j_1,\ldots,\Hat{\jmath}_h,\ldots,j_p)$ depends only on $h,s,t$; thus we may write it as $\sigma(h,s,t)$. Clearly $\sigma(h,s,t) = \sigma(h,t,s)$.

Second we prove that $\sigma(h,s,t)$ does not depend on $s$ and $t$. We also assume that $h = 1$ and let $s_1, s_2, s_3$ be three different numbers in $[1,q_h]$. Put $n_1 \defeq \sigma(1,s_1,s_2)$ and $n_2 \defeq \sigma(1,s_1,s_3)$. Suppose that $n_1 \neq n_2$. Since
\[\alpha_{n_1}(s_3,1,\ldots,1) = \alpha_{n_1}(s_1,1,\ldots,1) \neq \alpha_{n_1}(s_2,1,\ldots,1) \,,\]
we have $\sigma(1,s_2,s_3) = n_1$. So
\[\alpha_{n_2}(s_3,1,\ldots,1) = \alpha_{n_2}(s_2,1,\ldots,1) = \alpha_{n_2}(s_1,1,\ldots,1) \,,\]
which contradicts to the fact that $\sigma(1,s_1,s_3) = n_2$. Thus $\sigma(h,s,t)$ depends only on $h$, so we may write it as $\sigma(h)$.

We shall prove that $\sigma \colon [1,r] \to [1,r]$ is injective. Suppose that $\sigma(1) = \sigma(2) = n$. Let
\[1 \leqslant s_1 < s_2 \leqslant q_1, \ 1 \leqslant t_1 < t_2 \leqslant q_2, 1 \leqslant j_3 \leqslant q_3, \ldots, \ 1 \leqslant j_p \leqslant q_p\]
be integers. Then for any $l \neq n$,
\[\alpha_l(s_1,t_1,j_3,\ldots,j_p) = \alpha_l(s_2,t_1,j_3,\ldots,j_p) = \alpha_l(s_2,t_2,j_3,\ldots,j_p) \,.\]
Thus
\[\mathfrak{d}'\bigl(\alpha(s_1,t_1,j_3,\ldots,j_p), \alpha(s_2,t_2,j_3,\ldots,j_p)\bigr) \leqslant 1 \,,\]
which contradicts to \eqref{E:local6}.

Therefore we obtain an element $\sigma \in S_p$.

From above discussion we see that for any
\[(h;j_1,\ldots,\Hat{\jmath}_h,\ldots,j_p) \in [1,r] \times [1,q_1] \times \cdots \times [1,q_{h-1}] \times [1,q_{h+1}] \times \cdots \times [1,q_p] \,,\]
there exists an injective map
\[\tau_h(j_1,\ldots,\Hat{\jmath}_h,\ldots,j_p) \colon [1,q_h] \to [1,q'_{\sigma(h)}]\]
such that for all $s \in [1,q_h]$,
\[\alpha_{\sigma(h)}(j_1,\ldots,j_{h-1},s,j_{h+1},\ldots,j_h) = \tau_h(j_1,\ldots,\Hat{\jmath}_h,\ldots,j_p)(s) \,,\]
and for all $l \in [1,p] - \{\sigma(h)\}$, the value of $\alpha_l(j_1,\ldots,j_{h-1},s,j_{h+1},\ldots,j_h)$ does not depend on $s$. Now we prove that $\tau_h(j_1,\ldots,\Hat{\jmath}_h,\ldots,j_p)$ does not depend on $j_1,\ldots,\Hat{\jmath}_h,\ldots,j_p$. We assume that $h = 1$ and $t_1, t_2 \in [1, q_2]$ are two different numbers. Suppose that there is an $s \in [1,q_1]$ such that
\[n_1 \defeq \tau_1(t_1,j_3,\ldots,j_p)(s) \neq n_2 \defeq \tau_1(t_2,j_3,\ldots,j_p)(s) \,.\]
Then
\[\alpha_{\sigma(1)}(s,t_1,j_2,\ldots,j_p) = n_1 \neq n_2 = \alpha_{\sigma(1)}(s,t_2,j_2,\ldots,j_p) \,.\]
Thus $\sigma(2) = \sigma(1)$, which contradicts to the injectivity of $\sigma$. So we have a well-defined injection $\tau_h \colon [1,q_h] \to [1,q'_{\sigma(h)}]$. Hence $q_h \leqslant q'_{\sigma(h)}$. So we get
\[\sum^p_{i=1}q_i \leqslant \sum^p_{i=1}q'_{\sigma(i)} = \sum^p_{i'=1}q'_{i'} \,.\]
Applying above argument to $\varphi^{-1}$, we obtain $\sum\limits^p_{i'=1}q'_{i'} \leqslant \sum\limits^p_{i=1}q_i$. Hence for every $i \in [1,p]$, we have $q_i = q'_{\sigma(i)}$ and $\tau_h$ is a bijective.

Put
\[\tau_{r+1} = \cdots = \tau_p = \iDe \colon \{1\} \to \{1\} \,.\]
Then for every $(\Sequ{p}{j}) \in J$, we have
\[\alpha(\Sequ{p}{j}) = \bigl(\tau_{\sigma^{-1}(1)}(j_{\sigma^{-1}(1)}),\tau_{\sigma^{-1}(2)}(j_{\sigma^{-1}(2)}),\ldots,\tau_{\sigma^{-1}(p)}(j_{\sigma^{-1}(p)})\bigr) \,.\]
In other words,
\begin{equation}\label{E:local10}
\varphi(\mathfrak{a}_{j_1,j_2,\ldots,j_p}) = \mathfrak{b}_{\tau_{\sigma^{-1}(1)}(j_{\sigma^{-1}(1)}),\tau_{\sigma^{-1}(2)}(j_{\sigma^{-1}(2)}),\ldots,\tau_{\sigma^{-1}(p)}(j_{\sigma^{-1}(p)})}
\end{equation}
and
\begin{equation}\label{E:local11}
\varphi(\mathfrak{p}_{j_1,j_2,\ldots,j_p}) = \mathfrak{q}_{\tau_{\sigma^{-1}(1)}(j_{\sigma^{-1}(1)}),\tau_{\sigma^{-1}(2)}(j_{\sigma^{-1}(2)}),\ldots,\tau_{\sigma^{-1}(p)}(j_{\sigma^{-1}(p)})} \,.
\end{equation}

Let $h \in [1,r]$ and $s \in [1,q_h]$. Put
\[x \defeq x_{hs} \,, \qquad y \defeq y_{\sigma(h),\tau_h(s)} \,, \qquad e \defeq e_{hs} \,, \qquad e' \defeq e'_{\sigma(h),\tau_h(s)}\]
for shortness. Then we have
\[(x^e) = \bigcap_{\substack{j_1 \in [1,q_1], \ldots, j_{h-1} \in [1,q_{h-1}], \\ j_{h+1} \in [1,q_{h+1}], \ldots, j_p \in [1,q_p]}} \mathfrak{a}_{j_1,\ldots,j_{h-1},s,j_{h+1},\ldots,j_p} \ ,\]
By \eqref{E:local10}, we have $\varphi\bigl((x^e)\bigr) = \bigl(y^{e'}\bigr)$. So $\varphi(x^e) = uy^{e'}$ for some $u \in B^{\ast}$. Note that $x^e \in \mathfrak{M}_1^e - \mathfrak{M}_1^{e+1}$ and $uy^{e'} \in \mathfrak{M}_2^{e'} - \mathfrak{M}_2^{e'+1}$. So $e = e'$. On the other hand, we have
\[(x) + \mathfrak{N}_1 = \bigcap_{\substack{j_1 \in [1,q_1], \ldots, j_{h-1} \in [1,q_{h-1}], \\ j_{h+1} \in [1,q_{h+1}], \ldots, j_p \in [1,q_p]}} \mathfrak{p}_{j_1,\ldots,j_{h-1},s,j_{h+1},\ldots,j_p} \ ,\]
By \eqref{E:local11}, we have
\[\varphi\bigl((x)+\mathfrak{N}_1\bigr) = (y) + \mathfrak{N}_2 \,.\]
So $\varphi(x) = vy + w$ for some $v \in B$ and $w \in \mathfrak{N}_2$. We write $v$ and $w$ as the form \eqref{E:local0} and assume that $w$ does not contain $y$. Suppose that $v \in \mathfrak{M}_2$. Then $\varphi(x) \in \mathfrak{M}_2^2 + \mathfrak{N}_2$. So $x \in \mathfrak{M}_1^2 + \mathfrak{N}_1$, a contradiction. Thus $v \in B^{\ast}$, i.e., the constant term $c_0$ of $v$ is nonzero.

Suppose that $w \neq 0$. We write $w$ as
\[w = c_1L_1 + c_2L_2 + \cdots + c_sL_s + H \,,\]
where $L_1, L_2, \ldots, L_s$ are monic monomials occurred in $w$ with lowest degree $n$ ($\geqslant 1$), $c_1, c_2, \ldots, c_s$ are nonzero elements in $K$, and $H$ is the sum of monomials of degree greater that $n$ in $w$. Note that
\[uy^e = \varphi(x^e) = \varphi(x)^e = (vy+w)^e \,.\]
By Comparing the coefficients of $y^{e-1}L_1$ in the above equality, we get $0 = ec_0c_1$, a contradiction. So $w = 0$, i.e., $x = vy$.
\end{proof}

\begin{proof}[Proof of Theorem \thmref{667T1}]
By Lemma \thmref{667L3}, $p = p'$, $m = m'$, and for every $i \in [1,p]$ and $j \in [1,q_i]$,
\[\varphi(x_{ij}) = u_{ij}y_{ij} + w_{ij}\]
for some $u_{ij} \in B^{\ast}$ and $w_{ij} \in \mathfrak{m}B$. (Here to without loss of generality, we assume that $\sigma$, $\Sequ{p}{\tau}$ are identities.) We express $u_{ij}$ and $w_{ij}$ in the form of \eqref{E:local0} and assume that $w_{ij}$ does not contain $y_{ij}$. For every integer $h \geqslant 1$, put $\mathfrak{a}_{ih} \defeq \mathfrak{m}^h + (b_i)$. Assume that we have proved that $a_i \in \mathfrak{a}_{ih}$ and $w_{ij} \in \mathfrak{a}_{ih}B$. Then we have
\begin{align}
a_i &= \varphi \bigl(\prodT^{q_i}_{j=1}x_{ij}^{e_{ij}}\bigr) = \prod^{q_i}_{j=1}(u_{ij}y_{ij}+w_{ij})^{e_{ij}} \notag \\
\label{E:local12}
& \equiv \sum^{q_i}_{j=1}e_{ij}u'_{ij}y_{i1}^{e_{ij}} \cdots y_{i,j-1}^{e_{i,j-1}}y_{ij}^{e_{ij}-1}y_{i,j+1}^{e_{i,j+1}} \cdots y_{i,q_i}^{e_{i,q_i}}w_{ij} 
\pmod{\mathfrak{a}_{i,h+1}B} \,,
\end{align}
where $u'_{ij} \in B^{\ast}$. Now we apply Lemma \thmref{667L2} to the \mDash{R/\mathfrak{a}_{i,h+1}}{algebra} $B/\mathfrak{a}_{i,h+1}B$. By comparing the constant terms in \eqref{E:local12}, we have $a_i \in \mathfrak{a}_{i,h+1}$. Suppose that $w_{ij} \notin \mathfrak{a}_{i,h+1}B$. Then we have
\[w_{ij} \equiv \sum^{s_{ij}}_{l=1}c_{ijl}L_{ijl} + H_{ij} \pmod{\mathfrak{a}_{i,h+1}B} \,,\]
where $c_{ij1},c_{ij2},\ldots,c_{ij,s_{ij}} \in \mathfrak{a}_{ih} - \mathfrak{a}_{i,h+1}$, $L_{ij1},L_{ij2},\ldots,L_{ij,s_{ij}}$ are different monoic monomials in $w_i$ with lowest degree $t_{ij}$, and $H_{ij}$ are sums of monomials of degree greater $t_{ij}$ in $w_i$. By  comparing the coefficients of the term
\[y_{i1}^{e_{ij}} \cdots y_{i,j-1}^{e_{i,j-1}}y_{ij}^{e_{ij}-1}y_{i,j+1}^{e_{i,j+1}} \cdots y_{i,q_i}^{e_{i,q_i}}L_{ij1}\]
in \eqref{E:local12}, we get a contradiction. So we have
\[a_i \in \bigcap^{\infty}_{h=1}\bigl(\mathfrak{m}^h + (b_i)\bigr) = (b_i)\]
and
\[w_{ij} \in \bigcap^{\infty}_{h=1}\bigl(\mathfrak{m}^hB + Bb_i\bigr) = Bb_i \,.\]
The same reasoning for $\varphi^{-1}$ shows that $b_i \in (a_i)$. So $a_i = u_ib_i$ for some $u_i \in R^{\ast}$. Put $w_{ij} \defeq w'_{ij}b_i$ and
\[v_{ij} \defeq u_{ij} + y_{i1}^{e_{ij}} \cdots y_{i,j-1}^{e_{i,j-1}}y_{ij}^{e_{ij}-1}y_{i,j+1}^{e_{i,j+1}} \cdots y_{i,q_i}^{e_{i,q_i}}w'_{ij} \,.\]
Then $v_{ij} \in B^{\ast}$ and $\varphi(x_{ij}) = v_{ij}y_{ij}$. This complete the proof of Theorem \thmref{667T1}.
\end{proof}

The following Theorem is easy to prove.

\begin{noteTh}[667P0]
For each $i \in [1,p]$, let $\mathfrak{a}_i$ denote the kernel of multiplication by $a_i$ on $R$; and for each $j \in [1,q_i]$, let $\mathfrak{I}_{ij}$ denote the kernel of multiplication by $x_{ij}$ on $A$. Then
\begin{enumerate}
\item for each $i \in [1,p]$ and $j \in [1,q_i]$, we have
\[\mathfrak{I}_{ij} = \mathfrak{a}_i \cdot \bigl(x_{i1}^{e_{i1}} \cdots x_{i,j-1}^{e_{i,j-1}} x_{ij}^{e_{ij}-1} x_{i,j+1}^{e_{i,j+1}} \cdots x_{iq_i}^{e_{iq_i}}\bigr) \,;\]
\item for each $i \in [1,p]$, the canonical homomorphism of \mDash{A}{modules}
\[\bigoplus^{q_i}_{j=1} \mathfrak{I}_{ij} \to \sum^{q_i}_{j=1} \mathfrak{I}_{ij}\]
is an isomorphism.
\end{enumerate}
\end{noteTh}

\section{\RlCs}
In this section we define the concept of \rlC, which is more delicate than \lCs. Log structures induced by \rlCs are all locally isomorphic. But it is not true for \lCs. Also we introduce two assumptions on which the main results of this paper is built.

Let $f \colon X \to Y$ be a surjective, proper and \wnC morphism of locally noetherian schemes.

\subsection{The singular locus}\

\begin{noteLe}[668P1]
Let $\varphi \colon A \to B$ be a flat local homomorphism of noetherian local rings, $x$ and $y$ two nonzero elements in $A$. If there is a $v \in B^{\ast}$ such that $\varphi(y) = v \varphi(x)$, then there exists a $u \in A^{\ast}$ such that $y = ux$.
\end{noteLe}

\begin{proof}
Since $A$ and $B$ are local rings and $\varphi$ is flat, we see that $\varphi$ is faithfully flat. So $xA = \varphi^{-1}(xB)$ and $yA = \varphi^{-1}(yB)$. Thus $xA = yA$ if and only if $xB = yB$.
\end{proof}

\begin{noteLe}[668P2]
Let $x \in X$, $y \defeq f(x)$,
\[U \to \Spec P \Big/ \Bigl(\prodT^{q_1}_{j_1=1}T_{1j_1}^{e_{1j_1}} - a_1, \prodT^{q_2}_{j_2=1}T_{2j_2}^{e_{2j_2}} - a_2, \ldots, \prodT^{q_p}_{j_p=1}T_{pj_p}^{e_{pj_p}} - a_p\Bigr)\]
and
\[U' \to \Spec P' \Big/ \Bigl(\prodT^{q'_1}_{j'_1=1}(T'_{1j'_1})^{e'_{1j'_1}} - a'_1, \prodT^{q'_2}_{j'_2=1}(T'_{2j'_2})^{e'_{2j'_2}} - a'_2, \ldots, \prodT^{q'_{p'}}_{j'_{p'}=1}(T'_{p'j'_{p'}})^{e'_{p'j'_{p'}}} - a'_{p'}\Bigr)\]
be two \lCs of $f$ at $x$. Let $t_{ij_i}$ and $t'_{i'j'_{i'}}$ be the image of $T_{ij_i}$ and $T'_{i'j'_{i'}}$ in $\mathcal{O}_{X,\bar{x}}$ respectively. Then $p = p'$ and there exists a $\sigma \in S_p$ such that for each $i \in [1,p]$, we have
\begin{enumerate}
\item $q_i = q'_{\sigma(i)}$,
\item $a_i = u_ia'_{\sigma(i)}$ for some $u_i \in \mathcal{O}_{Y,\bar{y}}^{\ast}$,
\item there exists a $\tau_i \in S_{q_i}$ such that for each $j \in [1,q_i]$, $e_{ij} = e'_{\sigma(i),\tau_i(j)}$ and $t_{ij} = v_{ij}t'_{\sigma(i),\tau_i(j)}$ for some $v_{ij} \in \mathcal{O}_{X,\bar{x}}^{\ast}$.
\end{enumerate}
\end{noteLe}

\begin{proof}
We use notations in Definition \thmref{667E1}. Let $x' \in U$ and $x'' \in U'$ be the points as in Definition \thmref{667E1} (4) and $y' \in V$ the point as in Definition \thmref{667E1} (2). Put $U'' \defeq U \times_X U'$. Then there is a point $x_0 \in U''$ which maps onto both $x'$ and $x''$. Let $x'_1$ be a closed point in $\overline{\{x_0\}} \subseteq U''$ and let $x_1$ be the image of $x'_1$ on $X$. Then $x_1$ is a closed point in $\overline{\{x\}}$. So by considering the cospecialization map $\mathcal{O}_{X,\bar{x}_1} \to \mathcal{O}_{X,\bar{x}}$, we may assume that $x = x_1$ is a closed point. Then $\kappa(x)/\kappa(y)$ is a finite extension of fields. By \cite[Ch.~0, (10.3.1)]{AGro1}, there is a complete noetherian local ring $R'$ whose residue field is algebraically closed and a flat local homomorphism $\mathcal{O}_{S,\bar{y}} \to R'$. By taking base extension $\Spec R' \to S$ and applying Lemma \thmref{668P1}, we may assume that $\kappa(y)$ is algebraic closed. As $x$ is a closed point, $\kappa(x) = \kappa(y)$ is algebraic closed. Thus $\kappa(x) = \kappa(x')$. Let $\mathfrak{M}$ and $\mathfrak{m}$ be the maximal ideals of $\mathcal{O}_{\Spec P, x'}$ and $\mathcal{O}_{V,y'}$ respectively. There there is a canonical isomorphism:
\[L \defeq \mathfrak{M}/(\mathfrak{M}^2 + \mathfrak{m}P) \xrightarrow{\sim} \varOmega_{P/R} \otimes_P P/\mathfrak{M} \,.\]
As $\ldots, \overline{T}_{ij_i},\ldots \in L$ are linearly independent over $\kappa(x')$, we may select $\Sequ{n}{T} \in \mathfrak{M}$ such that $\{\ldots, \overline{T}_{ij_i}, \ldots, \overline{T}_k, \ldots\}$ is a basis of $L$. By taking a connected affine open neighborhood of $x'$ in $\Spec P$, we may assume that $\bigl\{\ldots, d(T_{ij_i}), \ldots, d(T_k), \ldots\bigr\}$ is a basis of $\varOmega_{P/R}$. Then $\{\ldots, T_{ij_i}, \ldots, T_k, \ldots\}$ is algebraically independent over $R$, and $P$ is \'{e}tale over $R[\ldots, T_{ij_i}, \ldots, T_k, \ldots]$. So we have an isomorphism of \mDash{\widehat{\mathcal{O}}_{S,\bar{y}}}{algebras}:
\[\widehat{\mathcal{O}}_{X,\bar{x}} \xrightarrow{\sim} \widehat{\mathcal{O}}_{S,\bar{y}}[[\ldots, T_{ij_i}, \ldots, T_k, \ldots]] \Big/ \Bigl(\ldots, \prodT^{q_i}_{j_i=1}T_{ij_i}^{e_{ij_i}} - a_i,\ldots \Bigr) \,.\]
Similarly we have
\[\widehat{\mathcal{O}}_{X,\bar{x}} \xrightarrow{\sim} \widehat{\mathcal{O}}_{S,\bar{y}}[[\ldots, T'_{i'j'_{i'}}, \ldots, T'_{k'}, \ldots]] \Big/ \Bigl(\ldots, \prodT^{q'_{i'}}_{j'_{i'}=1}(T'_{i'j'_{i'}})^{e'_{i'j'_{i'}}} - a'_{i'},\ldots \Bigr) \,.\]
So the lemma is valid by Theorem \thmref{667T1}.
\end{proof}

Let $x$ be a point on $X$ equipped with a \lC of the form \eqref{E:defn0e}. For each $i \in [1,p]$, let $\mathscr{I}_i$ be the ideal of $\mathcal{O}_U$ generated by
\begin{equation}\label{E:global0e}
\defset[1]{T_{i1}^{e_{i1}} \cdots T_{i,j-1}^{e_{i,j-1}} T_{ij}^{e_{ij}-1} T_{i,j+1}^{e_{i,j+1}} \cdots T_{i,q_i}^{e_{i,q_i}}}{j \in [1,q_i]} \,.
\end{equation}
By Lemma \thmref{668P2}, we see that
\[\bigl(\Sequ{p}{\mathscr{I}}\bigr)\]
are independent of the choice of local charts (up to a unique permutation of the subscripts in $S_p$). From \eqref{E:defn1e}, we know
\[D_i(U/V) = \Spec \bigl(\mathcal{O}_U/\mathscr{I}_i\bigr)\]
for each $i \in [1,p]$. So for $U/V$ we have a finite \mDash{V}{morphism}
\[D_{U/V} \defeq \coprod^p_{i=1} D_i(U/V) \to U \,.\]
Clearly
\[\defset[1]{D_{U/V} \to U}{\mbox{$U/V$ is a local chart for $f$}}\]
can be glued to a global finite \mDash{S}{morphism} $g \colon D \to X$. To consider the properties under base extension, we also use $D(f)$ or $D(X/S)$ to denote the scheme $D$ for preciseness.

Obviously the set-theoretic image of the finite morphism $D(f) \to X$ is the set of all points at which $f$ are not smooth.

Clearly we have

\begin{noteTh}[668P12]
Let
\[\xymatrix{X' \ar[r]^{f'} \ar[d] \ar@{}[dr]|-*+{\square} & S' \ar[d] \\ X \ar[r]_{f} & S}\]
be a Cartesian square of locally noetherian schemes. Then we have
\[D(f') = D(f) \times_S S' \,.\]
\end{noteTh}

For a point $y \in Y$, we define
\[D_{\bar{y}} \defeq D\bigl((X \times_S \Spec \mathcal{O}_{Y,\bar{y}}) / \Spec \mathcal{O}_{Y,\bar{y}} \bigr) = D(X/S) \times_S \Spec \mathcal{O}_{Y,\bar{y}} \,,\]
and let $\cpC(y)$ denote the set of connected components of $D_{\bar{y}}$.

\subsection{Reduced to local cases}\label{Se:global:2}

In this subsection, we assume that $S = \Spec R$, where $R$ is a strictly Henselian noetherian local ring, and $y_1$ is the closed point of $S$.

\begin{noteLe}[668P19]
Let $T$ be the spectrum of a Henselian local ring, $t$ the closed point of $T$, $Y$ a connected scheme, $g \colon Y \to T$ a \'{e}tale morphism, $y$ a point on $Y$ such that $g(y) = t$ and $\kappa(t) \to \kappa(y)$ is isomorphic. Then $g \colon Y \to T$ is an isomorphism.
\end{noteLe}

\begin{proof}
See \cite{AGro2}, (18.5.11) a) $\Longrightarrow$ c) and (18.5.18).
\end{proof}

\begin{noteLe}[668P4]
Let $T$ be the spectrum of a Henselian local ring, $t$ the closed point of $T$, $Y \to T$ a proper morphism. Then $Z \mapsto Z_t$ defines a bijection from the set of connected components of $Y$ to the set of connected components of $Y_t$.
\end{noteLe}

\begin{proof}
See \cite[(18.5.19)]{AGro2}.
\end{proof}

By Lemma \thmref{668P19}, if $x \in f^{-1}(y_1)$ and $U/V$ is a \lC of $f$ at $x$, then $V = S$. By Lemma \thmref{668P2}, there is a canonical map
\begin{equation}\label{E:global1e}
\omega \colon \cpC(y_1) \to R/R^{\ast}
\end{equation}
(here ``$/$'' means taking quotient of monoids) such that if $x \in f^{-1}(y_1)$, $U$ a \lC of the form \eqref{E:defn0e} at $x$, $i \in [1,p]$, and $C \in \cpC(y_1)$ is the connected component of $D(X/S)$ which contains the image of $D_i(U/S)$ on $D(X/S)$, then
\[\overline{a_i} = \omega(C) \,.\]

Now we consider the following conditions.
\begin{quotation}
($\ast$) For each point $x \in f^{-1}(y_1)$ and for each \lC $U$ at $x$, the images of $D_1(U/S), D_2(U/S), \ldots, D_p(U/S)$ in $D(X/S)$ are contained in different connected components of $D(X/S)$.
\end{quotation}

\begin{noteLe}[668P14]
If \rmB{$\ast$} holds, then every connected components of $D(X/S)$ is a closed subscheme of $X$.
\end{noteLe}

\subsection{Construction of \rlC}\label{Se:global:3}

We return to the general case that $S$ is only a locally noetherian scheme. For each point $y \in S$, we use
\[\omega_{\bar{y}} \colon \cpC(y) \to \mathcal{O}_{S,\bar{y}}/\mathcal{O}_{S,\bar{y}}^{\ast}\]
to denote the map as in \eqref{E:global1e}.

In the following, we require that $f \colon X \to S$ satisfies the following conditions.
\begin{quotation}
(\dag) For each point $y \in S$, $X \times_S \Spec \mathcal{O}_{S,\bar{y}} \to \Spec \mathcal{O}_{S,\bar{y}}$ satisfies the condition ($\ast$).
\end{quotation}

\begin{noteLe}
Let $y$ be a point on $S$. Fix an open affine neighborhood $W_y$ of $y$. We define a full subcategory $N(y)$ of the category of \'{e}tale neighborhoods of $\bar{y}$ as follows: for an \'{e}tale neighborhood $V$ of $\bar{y}$, $V \in N(y)$ if and only if it satisfies the following conditions:
\begin{enumerate}
\renewcommand{\theenumi}{\alph{enumi}}
\item the image of $V$ in $S$ is contained in $W_y$;
\item the inverse image of $y$ in $V$ contains only one point $y' \in V$;
\item $V$ is an affine scheme and every irreducible component of $V$ contains $y'$;
\item for every irreducible component $F$ of $D_V$, the image of $F$ on $V$ contains $y'$.
\end{enumerate}
Then we have
\begin{enumerate}
\item For any pair of objects $V$ and $V'$ in $N(y)$, there exists at most one morphism from $V'$ to $V$;
\item $N(y)$ is a \emph{local base} of $\bar{y}$, i.e., for every \'{e}tale neighborhood $V$ of $\bar{y}$, there is an object $V'$ in $N(y)$ and a morphism $V' \to V$ of \'{e}tale neighborhoods of $\bar{y}$.
\item for any morphism $V' \to V$ in $N(y)$, $D_{V'} \to D_V$ is dominant.
\end{enumerate}
\end{noteLe}

\begin{proof}
(1) is by \cite[I 5.4]{AGro3}.

(2) For every \'{e}tale neighborhood $V$ of $\bar{y}$, we may contract $V$ under the Zariski topology to obtain an object in $N(y)$.

(3) For each object $V$ in $N(y)$, since the image of $\Spec \mathcal{O}_{S,\bar{y}} \to V$ contains all generizations of $y'$ in $V$, the morphism $D_{\bar{y}} \to D_V$ is dominant. So for any morphism $V' \to V$ in $N(y)$, $D_{V'} \to D_V$ is dominant.
\end{proof}

\begin{noteRe}
By (1), we may define a partial order on $N(y)$ as follows. For any pair of objects $V$ and $V'$ in $N(y)$, $V' \geqslant V$ if and only if there is a morphism from $V'$ to $V$ in $N(y)$. Obviously $N(y)$ is directly ordered, i.e., for any pair of objects $V$ and $V'$ in $N(y)$, there exists an object $V''$ in $N(y)$ such that $V'' \geqslant V$ and $V'' \geqslant V'$.
\end{noteRe}

Let $y$ be a point on $S$. For each object $V$ in $N(y)$, let $q_V \colon D_{\bar{y}} \to D_V$ be the canonical morphism. By \cite[(8.4.2)]{AGro2}, there exists an object $V_0$ in $N(y)$ such that for all $V \geqslant V_0$, $F \mapsto \overline{q_V(F)}$ defines a bijection between $\cpC(y)$ and the set of connected components of $D_V$. By the definition of $N(y)$, the inverse image of $y$ in $V_0$ contains only one point $y'$. For every closed point $x$ in $f_{V_0}^{-1}(y')$, if $f_{V_0}$ is smooth at $x$, we select open neighborhoods $U_x$ of $x$ and $V_x$ of $y$ respectively such that $f_{V_0}(U_x) \subseteq V_x$ and $f_{V_0} \colon U_x \to V_x$ is smooth; otherwise we select a \lC $U_x/V_x$ at $x$. For every $U_x$, let $U'_x$ be its image in $X_{V_0}$. As
\[f_{V_0}^{-1}(y') \subseteq \bigcup U'_x\]
and $f_{V_0}^{-1}(y')$ is quasi-compact, there exists a finite number of closed points $\Sequ{m}{x}$ in $f_{V_0}^{-1}(y')$ such that
\[f_{V_0}^{-1}(y') \subseteq \bigcup^m_{i=1}U'_{x_i} \,.\]
Put
\[V' \defeq V_0 - f_{V_0}\left(X_{V_0} - \cupT^m_{i=1}U'_{x_i} \right) \,.\]
As $f_{V_0} \colon X_{V_0} \to V_0$ is proper, $V'$ is an open neighborhood of $y'$. There exists an object $V_1$ in $N(y)$ and a morphism of \'{e}tale neighborhoods of $\bar{y}'$,
\[V_1 \to V' \times_{V_0} V_{x_1} \times_{V_0} \cdots \times_{V_0} V_{x_m} \,.\]

Let $\Sequ{n}{C}$ be all connected components of $D_{\bar{y}}$. For each $i \in [1,n]$ and each $V \in N(y)$ with $V \geqslant V_1$, $C_i$ defines a connected component $C_i(V)$ of $D_V$. By Lemma \thmref{668P2}, there exists an element
\[b_i(V) \in \varGamma(V, \mathcal{O}_S/\mathcal{O}_S^{\ast})\]
such that for every point $z \in V$ and for every connected component $F$ of $D_{\bar{z}}$ which maps into $C_i(V)$,
\[\omega_{\bar{z}}(F) = b_i(V)_{\bar{z}} \,.\]
Obviously $b_i(V)$ depends only on $y$, $C_i$ and $V$. Let $Z_i(V)$ be the closed subscheme of $V$ defined by the ideal generated by $b_i(V)$. Clearly the inverse image of $y'$ ($\in V_0$) in $V$ is contained in all these subschemes $Z_i(V)$.

The following lemma can be directly verified.

\begin{noteLe}
Let $V' \geqslant V$ \rmB{$\geqslant V_1$} be two elements in $N(y)$. Then $C_i(V') = C_i(V) \times_V V'$ and $Z_i(V') = Z_i(V) \times_V V'$.
\end{noteLe}

\begin{noteLe}[668P10]
$C_i(V) \to V$ factors through $Z_i(V)$ and $C_i(V)$ is faithfully flat over $Z_i(V)$.
\end{noteLe}

\begin{proof}
It is by the following lemma.
\end{proof}

\begin{noteLe}
Let $R$ be a ring, $A = R[\Sequ{n}{T}]$ a polynomial ring over $R$, $\Sequ{n}{e}$ positive integers. Then
\[B \defeq A/(T_1^{e_1-1}T_2^{e_2} \cdots T_n^{e_n}, T_1^{e_1}T_2^{e_2-1} \cdots T_n^{e_n}, \ldots, T_1^{e_1}T_2^{e_2} \cdots T_n^{e_n-1})\]
is flat over $R$.
\end{noteLe}

\begin{proof}
Note that $B$ is a free \mDash{R}{module} with basis
\begin{gather*}
\defset[2]{T_1^{i_1}T_2^{i_2} \cdots T_n^{i_n}}{\parbox{25em}{either there exists an integer $k \in [1,n]$ such that $i_k < e_k - 1$ or there exist at least two integers $k \in [1,n]$ such that $i_k \leqslant e_k -1$}} \,. \qedhere
\end{gather*}
\end{proof}

\begin{noteNo}[668P9]
We defines a subset $N_0(y)$ of $N(y)$ as follows: for $V \in N(y)$, $V \in N_0(y)$ if and only if it satisfies that
\begin{enumerate}
\item $V \geqslant V_1$.
\item For each $i \in [1, n]$, there exists a section $a_i \in \varGamma(V,\mathcal{O}_V)$ such that
\[a_i \equiv b_i(V) \pmod{\mathcal{O}_V^{\ast}} \,.\]
\item If $z$ is the inverse image of $y$ in $V$, then for any $i \in [1, n]$ and any irreducible component $F$ of $Z_i(V)$, $F$ contains $z$.
\end{enumerate}
\end{noteNo}

It is easy to show that for each $V \in N(y)$, there exists an object $V'$ in $N_0(y)$ such that $V' \geqslant V$; and we have the following easy lemma.

\begin{noteLe}[668P5]
Let $z$ be a generization of $y$, $u \colon \mathcal{O}_{S,\bar{y}} \to \mathcal{O}_{S,\bar{z}}$ a cospecialization map, and $v \colon D_{\bar{z}} \to D_{\bar{y}}$ the morphism induced by $u$. 
\begin{enumerate}
\item For each $i$, $v^{-1}(C_i) = \emptyset$ if and only if $\bar{u}\bigl(\omega_{\bar{y}}(C_i)\bigr) = \bar{1}$ in $\mathcal{O}_{S,\bar{z}}/\mathcal{O}_{S,\bar{z}}^{\ast}$.
\item If $v^{-1}(C_i) \neq \emptyset$, then every connected component of $v^{-1}(C_i)$ is a connected component of $D_{\bar{z}}$.
\item All connected components of $D_{\bar{z}}$ can be obtained as in {\upshape(2)}.
\end{enumerate}
\end{noteLe}

We consider the following conditions about $f$:
\begin{quotation}
(\ddag) Let $y_1$ be a point on $S$, $y_0$ a generization of $y_1$, $u \colon \mathcal{O}_{S,\bar{y}_1} \to \mathcal{O}_{S,\bar{y}_0}$ a cospecialization map, and $v \colon D_{\bar{y}_0} \to D_{\bar{y}_1}$ the morphism induced by $u$. Then for every connected component $C$ of $D_{\bar{y}_1}$, $v^{-1}(C)$ is connected. (Here empty set is also regarded to be connected.)
\end{quotation}

\begin{noteLe}[668P7]
\rmB{\ddag} is satisfied if one of the following conditions holds:
\begin{enumerate}
\item $S$ is a spectrum of a field.
\item There exists a finite set $L$ of closed points in $S$ such that $f$ is smooth outside $L$.
\item $S$ is a spectrum of a discrete valuative ring and $f$ is smooth at the generic fiber.
\item $f$ is a \wnC morphism without powers.
\end{enumerate}
\end{noteLe}

\begin{proof}
(1), (2) and $(3)$ are trivial.

(4) Obviously for each $y \in S$, $X_y$ is geometrically reduced over $\kappa(y)$. Let $C$ be a connected component of $D_{\bar{y}_1}$ such that $v^{-1}(C) \neq \emptyset$. Fix an object $V \in N_0(y_1)$. Let $y'_1$ be the inverse image of $y_1$ in $V$. Then the cospecialization map $u \colon \mathcal{O}_{S,\bar{y}_1} \to \mathcal{O}_{S,\bar{y}_0}$ defines a point $y'_0$ on $V$ which maps to $y_0$ and is a generization of $y'_1$. $C$ defines a connected component $\bar{C}$ of $D_V$. By Lemma \thmref{668P10}, there is a closed subscheme $Z$ of $V$ such that $\bar{C}$ factors through $Z$ and $\bar{C}$ is proper and flat over $Z$. As $y'_1 \in Z$ and $v^{-1}(C) \neq \emptyset$, by Lemma \thmref{668P5} $y'_0 \in Z$. Obviously $\bar{C}_{y'_1}$ is geometrically connected and geometrically reduced over $\kappa(y'_1)$. Hence
\[\dim_{\kappa(y'_1)} \varGamma \bigl(\bar{C}_{y'_1}, \mathcal{O}_{\bar{C}_{y'_1}}\bigr) = 1 \,.\]
By \cite[(7.7.5)]{AGro1},
\[\defset[1]{z \in V}{\dim_{\kappa(z)} \varGamma \bigl(\bar{C}_z, \mathcal{O}_{\bar{C}_z}\bigr) \leqslant 1}\]
is an open neighborhood of $y'_1$. As $y'_0$ is a generization of $y'_1$, 
\[\dim_{\kappa(y_0)} \varGamma \bigl(\bar{C}_{y_0}, \mathcal{O}_{\bar{C}_{y_0}}\bigr) \leqslant 1 \,.\]
 So $\bar{C}_{y_0}$ is geometrically connected. By Lemma \thmref{668P4}, $v^{-1}(C)$ is connected.
\end{proof}

\begin{noteLe}
Assume that $f$ satisfies the condition \rmB{\ddag}. Let $y \in S$, $V \in N_0(y)$, $\Sequ{n}{C}$ the connected components of $D_{\bar{y}}$. For each $i \in [1,n]$, let $\bar{C}_i$ be the connected component of $D_V$ defined by $C_i$ and let $a_i \in \varGamma(V,\mathcal{O}_S)$ be a representative element of $b_i(V)$. Let $z'$ be a point on $V$ and $z$ its image on $S$. Then
\begin{enumerate}
\item Let $i \in [1,n]$ such that $(a_i)_{z'} \in \mathfrak{m}_{V,z'}$. Then $\bar{C}_i \times_V \Spec \mathcal{O}_{S,\bar{z}}$ is a connected component of $D_{\bar{z}}$ and its image under the map $\omega_{\bar{z}}$ is equal to $\overline{(a_i)_{z'}}$.
\item $\defset[1]{\bar{C}_i \times_V \Spec \mathcal{O}_{S,\bar{z}}}{i \in [1,n] \mbox{ and } (a_i)_{z'} \in \mathfrak{m}_{V,z'}}$ is the set of all connected components of $D_{\bar{z}}$.
\end{enumerate}
\end{noteLe}

\begin{proof}
(1) Obviously $z' \in Z_i(V)$. Let $F$ be an irreducible component of $Z_i(V)$ containing $z'$. Let $w'$ be the generic point of $F$ and $w$ its image on $S$. Suppose that $\bar{C}_i \times_V \Spec \mathcal{O}_{S,\bar{z}}$ is disconnected. Then by Lemma \thmref{668P5}, $\bar{C}_i \times_V \Spec \mathcal{O}_{S,\bar{w}}$ is disconnected. By the definition of $N_0(y)$, we see that $y \in F$. So $w$ is a generization of $y$. Let
\[u \colon \mathcal{O}_{S,\bar{y}} = \mathcal{O}_{S,\bar{y}'} \to \mathcal{O}_{S,\bar{w}'} = \mathcal{O}_{S,\bar{w}}\]
be a cospecialization map and $v \colon D_{\bar{w}} \to D_{\bar{y}}$ the induced morphism. Then \[v^{-1}(C_i) = \bar{C}_i \times_V \Spec \mathcal{O}_{S,\bar{w}} \,.\]
Since $f$ satisfies the condition (\ddag), $v^{-1}(C_i)$ is connected, a contradition.

(2) is a consequence of (1).
\end{proof}

\begin{noteDe}
Let $x \in X$ be a point and $y \defeq f(x)$. A \emph{\rlC} of $f$ at $x$ is of the form
\[(U,V; T_{11},\ldots,T_{1q_1}; \ldots; T_{p1},\ldots,T_{pq_p}; a_1,\ldots,a_n)\]
where
\[U \to \Spec P \Bigl/ \Bigl(\ldots,\prodT^{q_i}_{j_i=1}T_{ij_i}^{e_{ij_i}} - a_i, \ldots \Bigr)\]
is a \lC at $x$, $V \in N_0(y)$, $n \defeq \# \cpC(y) \geqslant p$ and $a_i \equiv b_i(V) \pmod{\mathcal{O}_V^{\ast}}$ for all $i \in [1,n]$. We also simply use $U/V$ or $U$ to denote a \rlC.
\end{noteDe}

\begin{noteRe}
If $f$ is smooth at $x$, then $p = 0$; and if $f$ is smooth at the fiber $X_y$, then $p = n = 0$.
\end{noteRe}

\begin{noteRe}
Obviously every \lC can be contracted in the sense of \'{e}tale topology to become a \rlC.
\end{noteRe}

\begin{noteRe}[668P6]
Let $\mathscr{M}_U$ be the log structure on $U$ associating to $\alpha_U \colon \mathbb{N}^m_U \to \mathcal{O}_U$ with $m \defeq \sum\limits^p_{i=1}q_i + n - p$, where if
\[\eta_{11},\ldots,\eta_{1q_1}, \ldots, \eta_{p1},\ldots,\eta_{pq_p}, \eta_{p+1},\ldots,\eta_n\]
is a basis of $\mathbb{N}^m$, then $\alpha_U(\eta_{ij_i}) = \bar{T}_{ij_i}$ for $i \in [1,p]$ and $j_i \in [1,q_i]$, and $\alpha_U(\eta_i) = a_i$ for $i \in [p+1,n]$. Let $\mathscr{N}_V$ be the log structrue on $V$ associating to $\beta_V \colon \mathbb{N}^n_V \to \mathcal{O}_V$, where if $\varepsilon_1,\ldots,\varepsilon_n$ is a basis of $\mathbb{N}^n$, then $\beta_V(\varepsilon_i) = a_i$ for all $i \in [1,n]$. Let $g \colon U \to V$ be the canonical morphism. Then there is a canonical morphism \[\varphi_{U/V} \colon g^{\ast}\mathscr{N}_V \to \mathscr{M}_U\]
defined by the map $\gamma \colon \mathbb{N}^n \to \mathbb{N}^m$, where $\gamma(\varepsilon_i) = \sum\limits^{q_i}_{j=1}e_{ij}\eta_{ij}$ for $i \in [1,p]$ and $\gamma(\varepsilon_i) = \eta_i$ for $i \in [p+1,n]$.
\end{noteRe}

\begin{noteRe}[668P8]
As $\overline{\mathscr{M}}_U = \mathbb{N}_U^m/\alpha_U^{-1}(\mathcal{O}_U^{\ast})$ does not depend on the choice of $a_i$ and $T_{ij_i}$, we may glue the sheaves $\overline{\mathscr{M}}_U$ to obtain a global sheaf $\mathscr{P}$ of monoids on $\etSite{X}$ and there is a canonical morphism
\[\theta \colon \mathscr{P} \to \mathcal{O}_X/\mathcal{O}_X^{\ast} \,.\]
Similarly we may glue the sheaves $\overline{\mathscr{N}}_V$ to obtain a global sheaf $\mathscr{Q}$ of monoids on $\etSite{S}$ and there is a canonical morphism
\[\vartheta \colon \mathscr{Q} \to \mathcal{O}_S/\mathcal{O}_S^{\ast} \,.\]
Moreover, there is a canonical morphism $\mathfrak{d} \colon f^{-1}\mathscr{Q} \to \mathscr{P}$ defined by $\gamma$ which makes the following diagram commutative:
\[\xymatrix@C+2em{f^{-1}\mathscr{Q} \ar[r]^-{f^{-1}\vartheta} \ar[d]_{\mathfrak{d}} & f^{-1}\bigl(\mathcal{O}_S/\mathcal{O}_S^{\ast}\bigr) \ar[d] \\ \mathscr{P} \ar[r]_-{\theta} & \mathcal{O}_X/\mathcal{O}_X^{\ast}}\]
\end{noteRe}

\begin{noteLe}[668P16]
$\mathscr{Q}$ is canonically isomorphic to the direct image of $\mathbb{N}_{D(f)}$ under morphism $D(f) \to S$.
\end{noteLe}

\section{Cohomology and Hypercoverings}
\newcommand{\hCoh}[1]{\mathbb{H}^{#1}}

In this section, we review some technique in \cite{JGir1}. A brief version can be found in \cite[\S2 and \S3]{SScBSi}.

\subsection{Cohomology}\label{Se:hyper:1}

Let $X$ be a scheme. We define a category $\mathfrak{U}(X)$ as follows: an object in $\mathfrak{U}(X)$ is a diagram
\begin{equation}\label{E:hyper0e}
\vcenter{\xymatrix{V \ar@<.4ex>[r]^{v_1} \ar@<-.4ex>[r]_{v_2} & U \ar[r]^{u} & X}}
\end{equation}
where $U$ and $V$ are schemes, $u \colon U \to X$ and $\bMor{v_1,v_2 \colon V}{U}$ are surjective \'{e}tale morphisms such that $u \circ v_1 = u \circ v_2$ and the induced morphism
\[(v_1,v_2)_X \colon V \to U \times_X U\]
is surjective (and obviously \'{e}tale); we also simply use $U/V$ to denote the object \eqref{E:hyper0e}; a morphism in $\mathfrak{U}(X)$ is a pair of morphisms
\[(f,g) \colon U'/V' \to U/V\]
which makes a commutative diagram
\[\xymatrix{V' \ar@<.4ex>[r]^{v'_1} \ar@<-.4ex>[r]_{v'_2} \ar[d]_{g} & U' \ar[d]^{f} \ar[r]^{u'} & X \ar@{=}[d] \\ V \ar@<.4ex>[r]^{v_1} \ar@<-.4ex>[r]_{v_2} & U \ar[r]_{u} & X}\]

Given an object of form \eqref{E:hyper0e} in $\mathfrak{U}(X)$. Put $(V/U)_0 \defeq U$, $(V/U)_1 \defeq V$, $p_{00} \defeq u$, $p_{10} \defeq v_1$, $p_{11} \defeq v_2$. Assume that for some integer $n \geqslant 2$, we have schemes $(V/U)_k$ for $k \in [1,n-1]$ and \'{e}tale morphisms
\[p_{ki} \colon (V/U)_k \to (V/U)_{k-1}\]
for $i \in [0,k]$ such that whenever $0 \leqslant i < j \leqslant k$, we have
\begin{equation}\label{E:hyper1e}
p_{k-1,i} \circ p_{kj} = p_{k-1,j-1} \circ p_{ki} \,.
\end{equation}
Put
\[P_n \defeq \underbrace{(V/U)_{n-1} \times_X (V/U)_{n-1} \times_X \cdots \times_X (V/U)_{n-1}}_{n+1 \text{ copies of } (V/U)_{n-1}}\]
and let $q_{ni} \colon P_n \to (V/U)_{n-1}$ be the \mDash{(i+1)}{th} projection. For each $0 \leqslant i < j \leqslant n$, let $K(n,i,j)$ be the equalizer of $p_{n-1,i} \circ q_{nj}$ and $p_{n-1,j-1} \circ q_{ni}$ in the category of schemes. As $(V/U)_{n-1}$ is \'{e}tale over $X$, $K(n,i,j)$ is an open subscheme of $P_n$. Put
\[(V/U)_n \defeq \bigcap_{0 \leqslant i < j \leqslant n} K(n,i,j)\]
and
\[p_{ni} \defeq q_{ni}|_{(V/U)_n} \colon (V/U)_n \to (V/U)_{n-1} \,.\]

Let $\mathscr{F}$ be an abelian sheaf on $\etSite{X}$. We define a cochain complex of abelian groups as follows: for each $n \in \mathbb{N}$, put
\[C^n(V/U,\mathscr{F}) \defeq \varGamma\bigl((V/U)_n,\mathscr{F}\bigr)\]
and let
\[d^n \defeq \sum^{n+1}_{i=0}(-1)^ip_{n+1,i}^{\ast} \colon C^n(V/U,\mathscr{F}) \to C^{n+1}(V/U,\mathscr{F})\]
be the differential. Let $H^n(V/U,\mathscr{F})$ be the corresponding cohomology group. We define
\[\hCoh{n}(X,\mathscr{F}) \defeq \varinjlim H^n(V/U,\mathscr{F}) \,,\]
where the colimit runs through all elements in $\mathfrak{U}(X)$.

Let
\[0 \to \mathscr{F}' \xrightarrow{f} \mathscr{F} \xrightarrow{g} \mathscr{F}'' \to 0\]
be an exact sequence of abelian sheaves on $\etSite{X}$. For $i=0,1$, we define
\[\delta^i \colon \hCoh{i}(X,\mathscr{F}'') \to \hCoh{i+1}(X,\mathscr{F}')\]
as follows.

Let $s \in Z^0(V/U,\mathscr{F}'')$ be a \mDash{0}{cocycle}. By refining $U$, we may choose a lifting $\tilde{s} \in \mathscr{F}(U)$ of $s$. Then we define
\[\delta^0(s) = p_{10}^{\ast}(\tilde{s}) - p_{11}^{\ast}(\tilde{s}) \,.\]

Let $t \in Z^1(V/U,\mathscr{F}'')$ be a \mDash{1}{cocycle}. By refining $V$, we may choose a lifting $\tilde{t} \in \mathscr{F}(V)$ of $t$. Then we define
\[\delta^1(t) = p_{20}^{\ast}(\tilde{t}) - p_{21}^{\ast}(\tilde{t}) + p_{22}^{\ast}(\tilde{t}) \,.\]

\begin{noteTh}
For each $n=0,1,2$, we have a natural equivalence
\[\hCoh{n}(X,\!\smdot[.5m]) \xrightarrow{\sim} H^n(\etSite{X},\!\smdot[.5em]) \,;\]
and for each short exact sequence the natural equivalences commute with the connecting functors $\delta$.
\end{noteTh}

\subsection{Gerbe}

Now we fix an abelian sheaf $\mathscr{F}$ on $\etSite{X}$.

\begin{noteDe}[669P3]
An \mDash{\mathscr{F}}{gerbe} $(\mathfrak{X},\omega)$ consists of the following two data:
\begin{enumerate}
\renewcommand{\theenumi}{\alph{enumi}}
\item a stack $\mathfrak{X}$ over $\etSite{X}$;
\item for each \'{e}tale \mDash{X}{scheme} $U$ and for each object $A$ in $\mathfrak{X}(U)$, an isomorphism of sheaves of groups:
\[\omega(A) \colon \mathscr{F}|_U \xrightarrow{\sim} \shAut_U(A) \,.\]
\end{enumerate}
These data satisfy the following conditions: 
\begin{enumerate}
\item for any \'{e}tale \mDash{X}{scheme} $U$, there exists an \'{e}tale covering $\{U_i \to U\}_{i \in I}$ in $\etSite{X}$ such that $\mathfrak{X}(U_i) \neq \emptyset$ for all $i \in I$;
\item for any \'{e}tale \mDash{X}{scheme} $U$ and any pair of objects $A$ and $B$ in $\mathfrak{X}(U)$, there exists an \'{e}tale covering $\{U_i \to U\}_{i \in I}$ in $\etSite{X}$ such that $A|_{U_i}$ and $B|_{U_i}$ are isomorphic in $\mathfrak{X}(U_i)$ for all $i \in I$;
\item for any \'{e}tale \mDash{X}{scheme} $U$, any element $g \in \mathscr{F}(U)$, and any isomorphim $\varphi \colon A \xrightarrow{\sim} B$ in $\mathfrak{X}(U)$, we have
\[\varphi \circ \omega(A)(g) = \omega(B)(g) \circ \varphi \,.\]
(So we may simply write $\varphi \circ g$ or $g \circ \varphi$ or even $g \cdot \varphi$ for above morphism.)
\end{enumerate}
\end{noteDe}

Fix an \mDash{\mathscr{F}}{gerbe} $\mathfrak{X}$. Choose an \'{e}tale covering $U \to X$ which admits an object $A \in \mathfrak{X}(U)$, and an \'{e}tale covering $V \to U \times_X U$ which admits an isomorphism
\[\phi \colon p_{10}^{\ast}(A) \xrightarrow{\sim} p_{11}^{\ast}(A)\]
in $\mathfrak{X}(V)$. Then there exists a cocycle $g \in Z^2(V/U,\mathscr{F})$ such that
\[g \circ p_{21}^{\ast}(\phi) = p_{22}^{\ast}(\phi) \circ p_{20}^{\ast}(\phi) \colon (p_{10} \circ p_{20})^{\ast}(A) \xrightarrow{\sim} (p_{11} \circ p_{22})^{\ast}(A) \,.\]
We define
\[[\mathfrak{X}] \defeq [g] \in H^2(\etSite{X},\mathscr{F}) \,.\]

\begin{noteLe}[669P1]
Let $V/U$ be an object in $\mathfrak{A}(X)$, $A$ an object in $\mathfrak{X}(U)$, and $\phi \colon p_{10}^{\ast}(A) \xrightarrow{\sim} p_{11}^{\ast}(A)$ an isomorphism in $\mathfrak{X}(V)$ satisfying the cocycle condition:
\[p_{21}^{\ast}(\phi) = p_{22}^{\ast}(\phi) \circ p_{20}^{\ast}(\phi) \colon (p_{10} \circ p_{20})^{\ast}(A) \xrightarrow{\sim} (p_{11} \circ p_{22})^{\ast}(A) \,.\]
Then there exists an object $B$ in $\mathfrak{X}(X)$ and an isomorphism $\varphi \colon B|_U \xrightarrow{\sim} A$ in $\mathfrak{X}(U)$ such that $\phi \circ p_{10}^{\ast}(\varphi) = p_{11}^{\ast}(\varphi)$. And $(B,\varphi)$ is unique up to isomorphism.
\end{noteLe}

\begin{noteTh}[669P0]
$\mathfrak{X}(X) \neq \emptyset$ if and only if $[\mathfrak{X}] = 0$ in $H^2(\etSite{X},\mathscr{F})$.
\end{noteTh}

Let $g \in Z^1(V/U,\mathscr{F})$, $A$ an object in $\mathfrak{X}(X)$. By Lemma \thmref{669P1}, there is an object $g(A)$ in $\mathfrak{X}(X)$ and an isomorphism $\phi_g \colon g(A)|_U \to A|_U$ in $\mathfrak{X}(U)$ such that $g \circ p_{10}^{\ast}(\phi_g) = p_{11}^{\ast}(\phi_g)$. It is easy to show that $(g,A) \mapsto g(A)$ defines an action of the group $H^1(\etSite{X},\mathscr{F})$ on the set of isomorphic classes of objects in $\mathfrak{X}(X)$.

Let $A$ and $B$ be objects in $\mathfrak{X}(X)$, $\phi \colon A|_U \xrightarrow{\sim} B|_U$ an isomorphism in $\mathfrak{X}(U)$ and $g \in \mathscr{F}(V)$ such that $g \circ p_{10}^{\ast}(\phi) = p_{11}^{\ast}(\phi)$. Then $g \in Z^1(V/U,\mathscr{F})$ and $g(B) \cong A$ in $\mathfrak{X}(X)$.

\begin{noteTh}[669P2]
If $\mathfrak{X}(X) \neq \emptyset$, then the group $H^1(\etSite{X},\mathscr{F})$ acts canonically and simply transitively on the set of isomorphic classes of objects in $\mathfrak{X}(X)$.
\end{noteTh}

\section{Local Cases}
\label{Se:loc}

Let $f \colon X \to S$ be a surjective, proper and \wnC morphism of locally noetherian schemes which satisfies the conditions (\dag) and (\ddag) in \S\ref{Se:global:3}, $\Sequ{n}{D}$ the connected components of $D(f)$. Assume that there exist global sections
\[\Sequ{n}{a} \in \varGamma(S,\mathcal{O}_S)\]
such that for any point $y \in S$ and any $i \in [1,n]$, the following two conditions holds:
\begin{enumerate}
\item if $a_{i,\bar{y}} \in \mathcal{O}_{S,\bar{y}}^{\ast}$, then $D_i \times_S \Spec \mathcal{O}_{S,\bar{y}} = \emptyset$;
\item if $a_{i,\bar{y}} \in \mathfrak{m}_{S,\bar{y}}$, then
\[D_{i,\bar{y}} \defeq D_i \times_S \Spec \mathcal{O}_{S,\bar{y}}\]
is a connected component of $D_{\bar{y}}$ and $\omega_{\bar{y}}(D_{i,\bar{y}}) = \overline{a_{i,\bar{y}}}$.
\end{enumerate}
Thus
\[\cpC(y) = \defset[1]{D_{i,\bar{y}}}{i \in [1,n], \ a_{i,\bar{y}} \in \mathfrak{m}_{S,\bar{y}}} \,.\]
Obviously when $S$ is a spectrum of a strictly Henselian local ring, then above conditions hold. Furthermore, for any point $y \in S$ and any $V \in N_0(y)$, $X_V \to V$ satisfies above conditions.

Let $\mathscr{N}$ be the log structure on $S$ defined by
\[\mathbb{N}^n_S \to \mathcal{O}_S \,, \qquad \varepsilon_i \mapsto a_i \,,\]
where $\Sequ{n}{\varepsilon}$ is a basis of $\mathbb{N}^n$.

For each $i \in [1,n]$, let $\mathscr{I}_i$ be the ideal sheaf on $X$ corresponding to the closed subscheme $D_i$ and let $\mathscr{K}_i$ denote the kernel of the multiplication by $a_i$ on $\mathcal{O}_S$. As $f$ is flat, the kernel of the multiplication by $a_i$ on $\mathcal{O}_X$ is equal to $\mathscr{K}_i \cdot \mathcal{O}_X$. Let $E_i$ be the closed subscheme of $X$ defined by $\mathscr{K}_i \cdot \mathscr{I}_i$ and put $E \defeq \coprod\limits^n_{i = 1}E_i$. We also use $E(f)$ or $E(X/S)$ to denote the scheme $E$.

Let $\mathscr{F}_i$ be the kernel of the morphism $\mathcal{O}_X^{\ast} \to \mathcal{O}_{E_i}^{\ast}$ on $\etSite{X}$. Then $\mathscr{F}_i = 1 + \mathscr{K}_i \cdot \mathscr{I}_i$. Put $\mathscr{F} \defeq \prod\limits^n_{i = 1}\mathscr{F}_i$. Then we have an exact sequence of abelian sheaves:
\[0 \to \mathscr{F} \to (\mathcal{O}_X^{\ast})^n \to \mathcal{O}_E^{\ast} \to 0 \,.\]
We also use $\mathscr{F}(f)$ or $\mathscr{F}(X/S)$ to denote this abelian sheaf $\mathscr{F}$.

Let $\mathscr{P}$, $\mathscr{Q}$, $\theta$, $\vartheta$ and $\mathfrak{d}$ be the notations defined in Remark \thmref{668P8}. Obviously there is a canonical morphism $\mathbb{N}_S^n \to \mathscr{Q}$. Let $\gamma$ denote the composite
\[\mathbb{N}_X^n \to f^{-1}\mathscr{Q} \xrightarrow{\mathfrak{d}} \mathscr{P} \,.\]

\bigskip

We define a stack $\mathfrak{X}$ on $\etSite{X}$ as follows. For each \'{e}tale \mDash{X}{scheme} $U$, an object in $\mathfrak{X}(U)$ is a pair $(\mathscr{M},\sigma)$, where $\mathscr{M}$ is a fine saturated log structure on $U$ and $\sigma \colon \mathscr{M} \to \mathscr{P}|_U$ is a morphism of sheaves of monoids which induces an isomorphism $\overline{\mathscr{M}} \xrightarrow{\sim} \mathscr{P}|_U$ and makes the following diagram commutative:
\[\xymatrix{\mathscr{M} \ar[r] \ar[d]_{\sigma} & \mathcal{O}_U \ar[d] \\ \mathscr{P}|_U \ar[r]_-{\theta|_U} & \mathcal{O}_U/\mathcal{O}_U^{\ast}}\]
If $U' \to U$ is a morphism  of \'{e}tale \mDash{X}{schemes}, $(\mathscr{M},\sigma) \in \mathfrak{X}(U)$ and $(\mathscr{M}',\sigma') \in \mathfrak{X}(U')$, then a morphism of $(\mathscr{M}',\sigma')$ to $(\mathscr{M},\sigma)$ in $\mathfrak{X}$ lying above $U' \to U$ is an isomorphism $\varphi \colon \mathscr{M}' \xrightarrow{\sim} \mathscr{M}|_{U'}$ of log structures such that $\sigma|_{U'} \circ \varphi = \sigma'$.

We shall prove that $\mathfrak{X}$ is an \mDash{\mathscr{F}}{gerbe} (see Lemma \thmref{670P3}). The proof needs the following three simple lemmas.

\begin{noteLe}[670P1]
Let $X$ be a scheme, $\mathscr{M}$ a fine saturated log structure on $X$ and $\bar{x}$ a geometric point on $X$. Then there exists an \'{e}tale neighborhood $U$ of $\bar{x}$ and a fine saturated chart $P_U \to \mathscr{M}|_U$ such that the induced map $P \to \overline{\mathscr{M}}_{\bar{x}}$ is a bijection.
\end{noteLe}

\begin{noteLe}[670P0]
Let $X$ be a scheme and $\alpha \colon \mathscr{M} \to \mathcal{O}_X$ a fine log structure on $X$. Put $\mathscr{P} \defeq \overline{\mathscr{M}}$ and let $\bar{\alpha} \colon \mathscr{P} \to \mathcal{O}_X/\mathcal{O}_X^{\ast}$ be the morphism induced by $\alpha$. We define an abelian sheaf $\mathscr{A}$ on $\etSite{X}$ as follows: for every \'{e}tale \mDash{X}{scheme} $U$, $\mathscr{A}(U)$ is the set of morphisms $\sigma \colon \Groth{\mathscr{P}}|_U \to \mathcal{O}_U^{\ast}$ of abelian sheaves such that for any \'{e}tale \mDash{U}{scheme} $V$ and any section $s \in \mathscr{P}(V)$, we have $\sigma_V(s) \cdot t = t$, where $t \in \varGamma(V,\mathcal{O}_X)$ is a lifting of $\bar{\alpha}(s)$. Then there is a canonical isomorphism from $\mathscr{A}$ to the sheaf of automorphisms of log structures of $\mathscr{M}$ which induce identities on $\mathscr{P}$ defined as follows: for any \'{e}tale \mDash{X}{scheme} $U$ and any section $\sigma \in \varGamma(U,\mathscr{A})$, $\omega_U(\sigma)_V(s) = s \cdot \sigma_V\bigl(\bar{s})$, where $V$ is an \'{e}tale \mDash{U}{scheme} and $s \in \mathscr{M}(V)$.
\end{noteLe}

\begin{noteLe}[128B0]
Let $X$ be a scheme and $P$ a fine monoid. For each $i=1,2$, let $\alpha_i \colon P_X \to \mathcal{O}_X$ be a morphism of sheaves of monoids, $\mathscr{M}_i$ the log structure associating to $\alpha_i$, $\iota_i \colon P_X \to \mathscr{M}_i$ the induced morphism. Assume that there exists a morphism $\delta \colon P_X \to \mathcal{O}_X^{\ast}$ of sheaves of monoids such that $\alpha_1 = \delta \cdot \alpha_2$. Then there exists a unique isomorphism $\rho \colon \mathscr{M}_1 \xrightarrow{\sim} \mathscr{M}_2$ of log structures such that $\rho \circ \iota_1 = \delta \cdot \iota_2$.
\end{noteLe}

\begin{noteLe}[670P3]
$\mathfrak{X}$ is an \mDash{\mathscr{F}}{gerbe}.
\end{noteLe}

\begin{proof}
We have to verify the conditions in Definition \thmref{669P3}. (1) is obvious. Datum (b) and Condition (3) is by Lemma \thmref{670P0} and Theorem \thmref{667P0}.

For the condition (2), let $U$ be an \'{e}tale \mDash{X}{scheme}, $\alpha_1 \colon \mathscr{M}_1 \to \mathcal{O}_U$ and $\alpha_2 \colon \mathscr{M}_2 \to \mathcal{O}_U$ two objects in $\mathfrak{X}(U)$, $x$ a point on $U$. Put $P \defeq \mathscr{P}_{\bar{x}}$. By Lemma \thmref{670P1}, for each $i=1,2$, there exists an \'{e}tale neighborhood $V_i$, and a chart $\beta_i \colon P_{V_i} \to \mathscr{M}_i|_{V_i}$ which inducs identity on $P = \mathscr{P}_{\bar{x}}$. Put $V_3 \defeq V_1 \times_U V_2$. Since both
\[(\alpha_i \circ \beta_i)_{\bar{x}} \colon P \to \mathcal{O}_{U,\bar{x}}\]
are liftings of
\[\theta_{\bar{x}} \colon P \to \mathcal{O}_{U,\bar{x}}/\mathcal{O}_{U,\bar{x}}^{\ast} \,,\]
we have an \'{e}tale neighborhood $V \to V_3$ of $\bar{x}$ and a morphism $u \colon P_V \to \mathcal{O}_V^{\ast}$ of sheaves of monoids such that $\delta_1 = u \cdot \delta_2$, where
\[\delta_i \defeq (\alpha_i \circ \beta_i)|_V \colon P_V \to \mathcal{O}_V \,.\]
As $P_V \to \mathscr{M}_i|_V$ are charts, by Lemma \thmref{128B0} there exists an isomorphism $\varphi \colon \mathscr{M}_1 \xrightarrow{\sim} \mathscr{M}_2$ of log structures such that $\varphi \circ \beta_1|_V = u \cdot \beta_2|_V$. Thus $\varphi$ induces identity on $\mathscr{P}|_V$. So $\varphi$ is an isomorphism in $\mathfrak{X}(V)$.
\end{proof}

Obviously there exists an \'{e}tale covering $U \to X$, an object $\mathscr{M}$ in $\mathfrak{X}(U)$, and a morphism $\rho \colon \mathbb{N}^n_U \to \mathscr{M}$ which is a lifting of $\gamma|_U \colon \mathbb{N}^n_U \to \mathscr{P}|_U$. So there exists an \'{e}tale covering $V \to U \times_X U$ and an isomorphism $\phi \colon p_{10}^{\ast}(\mathscr{M}) \xrightarrow{\sim} p_{11}^{\ast}(\mathscr{M})$ of log structures on $V$ (here the notations $(V/U)_k$ and $p_{ki} \colon (V/U)_k \to (V/U)_{k-1}$ are defined as in \S\ref{Se:hyper:1}). Hence there exists an element
\[u = (\Sequ{n}{u}) \in (\mathcal{O}_X^{\ast})^n(V)\]
such that
\[\phi \circ p_{10}^{\ast}(\rho) = u \cdot p_{11}^{\ast}(\rho)\]
and a cocycle $g \in Z^2(V/U,\mathscr{F})$ such that
\[g \circ p_{21}^{\ast}(\phi) = p_{22}^{\ast}(\phi) \circ p_{20}^{\ast}(\phi) \,.\]
We have
\[[\mathfrak{X}] = [g] \in H^2(\etSite{X},\mathscr{F}) \,.\]
By \eqref{E:hyper1e}, we have
\begin{align*}
p_{10} \circ p_{21} &= p_{10} \circ p_{20} \,, \\
p_{10} \circ p_{22} &= p_{11} \circ p_{20} \,, \\
p_{11} \circ p_{22} &= p_{11} \circ p_{21} \,.
\end{align*}
We also have
\begin{align*}
& \hspace{1.4em} g \cdot p_{21}^{\ast}(u) \cdot (p_{11} \circ p_{21})^{\ast}(\rho) \\
&= g \cdot p_{21}^{\ast}(\phi) \circ (p_{10} \circ p_{21})^{\ast}(\rho) \\
&= p_{22}^{\ast}(\phi) \circ p_{20}^{\ast}(\phi) \circ (p_{10} \circ p_{20})^{\ast}(\rho) \\
&= p_{20}^{\ast}(u) \cdot p_{22}^{\ast}(\phi) \circ (p_{11} \circ p_{20})^{\ast}(\rho) \\
&= p_{20}^{\ast}(u) \cdot p_{22}^{\ast}(\phi) \circ (p_{10} \circ p_{22})^{\ast}(\rho) \\
&= p_{20}^{\ast}(u) \cdot p_{22}^{\ast}(u) \cdot (p_{11} \circ p_{22})^{\ast}(\rho) \\
&= p_{20}^{\ast}(u) \cdot p_{22}^{\ast}(u) \cdot (p_{11} \circ p_{21})^{\ast}(\rho) \,.
\end{align*}
Thus $g = p_{20}^{\ast}(u) \cdot p_{22}^{\ast}(u) \cdot p_{21}^{\ast}(u)^{-1}$. Since the image of $g$ in $\mathcal{O}_E^{\ast}$ is equal to $1$, we see that $u$ determinates an element in $H^1(\etSite{E},\mathcal{O}_E^{\ast})$. So we obtain an invertible \mDash{\mathcal{O}_E}{module}, which depends only on the morphism $f \colon X \to S$. We denote it by $\mathscr{L}(f)$ or $\mathscr{L}(X/S)$.

\begin{noteDe}
A \emph{semistable log structure} on $X$ is an object $(\mathscr{M},\sigma)$ in $\mathfrak{X}(X)$ such that there is a morphism $\rho \colon \mathbb{N}_X^n \to \mathscr{M}$ which lifts the morphism $\gamma \colon \mathbb{N}_X^n \to \mathscr{P}$.
\end{noteDe}

\begin{noteTh}[670P2]\
\begin{enumerate}
\item There exists a semistable log structure on $X$ if and only if $\mathscr{L}(f) \cong \mathcal{O}_E$.
\item The semistable log structure on $X$ is unique (up to isomorphism) if it exists.
\end{enumerate}
\end{noteTh}

\begin{proof}
(1) If semistable log structures on $X$ exist, obviously $\mathscr{L}(f) \cong \mathcal{O}_E$.

Assume that $\mathscr{L}(f) \cong \mathcal{O}_E$. By above argument, $[\mathfrak{X}]$ is the image of $\mathscr{L}(f)$ under the connecting map
\[H^1(\etSite{E},\mathcal{O}_E^{\ast}) \to H^2(\etSite{X},\mathscr{F}) \,.\]
Thus $[\mathfrak{X}] = 0$. By Theorem \thmref{669P0}, there exists an element $\mathscr{M} \in \mathfrak{X}(X)$. Let $U \to X$ be an \'{e}tale covering such that there exists a lifting $\rho \colon \mathbb{N}_U^n \to \mathscr{M}|_U$ of $\gamma|_U \colon \mathbb{N}_U^n \to \mathscr{P}|_U$. Put $V \defeq U \times_S U$ and let $u \in (\mathcal{O}_X^{\ast})^n(V)$ such that $p_{10}^{\ast}(\rho) = u \cdot p_{11}^{\ast}(\rho)$. Let $\bar{u}$ be the image of $u$ in $\mathcal{O}_E^{\ast}(V)$. As $\mathscr{L}(f) \cong \mathcal{O}_E$ is represented by $[\bar{u}]$, there exists an element $v' \in \mathcal{O}_E(U)$ such that $\bar{u} = p_{10}^{\ast}(v') \cdot p_{11}^{\ast}(v')^{-1}$. By contracting $U$ suitably, we may choose a lifting $v \in (\mathcal{O}_X^{\ast})^n(U)$ of $v'$. Then
\[u = p_{10}^{\ast}(v) \cdot p_{11}^{\ast}(v)^{-1} \cdot g\]
for some $g \in \mathscr{F}(V)$. As $\partial^2(g) = \partial^2(u) = 1$, $g$ is a cocycle. Put
\[\rho_1 \defeq v^{-1} \cdot \rho \colon \mathbb{N}_U^n \to \mathscr{M}|_U \,.\]
Then
\[p_{10}^{\ast}(\rho_1) = g \cdot p_{11}^{\ast}(\rho_1) \,.\]
By Lemma \thmref{669P1}, there exists an object $\mathscr{M}'$ in $\mathfrak{X}(X)$ and an isomorphism $\phi \colon \mathscr{M}'|_U \to \mathscr{M}|_U$ in $\mathfrak{X}(U)$ such that $g^{-1} \cdot p_{10}^{\ast}(\phi) = p_{11}^{\ast}(\phi)$. Put $\rho_2 \defeq \phi^{-1} \circ \rho_1$. Then $p_{10}^{\ast}(\rho_2) = p_{11}^{\ast}(\rho_2)$. So there exists a morphism $\rho' \colon \mathbb{N}_X^n \to \mathscr{M}'$ such that $\rho'|_U = \rho_2$. Obviously $\rho'$ is a lift of $\gamma \colon \mathbb{N}_X^n \to \mathscr{P}$. Thus $\mathscr{M}'$ is a semistable log structure on $X$.

(2) Let $\mathscr{M}_1$ and $\mathscr{M}_2$ be two semistable log structures on $X$. For $i=1,2$, let $\rho_i \colon \mathbb{N}_X^n \to \mathscr{M}_i$ be liftings of $\gamma \colon \mathbb{N}_X^n \to \mathscr{P}$. By Theorem \thmref{669P2}, there exist \'{e}tale coverings $U \to X$ and $V \to U \times_X U$, a cocycle $g \in Z^1(V/U,\mathscr{F})$ and an isomorphism $\phi \colon \mathscr{M}_2|_U \xrightarrow{\sim} \mathscr{M}_1|_U$ such that $g \cdot p_{10}^{\ast}(\phi) = p_{11}^{\ast}(\phi)$. Let $\delta \in (\mathcal{O}_X^{\ast})^n(U)$ such that $\phi \circ \rho_2|_U = \delta \cdot \rho_1|_U$. Then we have
\begin{align*}
g \cdot p_{10}^{\ast}(\delta) \cdot \rho_1|_V &= g \cdot p_{10}^{\ast}(\delta \cdot \rho_1|_U) \\
&= g \cdot p_{10}^{\ast}(\phi \circ \rho_2|_U) \\
&= \bigl(g \cdot p_{10}^{\ast}(\phi)\bigr)  \circ \rho_2|_V \\
&= p_{11}^{\ast}(\phi) \circ \rho_2|_V \\
&= p_{11}^{\ast}(\delta) \cdot \rho_1|_V \,.
\end{align*}
So $g = p_{11}^{\ast}(\delta) \cdot p_{10}^{\ast}(\delta)^{-1}$, i.e., $[g] = 0$. Therefore $\mathscr{M}_1 \cong \mathscr{M}_2$ in $\mathfrak{X}(X)$.
\end{proof}

\begin{noteRe}
Note that the isomorphisms between semistable log structures may not be unique. So this kind of structure is not canonical.
\end{noteRe}

\begin{noteTh}
Assume that all $\Sequ{n}{a}$ are regular elements in $\varGamma(S,\mathcal{O}_S)$ (i.e., $(0:a_i) = 0$ for all $i \in [1,n]$). Then $\mathscr{L}(f) \cong \mathcal{O}_E$, i.e., there exists a semistable log structure on $X$.
\end{noteTh}

\begin{proof}
Note that
\[\mathscr{F} = \prod^n_{i=1}(1 + \mathscr{K}_i \cdot \mathscr{I}_i) = 0 \,.\]
So there exists an object $\mathscr{M}$ in $\mathfrak{X}(X)$. Obviouly there exists an \'{e}tale covering $\{U_{\lambda}\}_{\lambda \in \Lambda}$ of $X$ such that for each $\lambda \in \Lambda$, there exists a lifting $\rho_{\lambda} \colon \mathbb{N}_{U_{\lambda}}^n \to \mathscr{M}|_{U_{\lambda}}$ of $\gamma|_{U_{\lambda}} \colon \mathbb{N}_{U_{\lambda}}^n \to \mathscr{P}|_{U_{\lambda}}$ such that the composite morphism
\[\mathbb{N}_{U_{\lambda}}^n \xrightarrow{\rho_{\lambda}} \mathscr{M}|_{U_{\lambda}} \to \mathcal{O}_{U_{\lambda}}\]
is equal to
\[\mathbb{N}_{U_{\lambda}}^n \to \mathcal{O}_{U_{\lambda}} \,, \qquad \varepsilon_i \to a_i \,.\]
Since $\Sequ{n}{a}$ are regular elements in $\varGamma(S,\mathcal{O}_S)$ and $f \colon X \to S$ is flat, $\Sequ{n}{a}$ are regular elements in $\varGamma(X,\mathcal{O}_X)$ too. So on each $U_{\lambda\mu}$, we have $\rho_{\lambda}|_{U_{\lambda\mu}} = \rho_{\mu}|_{U_{\lambda\mu}}$. Thus $\{\rho_{\lambda}\}$ can be glued to obtain a global lifting $\rho \colon \mathbb{N}_X^n \to \mathscr{M}$ of $\gamma \colon \mathbb{N}_X^n \to \mathscr{P}$.
\end{proof}

\begin{noteTh}[670P4]
Let $\mathscr{M}$ be a semistable log structure on $X$ and $\rho \colon \mathbb{N}^n \to \mathscr{M}$ a lifting of $\gamma \colon \mathbb{N}^n_X \to \mathscr{P}$. By Lemma \thmref{128B0}, $\rho$ induces a morphism $\varphi \colon f^{\ast}\mathscr{N} \to \mathscr{M}$ of log structures. Then
\[(f,\varphi) \colon (X,\mathscr{M}) \to (Y,\mathscr{N})\]
is log smooth and integral.
\end{noteTh}

\begin{proof}
Let $x$ be a point on $X$ and $y \defeq f(x)$. Let
\[(U,V; T_{11},\ldots,T_{1q_1}; \ldots; T_{p1},\ldots,T_{pq_p}; a_1,\ldots,a_l)\]
be a \rlC of $f$ at $x$ and put $m \defeq \sum\limits^p_{i=1}q_i + l - p$. Let $g \colon U \to V$ be the induced morphism. Then
\[(g,\varphi_{U/V}) \colon (U,\mathscr{M}_U) \to (V,\mathscr{N}_V)\]
is log smooth, where $\mathscr{M}_U$, $\mathscr{N}_V$ and $\varphi_{U/V}$ are defined in Remark \thmref{668P6}. By Lemma \thmref{128B0}, we may contract $U/V$ suitably such that there exists two isomorphism $\sigma \colon \mathscr{M}_U \xrightarrow{\sim} \mathscr{M}|_U$ and $\tau \colon \mathscr{N}_V \xrightarrow{\sim} \mathscr{N}|_V$ of log structures. Put $\mathscr{M}_0 \defeq \mathscr{M}_U$, $\mathscr{N}_0 \defeq \mathscr{N}_V$, $\varphi_1 \defeq \varphi_{U/V}$, $\varphi_2 \defeq \sigma^{-1} \circ \varphi|_U \circ f^{\ast}(\tau)$. Let $\alpha \defeq \alpha_U \colon \mathbb{N}^m_U \to \mathscr{M}_0$, $\beta \defeq \beta_V \colon \mathbb{N}^l_V \to \mathscr{N}_0$ and $\gamma \colon \mathbb{N}^l \to \mathbb{N}^m$ be the notations as in Remark \thmref{668P6}. Now $(g,\varphi_1)$ is log smooth. We shall use the definition of log smoothness to prove that $(g,\varphi_2)$ is also log smooth. Let $(T_0,\mathscr{T}_0)$ and $(T,\mathscr{T})$ be fine log schemes, $(T_0,\mathscr{T}_0) \to (T,\mathscr{T})$ a thickening of order one (Cf.~\cite[3.5]{FKato2}), $(t_0,\psi_0) \colon (T_0,\mathscr{T}_0) \to (U,\mathscr{M}_0)$ and $(t,\psi) \colon (T,\mathscr{T}) \to (V,\mathscr{N}_0)$ be morphisms of log schemes which makes a commutative diagram:
\[\xymatrix@C+2em{(T_0,\mathscr{T}_0) \ar[r]^{(t_0,\psi_0)} \ar[d] & (U,\mathscr{M}_0) \ar[d]^{(g,\varphi_2)} \\ (T,\mathscr{T}) \ar[r]_{(t,\psi)} & (V,\mathscr{N}_0)}\]
For each $i=1,2$, let $\rho_i$ denote the composite morphism
\[\mathbb{N}^l_U \xrightarrow{g^{\ast}(\beta)} g^{\ast}\mathscr{N}_0 \xrightarrow{\varphi_i} \mathscr{M}_0 \,.\]
Then there exists a section $u \in (\mathcal{O}_U^{\ast})^n(U)$ such that $\rho_2 = u \cdot \rho_1$. The composite morphism
\[\mathbb{N}^l_T \xrightarrow{t^{\ast}(\beta)} t^{\ast}\mathscr{N}_0 \xrightarrow{\psi} \mathscr{T} \to \mathscr{T}_0\]
is equal to $\psi_0 \circ t_0^{\ast}(\rho_2)$. Let $v \in (\mathcal{O}_T^{\ast})^n(T)$ be a lift of
\[t_0^{\ast}(u) \in (\mathcal{O}_{T_0}^{\ast})^n(T_0) \,.\]
Then there is a morphism $\psi' \colon t^{\ast}\mathscr{N}_0 \to \mathscr{T}$ of log structures such that
\[\psi' \circ t^{\ast}(\beta) = v^{-1} \cdot \bigl(\psi \circ t^{\ast}(\beta)\bigr) \,.\]
So the composite morphism
\[\mathbb{N}^l_T \xrightarrow{t^{\ast}(\beta)} t^{\ast}\mathscr{N}_0 \xrightarrow{\psi'} \mathscr{T} \to \mathscr{T}_0\]
is equal to $\psi_0 \circ t_0^{\ast}(\rho_1)$, which shows that
\[\xymatrix@C+2em{(T_0,\mathscr{T}_0) \ar[r]^{(t_0,\psi_0)} \ar[d] & (U,\mathscr{M}_0) \ar[d]^{(g,\varphi_1)} \\ (T,\mathscr{T}) \ar[r]_{(t,\psi')} & (V,\mathscr{N}_0)}\]
is commutative. As $(g,\psi_1)$ is log smooth, by replace $T$ with an \'{e}tale covering, we may assume that there is a morphism
\[(h,\xi_1) \colon (T,\mathscr{T}) \to (U,\mathscr{M}_0)\]
of log schemes which makes the following diagram
\[\xymatrix@C+2em{(T_0,\mathscr{T}_0) \ar[r]^{(t_0,\psi_0)} \ar[d] & (U,\mathscr{M}_0) \ar[d]^{(g,\varphi_1)} \\ (T,\mathscr{T}) \ar[ur]^{(h,\xi_1)} \ar[r]_{(t,\psi')} & (V,\mathscr{N}_0)}\]
commutative. By the following Lemma \thmref{670P5}, there exists a section $w \in (\mathcal{O}_T^{\ast})^m(T)$ which makes a commutative diagram:
\[\xymatrix@C+2em{\mathbb{N}^l_T \ar[r]^-{v \cdot h^{\ast}(u)^{-1}} \ar[d]_{\gamma_T} & \mathcal{O}_T^{\ast} \\ \mathbb{N}^m_T \ar[ur]_-{w}}\]
Let $\xi_2 \colon h^{\ast}\mathscr{M}_0 \to \mathscr{T}$ be the morphism of log structures satisfying that $\xi_2 \circ h^{\ast}(\alpha) = w \cdot \bigl(\xi_1 \cdot h^{\ast}(\alpha)\bigr)$. Then we have a commutative diagram:
\[\xymatrix@C+2em{(T_0,\mathscr{T}_0) \ar[r]^{(t_0,\psi_0)} \ar[d] & (U,\mathscr{M}_0) \ar[d]^{(g,\varphi_2)} \\ (T,\mathscr{T}) \ar[ur]^{(h,\xi_2)} \ar[r]_{(t,\psi)} & (V,\mathscr{N}_0)}\]
Thus $(g,\psi_2)$ is log smooth.
\end{proof}

\begin{noteLe}[670P5]
Let $A$ be a ring, $I$ an ideal of $A$ such that $I^2 = (0)$, $u \in 1+I$, $\Sequ{n}{e}$ be positive integers which are invertible in $A$. Then there exists elements $\Sequ{n}{v} \in 1+I$ such that $u = \prod\limits^n_{i=1}v_i^{e_i}$.
\end{noteLe}

The following theorem is obvious.

\begin{noteTh}
Let $S'$ be a locally noetherian scheme and $S' \to S$ a flat morphism. Put $X' \defeq X \times_S S'$ and let $f' \colon X' \to S'$ be the projection. Then
\begin{enumerate}
\item $f'$ satisfies all these conditions of $f$ mentioned at the beginning of this section.
\item For each $i \in [1,n]$, let $b_i$ be the image of $a_i$ in $\varGamma(S', \mathcal{O}_{S'})$. Assume that $b_i$ is invertible for $i \in [m+1,n]$ and $b_i$ is not invertible for $[1,m]$. Then
\[E(f') = \Bigl(\coprodT^m_{i=1}E_i \Bigr) \times_S S' \,.\]
\item $\mathscr{L}(f')$ is isomorphic to the inverse image of $\mathscr{L}$ under the canonical morphism $E(f') \to E(f)$.
\end{enumerate}
\end{noteTh}

\section{Global Cases}
\label{Se:glo}

Let $X$ and $S$ be locally noetherian schemes, $f \colon X \to S$ a surjective proper \wnC morphism without powers such that $f$ satisfies the condition (\dag) in \S\ref{Se:global:3} and every fiber of $f$ is geometrically connected. By Lemma \thmref{668P7} (4), $f$ also satisfies the condition (\ddag) in \S\ref{Se:global:3}.

Let $\mathscr{P}$, $\mathscr{Q}$, $\theta$, $\vartheta$ and $\mathfrak{d}$ be the notations defined in Remark \thmref{668P8}.

For every point $y \in S$, we write
\begin{align*}
E_{\bar{y}} & \defeq E\bigl(X \times_S \Spec \mathcal{O}_{S,\bar{y}} \big/ \Spec \mathcal{O}_{S,\bar{y}}\bigr) \,, \\
\mathscr{L}_{\bar{y}} & \defeq \mathscr{L}\bigl(X \times_S \Spec \mathcal{O}_{S,\bar{y}} \big/ \Spec \mathcal{O}_{S,\bar{y}}\bigr) \,.
\end{align*}

\begin{noteLe}
Let $y \in S$. If $\mathscr{L}_{\bar{y}}$ is trivial, then there exists an element $V_0 \in N_0(y)$ such that for all elements $V \geqslant V_0$ in $N_0(y)$, $\mathscr{L}(X_V/V)$ is trivial.
\end{noteLe}

\begin{proof}
See \cite[(8.5.2.5)]{AGro2}.
\end{proof}

\begin{noteCo}
Let $y \in S$. If $\mathscr{L}_{\bar{y}}$ is trivial, then there exists an open neighborhood $V$ of $y$ such that for all $z \in V$, $\mathscr{L}_{\bar{z}}$ is trivial.
\end{noteCo}

\begin{noteLe}[671P0]
Let $f \colon X \to Y$ be a proper and flat morphism of locally noetherian schemes such that every fiber of $f$ is geometrically reduced and geometrically connected. Then the canonical morphism $\mathcal{O}_Y \to f_{\ast}\mathcal{O}_X$ is isomorphic.
\end{noteLe}

\begin{proof}
See \cite[(7.8.7) and (7.8.8)]{AGro1}.
\end{proof}

\begin{noteLe}[671P2]
Let $R$ be a noetherian local ring with maximal ideal $\mathfrak{m}$, $a \in \mathfrak{m}$, $\mathfrak{a} = (0: a)$, $n \geqslant 2$ an integer,
\[A = R[[\Sequ{n}{T}]]\]
a ring of power series over $R$, $I$ the ideal of $A$ generated by
\[T_2 \cdots T_n, \ldots, T_1 \cdots \hat{T}_i \cdots T_n, \ldots, T_1 \cdots T_{n-1}\]
and
\[J \defeq (T_1T_2 \cdots T_n - a) + \mathfrak{a} \cdot I \,.\]
Then $J \cap R = (0)$.
\end{noteLe}

\begin{proof}
Let $b \in J \cap R$ and put
\[b = (T_1T_2 \cdots T_n - a) \cdot F_0 + \sum^n_{i=1}T_1 \cdots \hat{T}_i \cdots T_n \cdot F_i \,,\]
where $F_i \in A$ for all $i \in [0,n]$ and for every $i \in [1,n]$, all coefficients of $F_i$ are contained in $\mathfrak{a}$. For each $i \in [1,n]$, put $F_i = G_i + T_i \cdot N_i$, where all monomials in $G_i$ do not contain $T_i$ and all coefficients of $G_i$ and $N_i$ are contained in $\mathfrak{a}$. Then we have
\[T_1 \cdots \hat{T}_i \cdots T_n \cdot F_i =T_1 \cdots \hat{T}_i \cdots T_n \cdot G_i + (T_1T_2 \cdots T_n - a) \cdot N_i \,.\]
Put $G_0 \defeq F_0 + \sum\limits^n_{i=1}N_i$. Then
\begin{equation}\label{E:glo0e}
b = (T_1T_2 \cdots T_n - a) \cdot G_0 + \sum^n_{i=1}T_1 \cdots \hat{T}_i \cdots T_n \cdot G_i \,.
\end{equation}
For each $q \in \mathbb{N}$, let $c_q$ denote the coefficient of $(T_1T_2 \cdots T_n)^q$ in $G_0$. By comparing the coefficient of $(T_1T_2 \cdots T_n)^q$ in \eqref{E:glo0e}, we have $b = -ac_0$ and $c_{q-1} = ac_q$ for all $q \geqslant 1$. Thus
\begin{gather*}
b \in \bigcap^{\infty}_{q=1}(a^q) \subseteq \bigcap^{\infty}_{q=1}\mathfrak{m}^q = (0) \,. \qedhere
\end{gather*}
\end{proof}

\begin{noteLe}[671P1]
For all points $y \in S$ and $V \in N_0(y)$, we have $(f_V)_{\ast}\bigl(\mathscr{F}(X_V/V)\bigr) = 1$.
\end{noteLe}

\begin{proof}
Let $\Sequ{n}{D}$ be the connected components of $D_V$, $\Sequ{n}{a} \in \mathcal{O}_S(V)$ the corresponding sections. For each $i \in [1,n]$, let $\mathscr{K}_i$ and $\mathscr{I}_i$ be ideal sheaves on $V$ and $X_V$ respectively defined in \S\ref{Se:loc}, and put $\mathscr{J}_i \defeq \mathscr{K}_i \cdot \mathscr{I}_i$. Then $\mathscr{F}(X_V/V) = \prod\limits^n_{i=1}(1+\mathscr{J}_i)$. So we have only to prove $(f_V)_{\ast}\mathscr{J}_i = (0)$. Let $W$ be an open subset of $V$ and $b \in \varGamma(X_W, \mathscr{I}_i)$. By Lemma \thmref{671P0}, we have $b \in \varGamma(W,\mathcal{O}_S)$. Suppose that $W_b \neq \emptyset$. As $b|_{W_b} = b|_{f^{-1}(W_b)}$ is invertible, $\mathscr{J}_i|_{f^{-1}(W_b)} = (1)$, so $\mathscr{K}_i|_{W_b} = (1)$ and $\mathscr{I}_i|_{f^{-1}(W_b)} = (1)$. As $\mathscr{K}_i = (0:a)\sptilde$, $a_i = 0$. Hence $(D_i)_{W_b} \neq \emptyset$, i.e., $\mathscr{I}_i|_{f^{-1}(W_b)} \neq (1)$, a contradiction. Hence $W_b = \emptyset$. Thus for any $w \in W$, $b_{\bar{w}}$ is contained in the maximal ideal of $\mathcal{O}_{S,\bar{w}}$. By Lemma \thmref{671P2}, $b_{\bar{w}} = 0$. So $b = 0$.
\end{proof}

\begin{noteDe}
A \emph{semistable log structure} for $f$ is of the form $(\mathscr{M}, \mathscr{N}, \sigma, \tau, \varphi)$, where $\mathscr{M}$ and $\mathscr{N}$ are fine saturated log structures on $X$ and $S$ respectively, $\varphi \colon f^{\ast}\mathscr{N} \to \mathscr{M}$ is a morphism of log structures on $X$, $\sigma \colon \mathscr{M} \to \mathscr{P}$ and $\tau \colon \mathscr{N} \to \mathscr{Q}$ are morphisms of sheaves of monoids, such that $\sigma$ and $\tau$ induce isomorphisms $\bar{\sigma} \colon \overline{\mathscr{M}} \xrightarrow{\sim} \mathscr{P}$ and $\bar{\tau} \colon \overline{\mathscr{N}} \xrightarrow{\sim} \mathscr{Q}$, and the following three diagrams are commutative:
\[\xymatrix{\mathscr{N} \ar[r] \ar[d]_-{\tau} & \mathcal{O}_S \ar[d] \\ \mathscr{Q} \ar[r]_-{\vartheta} & \mathcal{O}_S/\mathcal{O}_S^{\ast}} \qquad
\xymatrix{\mathscr{M} \ar[r] \ar[d]_-{\sigma} & \mathcal{O}_X \ar[d] \\ \mathscr{P} \ar[r]_-{\theta} & \mathcal{O}_X/\mathcal{O}_X^{\ast}} \qquad
\xymatrix{f^{-1}\mathscr{N} \ar[r] \ar[d]_{f^{-1}\tau} & \mathscr{M} \ar[d]^{\sigma} \\ f^{-1}\mathscr{Q} \ar[r]_-{\mathfrak{d}} & \mathscr{P}}\]
\end{noteDe}

The following two theorems are the main results of this papers.

\begin{noteTh}[671P3]\
\begin{enumerate}
\item There exists a semistable log structure for $f$ if and only if for every point $y \in S$, $\mathscr{L}_{\bar{y}}$ is trivial on $E_{\bar{y}}$.
\item Let $(\mathscr{M}_1, \mathscr{N}_1, \sigma_1, \tau_1, \phi_1)$ and $(\mathscr{M}_2, \mathscr{N}_2, \sigma_2, \tau_2, \phi_2)$ be two semistable log structures for $f$. Then there exists isomorphism $\varphi \colon \mathscr{M}_1 \xrightarrow{\sim} \mathscr{M}_2$ and $\psi \colon \mathscr{N}_1 \xrightarrow{\sim} \mathscr{N}_2$ of log structures such that $\varphi \circ \phi_1 = \phi_2 \circ f^{\ast}\psi$, $\sigma_2 \circ \varphi = \sigma_1$ and $\tau_2 \circ \psi = \tau_1$. Furthermore, such pair $(\varphi,\psi)$ is unique.
\end{enumerate}
\end{noteTh}

\begin{proof}
(2) Let $y$ be a point on $S$, $V \in N_0(y)$, and let
\[\Sequ{n}{a} \in \varGamma(V,\mathcal{O}_S)\]
be sections satisfying the condition (2) in Notation \thmref{668P9}. Clearly $\mathscr{N}_{\bar{y}} = \mathbb{N}^n$. By Lemma \thmref{670P1}, we may contract $V$ suitably to make both $\mathscr{N}_1$ and $\mathscr{N}_2$ isomorphic to the log structure associated to
\[\beta_0 \colon \mathbb{N}^n_V \to \mathcal{O}_V \,, \qquad \varepsilon_i \mapsto a_i \,,\]
where $\{\Sequ{n}{\varepsilon}\}$ is a basis of $\mathbb{N}^n$. In other words, we have an isomorphisms $\psi_0 \colon \mathscr{N}_1|_V \xrightarrow{\sim} \mathscr{N}_2|_V$ of log structures, and charts $\beta_i \colon \mathbb{N}^n \to \mathscr{N}_i|_V$ such that $\psi_0 \circ \beta_1 = \beta_2$, $\tau_2|_V \circ \psi_0 = \tau_1|_V$, and the diagram
\[\xymatrix@R-2ex{& \mathbb{N}^n \ar[dl] \ar[dd]^{\beta_i} \ar[dr]^{\beta_0} \\ \mathscr{Q}|_V & & \mathcal{O}_V \\ & \mathscr{N}_i|_V \ar[ul]^{\tau_i|_V} \ar[ur]}\]
is commutative. Since the composite morphisms
\[\mathbb{N}^n_{X_V} \xrightarrow{f^{\ast}_V(\beta_i)} f_V^{\ast}(\mathscr{N}_i|_V) \xrightarrow{\phi_i|_{X_V}} \mathscr{M}_i|_{X_V}\]
are liftings of
\[\mathbb{N}^n_{X_V} \to f^{-1}\mathscr{Q}|_{X_V} \xrightarrow{\mathfrak{d}|_{X_V}} \mathscr{P}|_{X_V} \,,\]
we see that $(\mathscr{M}_i|_{X_V}, \sigma_i|_{X_V})$ are semistable log structures on $X_V$. So, by Theorem \thmref{670P2} (2), there exists an isomorphism $\varphi_V \colon \mathscr{M}_1|_{X_V} \xrightarrow{\sim} \mathscr{M}_2|_{X_V}$ such that $\varphi_V \circ \sigma_1|_{X_V} = \sigma_2|_{X_V}$. Obviously there exists a section
\[v = (\Sequ{n}{v}) \in (\mathcal{O}_X^{\ast})^n(X_V)\]
such that
\[\varphi_V \circ \bigl((\phi_1|_{X_V}) \circ f^{\ast}_V(\beta_1)\bigr) = v \cdot \bigl((\phi_2|_{X_V}) \circ f^{\ast}_V(\beta_2)\bigr) \,.\]
Taking composite of both sides of above equality with the morphism $\mathscr{M}_2|_{X_V} \to \mathcal{O}_{X_V}$, we have $v_i a_i = a_i$ for each $i \in [1,n]$. Applying Lemma \thmref{671P0} to the morphism $X_V \to V$, we have $v_i \in \mathcal{O}_S^{\ast}(V)$. By Lemma \thmref{128B0}, there is an isomorphism $\psi_V \colon \mathscr{N}_1|_V \xrightarrow{\sim} \mathscr{N}_2|_V$ of log structures such that $\psi_V \circ \beta_1 = v \cdot \beta_2$. Thus $\varphi_V \circ (\phi_1|_{X_V}) = (\phi_2|_{X_V}) \circ f_V^{\ast}(\psi_V)$.

Suppose that there exists another pair of isomorphisms $\varphi' \colon \mathscr{M}_1|_{X_V} \xrightarrow{\sim} \mathscr{M}_2|_{X_V}$ and $\psi' \colon \mathscr{N}_1|_V \xrightarrow{\sim} \mathscr{N}_2|_V$ of log structures such that $\varphi' \circ (\phi_1|_{X_V}) = (\phi_2|_{X_V}) \circ f_V^{\ast}(\psi')$, $(\sigma_2|_{X_V}) \circ \varphi' = \sigma_1|_{X_V}$ and $\tau_2|_V \circ \psi' = \tau_1|_{X_V}$. By Lemma \thmref{670P3},
\[\varphi'^{-1} \circ \varphi|_{X_V} \in \varGamma\bigl(X_V, \mathscr{F}(X_V/V)\bigr) \,.\]
By Lemma \thmref{671P1}, $\varphi'^{-1} \circ \varphi|_{X_V} = \mathrm{id}$, i.e., $\varphi_V = \varphi'$. Thus $(\phi_2|_{X_V}) \circ f_V^{\ast}(\psi_V) = (\phi_2|_{X_V}) \circ f_V^{\ast}(\psi')$. It is easy to show that $\phi_2$ is injective. So $f_V^{\ast}(\psi_V) = f_V^{\ast}(\psi')$. Since $f_V$ is faithfully flat, we get $\psi_V = \psi'$.

Now we can glue these $(\varphi_V,\psi_V)$ to a pair of isomorphism of log structures $(\varphi,\psi)$.

(1) is by (2) and Theorem \thmref{670P2}.
\end{proof}

\begin{noteTh}
Let $(\mathscr{M}, \mathscr{N}, \sigma, \tau, \phi)$ be a semistable log structure for $f$. Then
\[(f,\phi) \colon (X,\mathscr{M}) \to (Y,\mathscr{N})\]
is log smooth and integral.
\end{noteTh}

\begin{proof}
The conclusion is a consequence of Theorem \thmref{670P4}.
\end{proof}

\section{Properties under Base Change}
\newcommand{\nA}{$(\mathrm{N}_1)$\xspace}
\newcommand{\nB}{$(\mathrm{N}_2)$\xspace}
\newcommand{\nC}{$(\mathrm{N}_3)$\xspace}

\begin{noteDe}[672P4]
Let $f \colon X \to S$ be a morphism of locally noetherian schemes.
\begin{enumerate}
\item We say that $f$ satisfies \nA if it is surjective, proper, \wnC without powers, and all fibers of $f$ are geometrically connected.
\item We say that $f$ satisfies \nB if it satisfies \nA and the condition \rmB{\dag} in \S\ref{Se:global:3}.
\item We say that $f$ satisfies \nC if it satisfies \nB and for every point $y \in S$, the invertible sheaf $\mathscr{L}_{\bar{y}}$ on $E_{\bar{y}}$ defined in \S\ref{Se:glo} is trivial.
\end{enumerate}
\end{noteDe}

\subsection{Properties under fibred products}

Let $S$, $X$ and $Y$ be locally noetherian schemes, $f \colon X \to S$ and $g \colon Y \to S$ two morphisms. For an \mDash{S}{scheme} $Z$ which satisfies \nB and a point $s$ on $S$, we use $E_{\bar{s}}(Z/S)$ and $\mathscr{L}_{\bar{s}}(Z/S)$ to denote the notations $E_{\bar{s}}$ and $\mathscr{L}_{\bar{s}}$ defined in \S\ref{Se:glo} for preciseness, and write
\[Z(\bar{s}) \defeq Z \times_S \Spec \mathcal{O}_{S,\bar{s}} \,.\]

\begin{noteTh}
Assume that $f$ and $g$  satisfies \nA. Then $X \times_S Y \to S$ satisfies \nA. Furthermore we have
\[D\bigl((X \times_S Y)/S\bigr) = \bigl(D(X/S) \times_S Y\bigr) \coprodT\, \bigl(X \times_S D(Y/S)\bigr) \,.\]
\end{noteTh}

\begin{noteTh}
Assume that $f$ and $g$ satisfies \nB. Then $X \times_S Y \to S$ satisfies \nB. Furthermore if $s \in S$, then
\begin{align*}
E_{\bar{s}}(X \times_S Y) &= \bigl(E_{\bar{s}}(X) \times_S Y(\bar{s})\bigr) \coprodT\, \bigl(X(\bar{s}) \times_S E_{\bar{s}}(Y)\bigr) \,, \\
\mathscr{L}_{\bar{s}}(X \times_S Y) &= \bigl(\mathscr{L}_{\bar{s}}(X) \otimes_S \mathcal{O}_{Y(\bar{s})}\bigr) \coprodT\, \bigl(\mathcal{O}_{X(\bar{s})} \otimes_S \mathscr{L}_{\bar{s}}(Y)\bigr) \,.
\end{align*}
\end{noteTh}

\begin{noteTh}
If $f$ and $g$ satisfies \nC, so does $X \times_S Y \to S$.
\end{noteTh}

\subsection{Properties under base extension}

Let
\[\xymatrix{X' \ar[r]^{p} \ar[d]_{f'} \CartS & X \ar[d]^{f} \\ S' \ar[r]_{q} & S}\]
be a Cartesian square of locally noetherian schemes.

\begin{noteTh}
If $f$ satisfies \nA, so does $f'$.
\end{noteTh}

\begin{noteLe}
Assume that $f$ satisfies \nA. Let $y'$ be a point on $S'$, $y \defeq q(y')$. Fix a \mDash{\kappa(y)}{embedding} of $\kappa(y)_s$ into $\kappa(y')_s$. By \cite[(18.8.8) (2)]{AGro2}, it induces a local homomorphism $u \colon \mathcal{O}_{S,\bar{y}} \to \mathcal{O}_{S',\bar{y}'}$ which makes a commutative diagram
\[\xymatrix{\mathcal{O}_{S,y} \ar[r]^-{q^{\#}_y} \ar[d] & \mathcal{O}_{S',y'} \ar[d] \\ \mathcal{O}_{S,\bar{y}} \ar[r]_{u} & \mathcal{O}_{S',\bar{y}'}}\]
Let $v \colon \Spec \mathcal{O}_{S',\bar{y}'} \to \Spec \mathcal{O}_{S,\bar{y}}$ be the morphism induced by $u$. Then
\begin{enumerate}
\item The diagram
\[\xymatrix{\Spec \mathcal{O}_{S',\bar{y}'} \ar[r]^-{v} \ar[d] & \Spec \mathcal{O}_{S,\bar{y}} \ar[d] \\ S' \ar[r]_-{q} & S}\]
is commutative.
\item The square
\[\xymatrix{D_{\bar{y}'} \ar[r]^-{\iDe \times v} \ar[d] & D_{\bar{y}} \ar[d] \\ \Spec \mathcal{O}_{S',\bar{y}'} \ar[r]_--{v} & \Spec \mathcal{O}_{S,\bar{y}}}\]
is Cartesian.
\item If $X \times_S \Spec \mathcal{O}_{S,\bar{y}} \to \Spec \mathcal{O}_{S,\bar{y}}$ satisfies the condition \rmB{$\ast$} in \S\ref{Se:global:2}, so does $X' \times_{S'} \Spec \mathcal{O}_{S',\bar{y}'} \to \Spec \mathcal{O}_{S',\bar{y}'}$.
\enumS
\end{enumerate}
We assume that $X \times_S \Spec \mathcal{O}_{S,\bar{y}} \to \Spec \mathcal{O}_{S,\bar{y}}$ satisfies the condition \rmB{$\ast$} in \S\ref{Se:global:2}.
\begin{enumerate}
\enumL
\item $v$ induces a canonical bijection
\[\varphi \colon \cpC(y) \xrightarrow{\sim} \cpC(y') \,, \qquad C \mapsto (\iDe \times v)^{-1}(C) \,.\]
\item We have a commutative diagram
\[\xymatrix{\cpC(y) \ar[r]^-{\varphi}_-{\sim} \ar[d]_{\omega_{\bar{y}}} & \cpC(y') \ar[d]^{\omega_{\bar{y}'}} \\ \mathcal{O}_{S,\bar{y}}/\mathcal{O}_{S,\bar{y}}^{\ast} \ar[r]_-{\bar{u}} & \mathcal{O}_{S',\bar{y}'}/\mathcal{O}_{S',\bar{y}'}^{\ast}}\]
\item There is a canonical closed immersion
\[E_{\bar{y}'} \hookrightarrow E_{\bar{y}} \times_{\Spec \mathcal{O}_{S,\bar{y}}} \Spec \mathcal{O}_{S',\bar{y}'} \,.\]
\item If $q$ is flat at $y'$, then above morphism is an isomorphism.
\item $\mathscr{L}_{\bar{y}'}$ is isomorphic the inverse image of $\mathscr{L}_{\bar{y}}$ under the morphism $E_{\bar{y}'} \to E_{\bar{y}}$.
\end{enumerate}
\end{noteLe}

\begin{noteTh}[672P0]
If $f$ satisfies \nB \rmB{resp.~\nC}, so does $f'$.
\end{noteTh}

\begin{noteTh}[672P1]\
\begin{enumerate}
\item Assume that $f$ satisfies \nB. Let $\mathscr{P}$, $\mathscr{Q}$, $\theta$, $\vartheta$ and $\mathfrak{d}$ be the notations for $f$ defined in Remark \thmref{668P8}, and $\mathscr{P}'$, $\mathscr{Q}'$, $\theta'$, $\vartheta'$ and $\mathfrak{d}'$ the corresponding notations for $f'$. Then there are two canonical isomorphisms of sheaves of monoids $\lambda \colon p^{-1}\mathscr{P} \xrightarrow{\sim} \mathscr{P}'$ and $\mu \colon q^{-1}\mathscr{Q} \xrightarrow{\sim} \mathscr{Q}'$ which makes the following three diagrams commutative.
\begin{gather*}
\xymatrix{p^{-1}\mathscr{P} \ar[d]_{p^{-1}\theta} \ar[r]^-{\lambda}_-{\sim} & \mathscr{P}' \ar[d]^{\theta'} \\ p^{-1}\bigl(\mathcal{O}_X/\mathcal{O}_X^{\ast}\bigr) \ar[r] & \mathcal{O}_{X'}/\mathcal{O}_{X'}^{\ast}} \qquad
\xymatrix{q^{-1}\mathscr{Q} \ar[d]_{q^{-1}\vartheta} \ar[r]^-{\mu}_-{\sim} & \mathscr{Q}' \ar[d]^{\vartheta'} \\ q^{-1}\bigl(\mathcal{O}_S/\mathcal{O}_S^{\ast}\bigr) \ar[r] & \mathcal{O}_{S'}/\mathcal{O}_{S'}^{\ast}} \\
\xymatrix@C+5em{*!<-3em,0ex>{p^{-1}\bigl(f^{-1}\mathscr{Q}\bigr) = f'^{-1}\big(q^{-1}\mathscr{Q}\bigr)} \ar[r]^-{f'^{-1}\mu}_-{\sim} \ar[d]_{p^{-1}\mathfrak{d}} & \mathscr{Q}' \ar[d]^{\mathfrak{d}'} \\ p^{-1}\mathscr{P} \ar[r]^-{\lambda}_-{\sim} & \mathscr{P}'}
\end{gather*}
\item Assume that $f$ satisfies \nC. Let $(\mathscr{M}, \mathscr{N}, \sigma, \tau, \varphi)$ and $(\mathscr{M}', \mathscr{N}', \sigma', \tau', \varphi')$ be the semistable log structures for $f$ and $f'$ respectively. Then there exists two isomorphisms $\zeta \colon p^{\ast}\mathscr{M} \xrightarrow{\sim} \mathscr{M}'$ and $\eta \colon q^{\ast}\mathscr{N} \xrightarrow{\sim} \mathscr{N}'$ of log structures which make the following three diagrams commutative.
\begin{gather*}
\xymatrix@C+5em{*!<-3em,0ex>{p^{\ast}\bigl(f^{\ast}\mathscr{N}\bigr) = f'^{\ast}\big(q^{\ast}\mathscr{N}\bigr)} \ar[r]^-{f'^{\ast}\eta}_-{\sim} \ar[d]_{p^{\ast}\varphi} & \mathscr{N}' \ar[d]^{\varphi'} \\ p^{\ast}\mathscr{M} \ar[r]^-{\zeta}_-{\sim} & \mathscr{M}'} \\
\xymatrix{p^{-1}\mathscr{M} \ar[r] \ar[d]_{p^{-1}\sigma} & p^{\ast}\mathscr{M} \ar@{-->}[ld] \ar[r]^-{\zeta}_-{\sim} & \mathscr{M}' \ar[d]^{\sigma'} \\ p^{-1}\mathscr{P} \ar[rr]^-{\lambda}_-{\sim} & & \mathscr{P}'} \qquad
\xymatrix{q^{-1}\mathscr{N} \ar[r] \ar[d]_{q^{-1}\tau} & q^{\ast}\mathscr{N} \ar@{-->}[ld] \ar[r]^-{\eta}_-{\sim} & \mathscr{N}' \ar[d]^{\tau'} \\ q^{-1}\mathscr{Q} \ar[rr]^-{\mu}_-{\sim} & & \mathscr{Q}'}
\end{gather*}
Moreover the pair $(\zeta,\eta)$ is unique. Simply speaking, the semistable log structure of $f'$ may be viewed as the inverse image of that of $f$.
\end{enumerate}
\end{noteTh}

The following theorem shows that above isomorphisms are functorial.

\begin{noteTh}[672P8]
Let
\[\xymatrix{X_2 \ar[r]^{p_2} \ar[d]_{f_2} \CartS & X_1 \ar[r]^{p_1} \ar[d]|{f_1} \CartS & X_0 \ar[d]^{f_0} \\ S_2 \ar[r]_{q_2} & S_1 \ar[r]_{q_1} & S_0}\]
be a commutative diagram of locally noetherian schemes with both squares Cartesian. Put $p_0 \defeq p_1 \circ p_2$ and $q_0 \defeq q_1 \circ q_2$.
\begin{enumerate}
\item Assume that $f$ satisfies \nB. For each $i = 1,2,3$, let $\mathscr{P}_i$, $\mathscr{Q}_i$, $\theta_i$, $\vartheta_i$ and $\mathfrak{d}_i$ be the notations for $f_i$ defined in Remark \thmref{668P8}. Let
\begin{align*}
\lambda_1 & \colon p_1^{-1}\mathscr{P}_0 \xrightarrow{\sim} \mathscr{P}_1 \,, & \lambda_2 & \colon p_2^{-1}\mathscr{P}_1 \xrightarrow{\sim} \mathscr{P}_2 \,, & \lambda_0 & \colon p_0^{-1}\mathscr{P}_0 \xrightarrow{\sim} \mathscr{P}_2 \,, \\
\mu_1 & \colon q_1^{-1}\mathscr{Q}_0 \xrightarrow{\sim} \mathscr{Q}_1 \,, & \mu_2 & \colon q_2^{-1}\mathscr{Q}_1 \xrightarrow{\sim} \mathscr{Q}_2 \,, & \mu_0 & \colon q_0^{-1}\mathscr{Q}_0 \xrightarrow{\sim} \mathscr{Q}_2 \,,
\end{align*}
be the isomorphisms defined in Theorem \thmref{672P1} {\upshape(1)}. Then 
\[\lambda_2 \circ p_2^{-1}(\lambda_1) = \lambda_0 \qquad \mbox{and} \qquad \mu_2 \circ q_2^{-1}(\mu_1) = \mu_0 \,.\]
\item Assume that $f$ satisfies \nC. For each $i = 1,2,3$, let $(\mathscr{M}_i, \mathscr{N}_i, \sigma_i, \tau_i, \varphi_i)$ be the semistable log structure for $f_i$. Let
\begin{align*}
\zeta_1 & \colon p_1^{\ast}\mathscr{M}_0 \xrightarrow{\sim} \mathscr{M}_1 \,, & \zeta_2 & \colon p_2^{\ast}\mathscr{M}_1 \xrightarrow{\sim} \mathscr{M}_2 \,, & \zeta_0 & \colon p_0^{\ast}\mathscr{M}_0 \xrightarrow{\sim} \mathscr{M}_2 \,, \\
\eta_1 & \colon q_1^{\ast}\mathscr{N}_0 \xrightarrow{\sim} \mathscr{N}_1 \,, & \eta_2 & \colon q_2^{\ast}\mathscr{N}_1 \xrightarrow{\sim} \mathscr{N}_2 \,, & \eta_0 & \colon q_0^{\ast}\mathscr{N}_0 \xrightarrow{\sim} \mathscr{N}_2 \,,
\end{align*}
be the isomorphisms defined in Theorem \thmref{672P1} {\upshape(2)}. Then 
\[\zeta_2 \circ p_2^{\ast}(\zeta_1) = \zeta_0 \qquad \mbox{and} \qquad \eta_2 \circ q_2^{\ast}(\eta_1) = \eta_0 \,.\]
\end{enumerate}
\end{noteTh}

\begin{noteTh}[672P7]
Let $S_0$, $S_1$ and $X_0$ be locally noetherian schemes, $f_0 \colon X_0 \to S_0$ and $q \colon S_1 \to S_0$ be two morphisms. Put $S_2 \defeq S_1 \times_{S_0} S_1$ and $S_3 \defeq S_1 \times_{S_0} S_1 \times_{S_0} S_1$. For each $i=1,2,3$, let $X_i \defeq X_0 \times_{S_0} S_i$ and $f_i \colon X_i \to S_i$ the second projections. Assume that both $S_2$ and $S_3$ are locally noetherian.
\[\xymatrix{X_3 \ar@<.4ex>[r] \ar[r] \ar@<-.4ex>[r] \ar[d]_{f_3} & X_2 \ar@<.3ex>[r] \ar@<-.3ex>[r] \ar[d]|{f_2} & X_1 \ar[r] \ar[d]|{f_1} & X_0 \ar[d]^{f_0} \\ S_3 \ar@<.4ex>[r] \ar[r] \ar@<-.4ex>[r] & S_2 \ar@<.3ex>[r] \ar@<-.3ex>[r] & S_1 \ar[r]_{q} & S_0}\]

\begin{enumerate}
\item Assume that $f_1$ satisfies \nB. Then $f_2$ and $f_3$ also satisfy \nB. For each $i = 1,2,3$, let $\mathscr{P}_i$, $\mathscr{Q}_i$, $\theta_i$, $\vartheta_i$ and $\mathfrak{d}_i$ be the notations for $f_i$ defined in Remark \thmref{668P8}. For each $i=1,2$, let $\lambda_i \colon \prM{i}^{-1}\mathscr{P}_1 \xrightarrow{\sim} \mathscr{P}_2$ and $\mu_i \colon \prM{i}^{-1}\mathscr{Q}_1 \xrightarrow{\sim} \mathscr{Q}_2$ be the isomorphisms corresponding to the \mDash{i}{th} projections defined in Theorem \thmref{672P8} {\upshape(1)}. Put
\begin{align*}
\lambda & \defeq \lambda_2^{-1} \circ \lambda_1 \colon \prM{1}^{-1}\mathscr{P}_1 \to \prM{2}^{-1}\mathscr{P}_1 \,, \\
\mu & \defeq \mu_2^{-1} \circ \mu_1 \colon \prM{1}^{-1}\mathscr{Q}_1 \to \prM{2}^{-1}\mathscr{Q}_1 \,.
\end{align*}
Then
\[\prM{23}^{-1}(\lambda) \circ \prM{12}^{-1}(\lambda) = \prM{13}^{-1}(\lambda) \qquad \mbox{and} \qquad \prM{23}^{-1}(\mu) \circ \prM{12}^{-1}(\mu) = \prM{13}^{-1}(\mu) \,.\]

\item Assume that $f_1$ satisfies \nC. Then $f_2$ and $f_3$ also satisfy \nC. For each $i = 1,2,3$, let $(\mathscr{M}_i, \mathscr{N}_i, \sigma_i, \tau_i, \varphi_i)$ be the semistable log structure for $f_i$. For each $i=1,2$, let $\zeta_i \colon \prM{i}^{\ast}\mathscr{M}_1 \xrightarrow{\sim} \mathscr{M}_2$ and $\eta_i \colon \prM{i}^{\ast}\mathscr{N}_1 \xrightarrow{\sim} \mathscr{N}_2$ be the isomorphisms corresponding to the \mDash{i}{th} projections defined in Theorem \thmref{672P8} {\upshape(2)}. Put
\begin{align*}
\zeta & \defeq \zeta_2^{-1} \circ \zeta_1 \colon \prM{1}^{\ast}\mathscr{M}_1 \xrightarrow{\sim} \prM{2}^{\ast}\mathscr{M}_1 \,, \\
\mu & \defeq \mu_2^{-1} \circ \mu_1 \colon \prM{1}^{\ast}\mathscr{N}_1 \xrightarrow{\sim} \prM{2}^{\ast}\mathscr{N}_1 \,.
\end{align*}
Then
\[\prM{23}^{\ast}(\zeta) \circ \prM{12}^{\ast}(\zeta) = \prM{13}^{\ast}(\zeta) \qquad \mbox{and} \qquad \prM{23}^{\ast}(\eta) \circ \prM{12}^{\ast}(\eta) = \prM{13}^{\ast}(\eta) \,.\]
\end{enumerate}
\end{noteTh}

\begin{proof}
For each $i=1,2,3$, let $\lambda'_i \colon \prM{i}^{-1}\mathscr{P}_1 \xrightarrow{\sim} \mathscr{P}_3$ be the isomorphism corresponding to the \mDash{i}{th} projection. For each $1 \leqslant i < j \leqslant 3$, let $\lambda''_{ij} \colon \prM{ij}^{-1}\mathscr{P}_2 \xrightarrow{\sim} \mathscr{P}_3$ be the isomorphism corresponding to $\prM{ij} \colon X_3 \to X_2$. By Theorem \thmref{672P8}, we have
\[\lambda''_{ij} \circ \prM{ij}^{-1}(\lambda_1) = \lambda'_i \qquad \mbox{and} \qquad \lambda''_{ij} \circ \prM{ij}^{-1}(\lambda_2) = \lambda'_j \,.\]
Thus
\[\prM{ij}^{-1}(\lambda) = \prM{ij}^{-1}(\lambda_2^{-1}) \circ \prM{ij}^{-1}(\lambda_1) = (\lambda'_j)^{-1} \circ \lambda'_j \,.\]
The other three equations are similar.
\end{proof}

\subsection{Properties under inverse limit}

\begin{noteTh}[672P5]
Let $S_0$ be a noetherian scheme, $(S_{\lambda},s_{\mu\lambda})_{\lambda,\mu \in \Lambda}$ an inverse system of noetherian affine \mDash{S_0}{schemes} such that all $s_{\mu\lambda}$ are affine morphisms. Let $(S,s_{\lambda})$ be its inverse limit. Assume that $S$ is also noetherian.
\begin{enumerate}
\item Let $f_0 \colon X_0 \to S_0$ be a morphism of finite type. Put $X \defeq X_0 \times_{S_0} S$ and $f \defeq (f_0)_S \colon X \to S$. For each $\lambda \in \Lambda$, put $X_{\lambda} \defeq X_0 \times_{S_0} S_{\lambda}$ and $f_{\lambda} \defeq (f_0)_{S_{\lambda}} \colon X_{\lambda} \to S_{\lambda}$. Then $f$ satisfies \nA \rmB{resp.~\nB or \nC} if and only if there exists an index $\lambda_0 \in \Lambda$ such that for any $\lambda \geqslant \lambda_0$, $f_{\lambda}$ satisfies \nA \rmB{resp.~\nB or \nC}. 
\item Let $f \colon X \to S$ be a morphism which satisfies \nA \rmB{resp.~\nB or \nC}. Then there exists an index $\lambda_0 \in \Lambda$ and a morphism $f_{\lambda_0} \colon X_{\lambda_0} \to S_{\lambda_0}$ which satisfies \nA \rmB{resp.~\nB or \nC} such that $X$ is \mDash{S}{isomorphic} to $X_{\lambda_0} \times_{S_{\lambda_0}} S$.
\end{enumerate}
\end{noteTh}

\begin{proof}
See \cite[\S8]{AGro2}. Note that every \lC (resp.~\rlC) of $f$ can be descended to some index $\lambda_0 \in \Lambda$ and $X$ can be covered by a finite number of \lC (resp.~\rlC) of $f$.
\end{proof}

\begin{noteCo}[672P6]
Let $f \colon X \to S$ be a morphism of finite type of locally noetherian schemes. Then $f$ satisfies \nA \rmB{resp.~\nB or \nC} if and only if for every point $y \in S$, $X \times_S \Spec \mathcal{O}_{S,\bar{y}} \to \Spec \mathcal{O}_{S,\bar{y}}$ satisfies \nA \rmB{resp.~\nB or \nC}.
\end{noteCo}

\subsection{Properties under flat descent}

\begin{noteLe}
Let $X' \to X$ be a faithfully flat morphism locally of finite presentation of schemes, $U \to X$ an \'{e}tale morphism of schemes, $\mathscr{M}$ a fine saturated log structure on $X$. Put $X'' \defeq X' \times_X X'$, and let $\mathscr{M}'$ and $\mathscr{M}''$ be the pull-back of $\mathscr{M}$ on $X'$ and $X''$ respectively. Then
\[\xymatrix{\mathscr{M}(U) \ar[r] & \mathscr{M}'(U \times_X X') \ar@<.4ex>[r] \ar@<-.4ex>[r] & \mathscr{M}''(U \times_X X'')}\]
is exact.
\end{noteLe}

\begin{proof}
\cite[Lemma 1.1.3]{AOMCO}.
\end{proof}

\begin{noteLe}[672P2]
Let $p \colon X' \to X$ be a faithfully flat morphism locally of finite presentation of schemes. Put $X'' \defeq X' \times_X X'$ and $X''' \defeq X' \times_X X' \times_X X'$. Let $\mathscr{M}$ be a fine saturated log structure on $X'$ and $\phi \colon \prM{1}^{\ast}\mathscr{M} \xrightarrow{\sim} \prM{2}^{\ast}\mathscr{M}$ an isomorphism of log structures on $X''$ such that on $X'''$ we have
\[\prM{13}^{\ast}(\phi) = \prM{23}^{\ast}(\phi) \circ \prM{12}^{\ast}(\phi) \,.\]
Then there exists a unique (up to isomorphism) pair $(\mathscr{N},s)$ on $X$, where $\mathscr{N}$ is a fine saturated log structure on $X$ and $s \colon p^{\ast}\mathscr{N} \xrightarrow{\sim} \mathscr{M}$ is an isomorphism of log structures on $X'$, such that the following diagram
\[\xymatrix{\prM{1}^{\ast}p^{\ast}\mathscr{N} \ar[r]^-{\prM{1}^{\ast}(s)} \ar@{=}[d] & \prM{1}^{\ast}\mathscr{M} \ar[d]^{\phi} \\ \prM{2}^{\ast}p^{\ast}\mathscr{N} \ar[r]^-{\prM{2}^{\ast}(s)} & \prM{2}^{\ast}\mathscr{M}}\]
is commutative.
\end{noteLe}

\begin{proof}
\cite[Theorem 1.1.5]{AOMCO}.
\end{proof}

\begin{noteLe}[672P3]
Let $f \colon (X,\mathscr{M}) \to (Y,\mathscr{N})$ be a morphism of fine saturated log schemes, $\mathscr{C}$ the cokernel of $f^{\ast}\!\mathscr{N} \to \mathscr{M}$, $\bar{x}$ a geometric point on $X$, $\bar{y} \defeq f(\bar{x})$, $Q \to \mathscr{N}$ a chart with $Q$ fine saturated. Assume that
\begin{enumerate}
\item $Q \to \overline{\mathscr{N}}_{\bar{y}}$ is isomorphic;
\item $\overline{\mathscr{N}}_{\bar{y}} \to \overline{\mathscr{M}}_{\bar{x}}$ is injective;
\item the torsion part of $\Groth{\mathscr{C}_{\bar{x}}}$ is a finite group of order invertible in $\kappa(\bar{x})$.
\end{enumerate}
Then there exists a chart 
\[(P \to \mathscr{M}|_U, Q \to \mathscr{N}|_V, Q \to P)\]
of $f$ at $\bar{x}$ such that $P \to \overline{\mathscr{M}}_{\bar{x}}$ is isomorphic.
\end{noteLe}

\begin{proof}
\cite[Theorem 2.13]{AOgus1}.
\end{proof}

\begin{noteLe}[128A3]
If $S$ is a scheme, $x$ is a point on $S$, and $\mathscr{F}$ is an \mDash{\mathcal{O}_S}{module}, we define $\mathscr{F}(x) \defeq \mathscr{F}_{\bar{x}} \otimes \kappa(\bar{x})$.

Let $X^{\dag} = (X,\mathscr{M})$ and $Y^{\dag} = (Y,\mathscr{N})$ be two fine log schemes, $f \colon X^{\dag} \to Y^{\dag}$ a morphism of log schemes, $x \in X$ a point, $\mathscr{C}$ the cokernel of the morphism $f^{\ast}\mathscr{N} \to \mathscr{M}$. Then there  is a commutative diagram
\[\xymatrix@C+2.5em{& \Groth{\mathscr{M}}_{\bar{x}} \ar[r] \ar[rd]^-{\pi} \ar[d]_{d\mathrm{log}} & \Groth{\mathscr{C}}_{\bar{x}} \ar[d]^{\bar{\pi}} \\ *!<-2em,0ex>{\varOmega_{X/Y}(x)} \ar[r] & \varOmega_{X^{\dag}/Y^{\dag}}(x) \ar[r]_-{\rho_{X^{\dag}/Y^{\dag},\bar{x}}} & *!<-.2em,0ex>{\kappa(\bar{x}) \otimes_{\mathbb{Z}} \Groth{\mathscr{C}}_{\bar{x}} \,,} }\]
where $\bar{\pi}(m) = 1 \otimes m$ for each $m \in \Groth{\mathscr{C}}_{\bar{x}}$, and the bottom row is an exact sequence of linear space over $\kappa(\bar{x})$.

The homomorphism $\rho_{X^{\dag}/Y^{\dag},\bar{x}}$ is sometimes called the \emph{Poincar\'{e} residue mapping} at $x$.
\end{noteLe}

\begin{proof}
See \cite[Proposition 2.22]{AOgus1}.
\end{proof}

\begin{noteLe}[181P5]
Let $f \colon (X,\mathscr{M}) \to (Y,\mathscr{N})$ be a log smooth morphism of fine saturated log schemes, $x \in X$, $y \defeq f(x)$, and let
\[(P \to \mathscr{M}, Q \to \mathscr{N}, Q \xrightarrow{\iota} P)\]
be a chart of $f$. Assume that $\iota \colon Q \to P$ is injective, $P \xrightarrow{\sim} \overline{\mathscr{M}}_{\bar{x}}$ and $Q \xrightarrow{\sim} \overline{\mathscr{N}}_{\bar{y}}$ are isomorphisms. Put $Z \defeq Y \times_{\Spec \mathbb{Z}[Q]} \Spec \mathbb{Z}[P]$ and $g \colon X \to Z$ the induced morphism. Then $g$ is smooth at $\bar{x}$.
\end{noteLe}

\begin{proof}
In this proof, for a log scheme $(X,\mathscr{M})$, we use simply $X^{\dag}$ to denote it. First we have
\[\varOmega_{Z^{\dag}/Y^{\dag}} = \mathcal{O}_Y \otimes_{\mathbb{Z}} (\Groth{P}/\Groth{Q}) \,.\]
Hence
\[\kappa(\bar{x}) \otimes (g^{\ast}\varOmega_{Z^{\dag}/Y^{\dag}}) \xrightarrow{\sim} \kappa(\bar{x}) \otimes_{\mathbb{Z}} (\Groth{P}/\Groth{Q}) \,.\]
Let $\mathscr{C}$ denote the cokernel of the morphism $f^{\ast}\mathscr{N} \to \mathscr{M}$. Then the composite
\[\kappa(\bar{x}) \otimes (g^{\ast}\varOmega_{Z^{\dag}/Y^{\dag}}) \to \kappa(\bar{x}) \otimes \varOmega_{X^{\dag}/Y^{\dag}} \xrightarrow{\rho_{X^{\dag}/Y^{\dag},\bar{x}}} \kappa(\bar{x}) \otimes_{\mathbb{Z}} \Groth{\mathscr{C}}_{\bar{x}}\]
is isomorphic, where $\rho_{X^{\dag}/Y^{\dag},\bar{x}}$ is defined in Lemma \thmref{128A3}. Hence
\[\kappa(\bar{x}) \otimes (g^{\ast}\varOmega_{Z^{\dag}/Y^{\dag}}) \to \kappa(\bar{x}) \otimes \varOmega_{X^{\dag}/Y^{\dag}}\]
is injective. Hence $g^{\ast}\varOmega_{Z^{\dag}/Y^{\dag}} \to \varOmega_{X^{\dag}/Y^{\dag}}$ has a left inverse in some open neighborhood $U$ of $x$. So $g|_U \colon U \to Z$ is smooth.
\end{proof}

\begin{noteTh}[181P6]
Let
\[\xymatrix{X' \ar[r]^{f'} \ar[d]_{p} \CartS & S' \ar[d]^{q} \\ X \ar[r]_{f} & S}\]
be a Cartesian square of locally noetherian schemes such that $q$ is faithfully flat and locally of finite type, and $f'$ satisfies \nC. Then $f$ satisfies \nC too.
\end{noteTh}

\begin{proof}
Let $\mathscr{P}'$, $\mathscr{Q}'$, $\theta'$, $\vartheta'$ and $\mathfrak{d}'$ be the notations for $f'$ defined in Remark \thmref{668P8}; and let $(\mathscr{M}', \mathscr{N}', \sigma', \tau', \varphi')$ be the semistable log structure of $f'$. By Theorem \thmref{672P7} and Lemma \thmref{672P2}, there exist fine saturated log structures $\mathscr{M}$ on $X$ and $\mathscr{N}$ on $S$, a morphism of log structures $\varphi \colon f^{\ast}\mathscr{N} \to \mathscr{M}$, two isomorphisms of log structures $p^{\ast}\mathscr{M} \xrightarrow{\sim} \mathscr{M}'$ and $q^{\ast}\mathscr{N} \xrightarrow{\sim} \mathscr{N}'$ which makes a commutative diagram.
\[\xymatrix@C+5em{*!<-2.5em,0ex>{f'^{\ast}q^{\ast}\mathscr{N} = p^{\ast}f^{\ast}\mathscr{N}} \ar[r]^-{p^{\ast}(\varphi)} \ar[d]_{\cong} & p^{\ast}\mathscr{M} \ar[d]^{\cong} \\ f'^{\ast}\mathscr{N}' \ar[r]_-{\varphi'} & \mathscr{M}'}\]
Let $x$ be a point on $X$ and $y \defeq f(x)$. Let $x' \in p^{-1}(x)$ and put $y' \defeq f'(x')$. Then $q(y') = y$. We have
\[\overline{\mathscr{M}}_{\bar{x}} \cong \overline{\mathscr{M}'}_{\bar{x}'} \cong \mathbb{N}^m \qquad \mbox{and} \qquad \overline{\mathscr{N}}_{\bar{y}} \cong \overline{\mathscr{N}'}_{\bar{y}'} \cong \mathbb{N}^n\]
for some $m,n \in \mathbb{N}$. Furthermore the homomorphism
\[d \defeq \bar{\varphi}_{\bar{x}} \colon \overline{\mathscr{N}}_{\bar{y}} \to \overline{\mathscr{M}}_{\bar{x}}\]
is defined as $d(\varepsilon_i) = \sum\limits^{s_i}_{j=1}\eta_{ij}$ for $i \in [1,r]$ and $d(\varepsilon_i) = \eta_i$ for $i \in [r+1,n]$, where $\{\Sequ{n}{\varepsilon}\}$ is a basis of $\mathbb{N}^n$,
\[\{\eta_{11},\ldots,\eta_{1s_1},\ldots,\eta_{r1},\ldots,\eta_{rs_r},\eta_{r+1},\ldots,\eta_n\}\]
is a basis of $\mathbb{N}^m$, and $m = \sum\limits^r_{i=1}s_i + n - r$. By Lemma \thmref{672P3}, there exists a chart
\[(\mathbb{N}^m \to \mathscr{M}|_U, \mathbb{N}^n \to \mathscr{N}|_V, \mathbb{N}^n \xrightarrow{d} \mathbb{N}^m)\]
of $f$ at $\bar{x}$ such that $\mathbb{N}^m \xrightarrow{\sim} \overline{\mathscr{M}}_{\bar{x}}$ and $\mathbb{N}^n \xrightarrow{\sim} \overline{\mathscr{N}}_{\bar{y}}$ are isomorphic. Put $U' \defeq U \times_X X'$ and $V' \defeq V \times_S S'$. By Lemma \thmref{181P5},
\[U' \to V' \times_{\Spec \mathbb{Z}[\mathbb{N}^n]} \Spec \mathbb{Z}[\mathbb{N}^m]\]
is smooth at $\bar{x}'$. By \cite[(17.7.1)]{AGro2},
\[U \to V \times_{\Spec \mathbb{Z}[\mathbb{N}^n]} \Spec \mathbb{Z}[\mathbb{N}^m]\]
is smooth at $\bar{x}$. Therefore we may contract $U/V$ suitable to obtain a \lC at $x$. Furthermore it is easy to verify that $(\mathscr{M},\mathscr{N},\varphi)$ is just the semistable log structure of $f$.
\end{proof}

\section{Semistable Curves}
\begin{noteDe}[182P0]
Let $k$ be a separably closed field. A \emph{semistable curve} over $k$ is a connected proper \mDash{1}{equidimensional} \mDash{k}{scheme} $X$ such that for any closed point $x \in X$, either $X$ is smooth at $x$ over $k$, or $\widehat{\mathcal{O}}_{X,x}$ is \mDash{k}{isomorphic} to $k[[T_1,T_2]]\big/(T_1T_2)$.
\end{noteDe}

\begin{noteLe}[182P9]
Let $k$ be a separably closed field and $X$ a semistable curve over $k$. Then
\begin{enumerate}
\item $X$ is reduced.
\item $X$ has only a finite number of singular points and all singular points are \mDash{k}{rational}.
\end{enumerate}
\end{noteLe}

\begin{noteDe}[182P8]
Let $S$ be a scheme. A \emph{semistable curve} over $S$ is an \mDash{S}{scheme} $f \colon X \to S$ such that $f$ is proper, faithfully flat, of finite presentation, and every geometric fiber is a semistable curve in the sense of Definition \thmref{182P0}.
\end{noteDe}

\begin{noteRe}
The notation of \emph{semistable curve} here is slightly weaker than the notation of \emph{stable curve} in \cite{PDDM}.
\end{noteRe}

\begin{noteLe}[182P2]
Let $S$ be a locally noetherian scheme, $f \colon X \to S$ a morphism which satisfies \nA in Definition \thmref{672P4} and is of relative dimension $1$. Then $f$ also satisfies \nC in Definition \thmref{672P4}.
\end{noteLe}

\begin{proof}
By Corollary \thmref{672P6}, we may assume that $S$ is the spectrum of a strictly Henselian ring $R$. Let $\mathfrak{m}$ be the maximal ideal of $R$ and $y_1$ the closed point of $S$. If $f$ is smooth, the question is trivial. Assume that $f$ is not smooth. Then $D(f) \neq \emptyset$. Let $\Sequ{n}{D}$ be the connected components of $D(f)$.  For each $i \in [1,n]$, select an element $a_i \in \mathfrak{m}$ such that $\bar{a}_i = \omega_{\bar{y}_1}(D_i)$. Then $D_i = \Spec R/(a_i)$. Let $z_i$ be the closed point of $D_i$ and let $l_i \colon D_i \to X$ be the inclusion. Then there exists a \rlC at $l_i(z_i)$ of the form
\[(U_i,S; T_{i1},T_{i2}; \Sequ{n}{a})\]
satisfying that $U_i \times_X D_j = \emptyset$ for all $j \in [1,n] - \{i\}$. Let $\mathscr{M}_i$ be the log structure on $U_i$ associate to $\alpha'_i \colon \mathbb{N}^{n+1}_{U_i} \to \mathcal{O}_{U_i}$, where if \[\eta_{i1},\eta_{i2},\eta_1,\ldots,\Hat{\eta}_i,\ldots,\eta_n\]
is a basis of $\mathbb{N}^{n+1}$, then $\alpha'_i(\eta_{ij}) = T_{ij}$ for $j=1,2$, and $\alpha'_i(\eta_k) = a_k$ for $k \in [1,n]-\{i\}$. Let $\alpha_i \colon \mathbb{N}^{n+1}_{U_i} \to \mathscr{M}_i$ be the induced morphism. Let $\Sequ{n}{\varepsilon}$ be a basis of $\mathbb{N}^n$. We define three homomorphisms of monoids
\[\partial_i, \partial_{i1},\partial_{i2} \colon \mathbb{N}^n \to \mathbb{N}^{n+1}\]
as follows: for $k \in [1,n] -\{i\}$,
\[\partial_i(\varepsilon_k) = \partial_{i1}(\varepsilon_k) = \partial_{i2}(\varepsilon_k) = \eta_k \,,\]
and
\[\partial_i(\varepsilon_i) = \eta_{i1}+\eta_{i2} \,, \qquad \partial_{i1}(\varepsilon_i) = \eta_{i1} \,, \qquad \partial_{i2}(\varepsilon_i) = \eta_{i2} \,.\]
Then
\[\rho_i \defeq \alpha_i \circ \partial_i \colon \mathbb{N}^n_{U_i} \to \mathscr{M}_i\]
is a lifting of $\gamma|_{U_i} \colon \mathbb{N}_{U_i}^n \to \mathscr{P}|_{U_i}$, where Let $\gamma \colon \mathbb{N}_X^n \to \mathscr{P}$ be the notation defined in \S\ref{Se:loc}. Put $U_0 \defeq X - D(f)$ and $\mathscr{M}_0$ the log structure on $U_0$ induced by
\[\mathbb{N}_{U_0}^n \to \mathcal{O}_{U_0} \,, \qquad \varepsilon_i \mapsto a_i \,.\]
Let $\rho_0 \colon \mathbb{N}_{U_0}^n \to \mathscr{M}_0$ be the induced morphism. Then $\rho_0$ is a lifting of $\gamma|_{U_0} \colon \mathbb{N}_{U_0}^n \to \mathscr{P}|_{U_0}$.

For any $i \in [1,n]$ and any point $x$ on $U_0 \times_X U_i$, there exists an \'{e}tale neighborhood $W$ of $\bar{x}$ such that $T_{i1}|_W$ or $T_{i2}|_W$ is invertible. Without lose of generality, we may assume that $T_{i2}|_W$ is invertible. Then $T_{i1}|_W = (T_{i2}|_W)^{-1} \cdot a_i$ and
\[\alpha_i|_W \circ \partial_{i1} \colon \mathbb{N}^n \to \mathscr{M}_i|_W\]
is a chart of $\mathscr{M}_i|_W$. So there exists an isomorphism $\phi \colon \mathscr{M}_i|_W \xrightarrow{\sim} \mathscr{M}_0|_W$ of log structures such that
\[\phi \circ \alpha_i|_W \circ \partial_{i1} = (\Sequ{n}{u}) \cdot \rho_0|_W \,,\]
where $u_i = (T_{i2}|_W)^{-1}$ and $u_k = 1$ for $k \in [1,n] - \{i\}$. Thus $\phi \circ \rho_i|_W = \rho_0|_W$.

For any pair of integers $1 \leqslant i < j \leqslant n$ and any point $x$ on $U_i \times_X U_j$, there exists an \'{e}tale neighborhood $W$ of $\bar{x}$ such that $T_{i1}|_W$ or $T_{i2}|_W$ is invertible, and $T_{j1}|_W$ or $T_{j2}|_W$ is invertible. Without lose of generality, we may assume that $T_{i2}|_W$ and $T_{j2}|_W$ are invertible. Then for $s = i,j$, $T_{s1}|_W = (T_{s2}|_W)^{-1} \cdot a_s$ and
\[\alpha_s|_W \circ \partial_{s1} \colon \mathbb{N}^n \to \mathscr{M}_s|_W\]
is a chart of $\mathscr{M}_s|_W$. So there exists an isomorphism $\phi \colon \mathscr{M}_i|_W \xrightarrow{\sim} \mathscr{M}_j|_W$ of log structures such that
\[\phi \circ \alpha_i|_W \circ \partial_{i1} = (\Sequ{n}{u}) \cdot (\alpha_j|_W \circ \partial_{j1}) \,,\]
where $u_i = (T_{i2}|_W)^{-1}$, $u_j = T_{j2}|_W$, and $u_k = 1$ for $k \in [1,n] - \{i,j\}$. Thus $\phi \circ \rho_i|_W = \rho_j|_W$.

Now we translate above analysis into the language in \S\ref{Se:loc}. We obtains an \'{e}tale covering $U \to X$, an object $\mathscr{M}$ in $\mathfrak{X}(U)$, and a morphism $\rho \colon \mathbb{N}^n_U \to \mathscr{M}$ which is a lifting of $\gamma|_U \colon \mathbb{N}^n_U \to \mathscr{P}|_U$, an \'{e}tale covering $V \to U \times_X U$, and an isomorphism $\phi \colon p_{10}^{\ast}(\mathscr{M}) \xrightarrow{\sim} p_{11}^{\ast}(\mathscr{M})$ of log structures on $V$, such that $\phi \circ p_{10}^{\ast}(\rho) = p_{11}^{\ast}(\rho)$. So $\mathscr{L}_{\bar{y}_1}$ is trivial.
\end{proof}

\begin{noteLe}[182P1]
Let $R, A$ be two complete noetherian local rings with maximal ideals $\mathfrak{m}, \mathfrak{M}$ and residue fields $k, K$ respectively, $R \to A$ a flat local homomorphism, $R[[T_1,T_2]]$ a ring of power series over $R$ with variables $T_1$ and $T_2$. Assume that exists two elements $x_{11},x_{12} \in \mathfrak{M}$ such that the homomorphism of \mDash{k}{algebras}
\[k[[T_1,T_2]]\big/(T_1T_2) \xrightarrow{\sim} A/\mathfrak{m}A \,, \qquad T_i \mapsto \bar{x}_{1i} \ (i = 1,2)\]
is an isomorphism. Then there exists an element $a \in \mathfrak{m}$, two elements $x_1,x_2 \in \mathfrak{M}$ such that the homomorphism of \mDash{R}{algebras}
\[R[[T_1,T_2]]\big/(T_1T_2-a) \xrightarrow{\sim} A \,, \qquad T_i \mapsto x_i \ (i = 1,2)\]
is an isomorphism.
\end{noteLe}

\begin{proof}
Put $a_1 \defeq 0$ and $P_1 \defeq k[[T_1,T_2]]\big/(T_1T_2)$. For each $n \in \mathbb{N}$, put $R_n \defeq R/\mathfrak{m}^n$ and $A_n \defeq A/\mathfrak{m}^nA$. Let $\psi_1 \colon P_1 \xrightarrow{\sim} A_1$ be the isomorphism defined in the lemma. Assume that we have found $a_n \in R$ and $x_{n1},x_{n2} \in A$ such that
\[\psi_n \colon P_n \defeq (R/\mathfrak{m}^n)[[T_1,T_2]]\big/(T_1T_2-\bar{a}_n) \xrightarrow{\sim} A_n \,, \qquad T_i \mapsto \bar{x}_{ni} \ (i = 1,2)\]
is an isomorphism. Then
\[z \defeq x_{n1}x_{n2} - a_n \in \mathfrak{m}^nA \,.\]
Put $z = \sum\limits^n_{j=1}b_jz_j$, where $b_j \in \mathfrak{m}^n$ and $z_j \in A$. Obviously we have $A = R + \mathfrak{M}$ and $\mathfrak{M} = \mathfrak{m}A + x_{n1}A + x_{n2}A$. So for each $j \in [1,n]$, we may write $z_j$ as
\[z_j = c_j + d_ju_j + x_{n1}v_j + x_{n2}w_j \,,\]
where $c_j \in R$, $d_j \in \mathfrak{m}$, $u_j, v_j, w_j \in A$. Put
\begin{align*}
x_{n+1,1} & \defeq x_{n1} - \sumT^n_{j=1}b_jw_j \,, \\
x_{n+1,2} & \defeq x_{n2} - \sumT^n_{j=1}b_jv_j \,, \\
a_{n+1} & \defeq a_n + \sumT^n_{j=1}b_jc_j \,.
\end{align*}
Then we have
\[x_{n+1,1}x_{n+1,2} - a_{n+1} = \sumT^n_{j=1}b_jd_ju_j + \Bigl(\sumT^n_{j=1}b_jv_j \Bigr) \Bigl(\sumT^n_{j=1}b_jw_j \Bigr) \in \mathfrak{m}^{n+1}A \,.\]
Put
\[P_{n+1} \defeq (R/\mathfrak{m}^{n+1})[[T_1,T_2]]\big/(T_1T_2 - \bar{a}_{n+1})\]
and let $\psi_{n+1} \colon P_{n+1} \to A_{n+1}$ be the homomorphism of \mDash{R}{algebras} defined by $\psi_{n+1}(T_i) = \bar{x}_{n+1,i}$ for $i=1,2$. Obviously $\psi_{n+1}$ is surjective. We shall prove that $\psi_{n+1}$ is injective. Assume that $\mathfrak{I} \defeq \Ker(\psi_{n+1}) \neq 0$. If $\mathfrak{m} = \mathfrak{m}^{n+1}$, then $\psi_{n+1} = \psi_1$ is an isomorphism. So we may assume that $\mathfrak{m}^{n+1} \neq \mathfrak{m}$. Since the diagram
\[\xymatrix{P_{n+1} \ar[r]^-{\psi_{n+1}} \ar[d] & A_{n+1} \ar[d] \\ P_n \ar[r]^-{\psi_n} & A_n}\]
is commutative, $\mathfrak{I}$ is contained in $(\mathfrak{m}^n/\mathfrak{m}^{n+1}) \cdot P_{n+1}$. Since both $P_{n+1}$ and $A_{n+1}$ are flat over $R_{n+1}$, so is $\mathfrak{I}$. As $\mathfrak{I}$ is a finitely generated \mDash{P_{n+1}}{module}, by \cite[Ch. 2, COROLLARY of (4.A)]{Mats1}, $\mathfrak{I}$ is faithfully flat over $R_{n+1}$. And because $\mathfrak{m}/\mathfrak{m}^{n+1} \neq 0$,
\[(\mathfrak{m}/\mathfrak{m}^{n+1}) \cdot \mathfrak{I} \cong (\mathfrak{m}/\mathfrak{m}^{n+1}) \otimes_R \mathfrak{I} \neq 0\]
by \cite[Ch. 2, (4.A)]{Mats1}. But it contradicts that
\[(\mathfrak{m}/\mathfrak{m}^{n+1}) \cdot \mathfrak{I} \subseteq (\mathfrak{m}/\mathfrak{m}^{n+1}) \cdot (\mathfrak{m}^n/\mathfrak{m}^{n+1}) \cdot P_{n+1} = 0 \,.\]
Thus $\psi_{n+1}$ is an isomorphism.

Clearly $\{a_n\}$ is a Cauchy sequence in $R$, and $\{x_{n1}\}$ and $\{x_{n2}\}$ are Cauchy sequences in $A$. Since $R$ and $A$ are complete, we may let $a \defeq \lim a_n$ and $x_i \defeq \lim x_{ni}$ for $i=1,2$.
\end{proof}

\begin{noteLe}[182P3]
Let $R$ be a ring, $A$ and $B$ two \mDash{R}{algebras}, $\mathfrak{a}$, $I$ and $J$ be ideals of $R$, $A$ and $B$ respectively such that $\mathfrak{a}A \subseteq I$ and $\mathfrak{a}B \subseteq J$. Let $C$ and $D$ denote the topological \mDash{R}{algebras} $A \otimes_R B$ and $\widehat{A} \otimes_{\widehat{R}} \widehat{B}$ equipped with \mDash{(IC+JC)}{adic} and \mDash{(ID + JD)}{adic} topologies. Then $\widehat{C} \cong \widehat{D}$.
\end{noteLe}

\begin{proof}
Note that for all $n \in \mathbb{N}$,
\[I^nC + J^nC \subseteq (IC+JC)^n\]
and
\[(IC+JC)^{2n} \subseteq I^nC + J^nC \,.\]
Thus
\[\widehat{C} \cong \varprojlim C/(I^nC + J^nC) \cong \varprojlim (A/I^n) \otimes_{R/\mathfrak{a}^n} (B/J^n) \,.\]
Similarly we have
\[\widehat{D} \cong \varprojlim (\widehat{A}/I^n\widehat{A}) \otimes_{\widehat{R}/\mathfrak{a}^n\widehat{R}} (\widehat{B}/J^n\widehat{B}) \cong \varprojlim (A/I^n) \otimes_{R/\mathfrak{a}^n} (B/J^n) \,.\]
Hence $\widehat{C} \cong \widehat{D}$.
\end{proof}

\begin{noteLe}[149T5]
Let $S$ be a scheme of finite type over a field or an excellent dedekind domain, $X_1$ and $X_2$ two \mDash{S}{schemes} of finite type, $x_1 \in X_1$ and $x_2 \in X_2$ two points which map onto the same point $s$ on $S$. Assume that $\widehat{\mathcal{O}}_{X_1,x_1}$ and $\widehat{\mathcal{O}}_{X_2,x_2}$ are \mDash{\mathcal{O}_{S,s}}{isomorphic}. Then there exists an \mDash{S}{scheme} $U$, a point $u \in U$, two \'{e}tale \mDash{S}{morphisms} $\varphi_1 \colon U \to X_1$ and $\varphi_2 \colon U \to X_2$, such that $\varphi_i(u) = x_i$ and $\kappa(x_i) \xrightarrow{\sim} \kappa(u)$ for $i=1,2$.
\end{noteLe}

\begin{proof}
See \cite[(2.6)]{MArt1}.
\end{proof}

\begin{noteRe}
Note that $\mathbb{Z}$ is an excellent dedekind domain. So to use this lemma, we usually apply the inverse limit of schemes to descend the base scheme to become of finite type over $\mathbb{Z}$.
\end{noteRe}

\begin{noteLe}
Let $A$ be a strictly Henselian noetherian local ring, $S \defeq \Spec A$, $s$ the closed point of $S$, $f \colon X \to S$ a faithfully flat, proper morphism such that $X_s$ is a semistable curve over $\kappa(s)$. Then $X$ is a semistable curve over $S$ and satisfies \nC in Definition \thmref{672P4}.
\end{noteLe}

\begin{proof}
Let $\mathfrak{m}$ be the maximal ideal of $A$ and $k \defeq A/\mathfrak{m}$ the residue field. If $X_s$ is smooth over $k$, then $X$ is smooth over $S$ and the lemma is valid. So we may assume that $X_s$ is not smooth over $k$. Let $\Sequ{n}{x}$ be all singular points of $X_s$. By Lemma \thmref{182P9}, all $x_i$ are \mDash{k}{rational}. So $x_i$ defines a closed immersion $\gamma_i \colon \Spec k \to X_s$. As $\widehat{\mathcal{O}}_{X_s,x_i} \cong k[[T_1,T_2]]\big/(T_1T_2)$, by Lemma \thmref{149T5}, there exists a \mDash{k}{scheme} $V_i$ of finite type, a point $y_i$ on $V_i$, two \'{e}tale \mDash{k}{morphisms} $p_i \colon V_i \to X_s$ and
\[q_i \colon V_i \to \Spec k[T_1,T_2]/(T_1T_2)\]
such that $p_i(y_i) = x_i$ and $q_i(y_i) = 0$, and $\kappa(y_i) = \kappa(x_i) = k$. So $y_i$ is a \mDash{k}{rational} point on $V_i$, and it defines a closed immersion $\delta_i \colon \Spec k \to V_i$. Let $U$ denote the set of points at which $X$ is smooth over $S$. Then $U$ is open in $X$ and
\[U_s = X_s - \{\Sequ{n}{x}\} \,.\]

By \cite[\S(8.8), \S(8.6), (8.10.5), (11.2.6), (17.7.8)]{AGro2}, there exists a finitely generated \mDash{\mathbb{Z}}{subalgebra} $R$ of $A$, a proper and faithfully flat \mDash{R}{scheme} $Y$ such that $X$ is \mDash{S}{isomorphic} to $Y \otimes_R A$, an open subscheme $U'$ of $Y$ such that $U'$ is smooth over $R$ and $U$ is the inverse image of $U'$ under the morphism $X \to Y$, an ideal $\mathfrak{a}$ of $R$ such that $\mathfrak{a}A = \mathfrak{m}$; and for each $i \in [1,n]$, an \mDash{(R/\mathfrak{a})}{scheme} $V'_i$ of finite type, two closed \mDash{(R/\mathfrak{a})}{immersions} $\gamma'_i \colon \Spec(R/\mathfrak{a}) \to Y \otimes_R (R/\mathfrak{a})$ and $\delta'_i \colon \Spec(R/\mathfrak{a}) \to V'_i$ such that
\begin{equation}\label{E:curve1}
\begin{split}
\bigl(U' \otimes_R (R/\mathfrak{a})\bigr) & \coprodT \underbrace{\Spec(R/\mathfrak{a}) \coprodT \cdots \coprodT \Spec(R/\mathfrak{a})}_n \\
& \xrightarrow{\gamma_0 \coprodT \gamma'_1 \coprodT \cdots \coprodT \gamma'_n} Y \otimes_R (R/\mathfrak{a})
\end{split}
\end{equation}
is surjective where $\gamma_0 \colon U' \otimes_R (R/\mathfrak{a}) \to Y \otimes_R (R/\mathfrak{a})$ is the inclusion, two \'{e}tale \mDash{(R/\mathfrak{a})}{morphisms} $p'_i \colon V'_i \to Y \otimes_R (R/\mathfrak{a})$ and
\[q'_i \colon V'_i \to \Spec (R/\mathfrak{a})[T_1,T_2]/(T_1T_2) \,,\]
finally a commutative diagram
\begin{equation}\label{E:curve2}
\vcenter{\xymatrix{& & \Spec k \ar[lld]_(.62){\gamma_i} \ar[ld]|-{\delta_i} \ar[dd] \ar[rd]^-{0} \\
X_s \ar[dd] & V_i \ar[l]^(.32){p_i} \ar[rr]_(.46){q_i} \ar[dd] & & \Spec k[T_1,T_2]/(T_1T_2) \ar[dd] \\
& & \Spec (R/\mathfrak{a}) \ar[lld]_(.6){\gamma'_i} \ar[ld]|-{\delta'_i} \ar[rd]^-{0} \\
Y \otimes_R (R/\mathfrak{a}) & V'_i \ar[l]^-{p'_i} \ar[rr]_(.47){q'_i} & & \Spec (R/\mathfrak{a})[T_1,T_2]/(T_1T_2) }}
\end{equation}
with all vertical squares Cartesian.

Put $\mathfrak{p} \defeq \mathfrak{m} \cap R$ and $k' \defeq \kappa(\mathfrak{p})$. Then $\mathfrak{a} \subseteq \mathfrak{p}$. Put $A_0 \defeq \Hensel{(R_{\mathfrak{p}})}$, $S_0 \defeq \Spec A_0$, $\hat{S} \defeq \Spec \widehat{A}_0$, $X_0 \defeq Y \times_{\Spec R} S_0$, $\hat{X} \defeq Y \times_{\Spec R} \hat{S}$. Let $s'$ be the point on $\Spec R$ defined by $\mathfrak{p}$. Then $s'$ is the image of $s$ under $S \to \Spec R$. Let $s_0$ and $\hat{s}$ be the closed points of $S_0$ and $\hat{S}$ respectively. As
\[\kappa(\hat{s}) = \kappa(s_0) = \kappa(s') = k' \,,\]
we may regard that
\[\hat{X}_{\hat{s}} = (X_0)_{s_0} = Y_{s'} \,.\]
 For each $i \in [1,n]$, $\gamma'_i$ induces a \mDash{k'}{rational} point $x'_i$ on $Y_{s'}$ which is the image of $x_i$ under $X \to Y$. By the bottom part of Diagram \eqref{E:curve2}, we have
 \begin{equation}\label{E:curve3}
 \widehat{\mathcal{O}}_{Y_{s'}, x'_i} \cong k'[[T_1,T_2]]\big/(T_1T_2) \,.
 \end{equation}
By the surjective morphism \eqref{E:curve1}, we know that $\Sequ{n}{x'}$ are all singular points of $Y_{s'}$ over $k'$. By Lemma \thmref{182P1} and Lemma \thmref{182P3},
\begin{equation}\label{E:curve4}
\widehat{\mathcal{O}}_{Y,x'_i} = \widehat{\mathcal{O}}_{\hat{X},x'_i} = \widehat{\mathcal{O}}_{X_0,x'_i} \cong \widehat{A}_0[[T_1,T_2]]/(T_1T_2 - a_i)
\end{equation}
for some $a_i \in \mathfrak{p}\widehat{A}_0$. Let $R'$ be the \mDash{R}{subalgebra} of $\widehat{A}_0$ generated by $\Sequ{n}{a}$. Put $\mathfrak{q} \defeq \mathfrak{p}\widehat{A}_0 \cap R'$, $T \defeq \Spec R'$, $Y' \defeq Y \times_{\Spec R} T$. Let $t \in T$ be the point defined by $\mathfrak{q}$. Then $\kappa(t) = k'$ and $Y'_t = Y_{s'}$. By \eqref{E:curve4} and Lemma \thmref{182P3}, we have
\[\widehat{\mathcal{O}}_{Y',x'_i} \cong \widehat{\mathcal{O}}_{T,t}[[T_1,T_2]]/(T_1T_2 - a_i) \,.\]
By Lemma \thmref{149T5}, there exists a \mDash{T}{scheme} $W'_i$ of finite type, two \'{e}tale \mDash{T}{morphisms} $p''_i \colon W'_i \to Y'$ and
\[q''_i \colon W'_i \to \Spec R'[T_1,T_2]/(T_1T_2 - a_i) \,,\]
a point $z'_i$ on $W'_i$ such that $p''_i(z'_i) = x'_i$, $q''_i(z'_i)$ is the point defined by the prime ideal generated by $\mathfrak{q} \cup \{T_1,T_2\}$, and $\kappa(z'_i) = k'$. Put $W_i \defeq W'_i \times_T \hat{S}$,
\begin{align*}
\hat{p}_i & \defeq p''_i \times \iDe_{\hat{S}} \colon W_i \to X_1 \,, \\
\hat{q}_i & \defeq q''_i \times \iDe_{\hat{S}} \colon W_i \to \Spec A_1[T_1,T_2]/(T_1T_2 - a_i) \,.
\end{align*}
As $z'_i \in W'_i$ and $\hat{s} \in \hat{S}$ both map onto $t \in T$, there is a point $z_i \in W_i$ which maps onto both $z'_i$ and $\hat{s}$. Then $\hat{p}_i(z_i) = x'_i$ and $\hat{q}_i(z_i)$ is the point defined by the prime ideal generated by $\mathfrak{p} \cup \{T_1,T_2\}$. Thus $W_i$ may be contracted to a local chart of $x'_i$. Therefore $\hat{X} \to \hat{S}$ satisfies \nA. By Lemma \thmref{182P2} $\hat{X} \to \hat{S}$ satisfies \nC. Let $\hat{D}_i$ denote the connected component of
\[\Spec \mathcal{O}_{W_i} \big/ \bigl((\hat{q}_i)^{\#}(T_1),(\hat{q}_i)^{\#}(T_2)\bigr)\]
containing $z_i$. Then $\hat{D}_i$ is \'{e}tale over $\Spec \widehat{A}_0/(a_i)$. Since $\kappa(z_i) = k'$ and $\widehat{A}_0/(a_i)$ is complete, a fortiori Henselian, by Lemma \thmref{668P19} $\hat{D}_i = \Spec \widehat{A}_0/(a_i)$. Thus $\hat{D}_i \to \hat{S}$ is a closed immersion. Since $\hat{X} \to \hat{S}$ is separated, the composite morphism
\[\hat{D}_i \hookrightarrow W_i \xrightarrow{p'''_i} \hat{X}\]
is a closed immersion. So we may regard $\hat{D}_i$ as a closed subscheme of $\hat{X}$. Since $\Sequ{n}{x'}$ are all singular points of $\hat{X}_{\hat{s}}$ over $k'$, we have
\[D(\hat{X}/\hat{S}) = \coprod^n_{i=1}\hat{D}_i \,.\]
Note that as subsets of $\hat{X}$, $\hat{D}_i \cap \hat{D}_j = \emptyset$ for all $1 \leqslant i < j \leqslant n$. Thus $D(\hat{X}/\hat{S})$ is also a closed subscheme of $\hat{X}$.

In the following we shall descend $a_i$ to elements in $A_0$. Put 
\[A_1 \defeq \widehat{\SH{(R_{\mathfrak{p}})}} \,, \qquad A_2 \defeq A_1 \otimes_{A_0} A_1 \,, \qquad A_3 \defeq A_1 \otimes_{A_0} A_1 \otimes_{A_0} A_1 \,.\]
For each $i \in [1,3]$, put $S_i \defeq \Spec A_i$ and $X_i \defeq Y \times_{\Spec R} S_i$. Obviously we may regard $\hat{A}_0$ as a subring of $A_1$. So there is a canonical morphism $S_1 \to \hat{S}$. Thus $X_1 \to S_1$ satisfies \nC. For each $i \in [1,n]$, put $D_i \defeq \hat{D}_i \times_{\hat{S}} S_1$. Then $\Sequ{n}{D}$ are all connected components of $D(X_1/S_1)$. Note that $S_2$ and $S_3$ might not be noetherian. We shall use the trick of inverse limits of schemes to avoid this difficult. By Theorem \thmref{672P5} and \cite[\S(8.6)]{AGro2}, there exists a finitely generated \mDash{A_0}{subalgebra} $A'$ of $A_1$ which contains $\Sequ{n}{a}$ such that if let $S' \defeq \Spec A'$ and $X' \defeq X_0 \times_{S_0} S'$, then $X' \to S'$ satisfies \nC; and closed subschemes $\Sequ{n}{D'}$ of $X'$ such that $D'_i \times_{S'} S_1 = D_i$ for all $i \in [1,n]$ and $D(X'/S') = \coprod\limits^n_{i=1}D'_i$.
Put
\begin{align*}
A'' & \defeq A' \otimes_{A_0} A' & S'' & \defeq \Spec A'' & X'' & \defeq X_0 \times_{S_0} S'' \\
A''' & \defeq A' \otimes_{A_0} A' \otimes_{A_0} A' & S''' & \defeq \Spec A''' & X''' & \defeq X_0 \times_{S_0} S''' \,.
\end{align*}
Let $g_1 \colon X_1 \to X'$, $g_2 \colon X_2 \to X''$, $g_3 \colon X_3 \to X'''$, $h_1 \colon S_1 \to S'$, $h_2 \colon S_2 \to S''$, $h_3 \colon S_3 \to S'''$ be the canonical morphisms. Then we have a commutative diagram with all squares Cartesian.
\[
\newcommand{\bM}{\ar@<.3ex>[rr] \ar@<-.3ex>[rr]}
\newcommand{\cM}{\ar@<.4ex>[rr] \ar[rr] \ar@<-.4ex>[rr]}
\xymatrix@C-1em@R-1ex{& X''' \cM \ar[dd] & & X'' \bM \ar[dd] & & X' \ar[rr] \ar[dd] & & X_0 \ar[dd] \\
X_3 \cM \ar[ur]|{g_3} \ar[dd] & & X_2 \bM \ar[ur]|{g_2} \ar[dd] & & X_1 \ar[ur]|{g_1} \ar[rr] \ar[dd] & & X_0 \ar[dd] \ar@{=}[ur] \\
& S''' \cM & & S'' \bM & & S' \ar[rr] & & S_0 \\
S_3 \cM \ar[ur]|{h_3} & & S_2 \bM \ar[ur]|{h_2} & & S_1 \ar[ur]|{h_1} \ar[rr] & & S_0 \ar@{=}[ur]}\]
Let $s_1$ be the closed point of $S_1$ and put $s' \defeq h_1(s_1) \in S'$. By Theorem \thmref{668P12}, we have
\[\prM{1}^{\ast}D(X'/S') = \prM{2}^{\ast}D(X'/S') = D(X''/S'') \,.\]
Pulling back to $X_2$, we have
\[\prM{1}^{\ast}D(X_1/S_1) = \prM{2}^{\ast}D(X_1/S_1) \,.\]
By \cite[VIII, 1.9]{AGro3}, there exists a closed subscheme $C$ of $X_0$ such that $D(X_1/S_1) = C \times_{X_0} X_1$. Let $\Sequ{n'}{C}$ be all connected components of $X_0$. By Lemma \thmref{668P4}, $n' = n$, and by rearranging the order of $\Sequ{n'}{C}$, we may assume that $C_i \times_{X_0} X_1 = D_i$ for all $i \in [1,n]$. By \cite[\S(8.6)]{AGro2} and by replacing $B$ with a suitably large finitely generated \mDash{A_0}{subalgebra} of $A_1$, we may assume that $C_i \times_{X_0} X' = D'_i$ for all $i \in [1,n]$. Then these data $X' \to S'$, $\Sequ{n}{D'}$, $\Sequ{n}{a}$ satisfy conditions in the begin of \S\ref{Se:loc}. Let $\mathscr{P}'$, $\mathscr{Q}'$, $\theta'$, $\vartheta'$ and $\mathfrak{d}'$ be the notations for $f'$ defined in Remark \thmref{668P8}; and let $(\mathscr{M}', \mathscr{N}', \sigma', \tau', \varphi')$ be a semistable log structure of $f'$. Let $\rho' \colon \mathbb{N}_{S'}^n \to \mathcal{O}_{S'}$ be the homomorphism of monoids defined by $\rho'(\varepsilon_i) = a_i$, where $\Sequ{n}{\varepsilon}$ is a basis of $\mathbb{N}^n$. Then there is a commutative diagram
\[\xymatrix{\mathbb{N}_{S'}^n \ar[r]^-{\rho'} \ar[d]_{\gamma} & \mathcal{O}_{S'} \ar[d] \\ \mathscr{Q} \ar[r]_-{\vartheta} & \mathcal{O}_{S'}/\mathcal{O}_{S'}^{\ast}}\]
where $\gamma \colon \mathbb{N}_{S'}^n \to \mathscr{Q}$ is the canonical morphism. As $\gamma_{\bar{s}'}$ is an isomorphism, by Lemma \thmref{670P1} there exists an affine \'{e}tale neighborhood $N$ of $\bar{s}'$ such that $\gamma|_N$ lifts to a chat $ \mathbb{N}_N^n \to \mathscr{N}'|_N$. Note that $A_1$ is a strictly Henselian local ring. By \cite[(18.8.1)]{AGro2}, $S_1 \to S'$ factors through $N$. So by replacing $S'$ with $N$, we obtains that $\gamma$ lifts to a chat $\rho \colon \mathbb{N}_{S'}^n \to \mathscr{N}'$ and the composite morphism
\[\mathbb{N}_{S'}^n \xrightarrow{\rho} \mathscr{N}' \to \mathcal{O}_{S'}\]
is equal to $\rho'$. By Theorem \thmref{672P7}, there is an isomorphism $\varphi \colon \prM{1}^{\ast}\mathscr{N}' \xrightarrow{\sim} \prM{2}^{\ast}\mathscr{N}'$ of log structures on $X''$ such that $\prM{23}^{\ast}(\varphi) \circ \prM{12}^{\ast}(\varphi) = \prM{13}(\varphi)$ on $X'''$. By Lemma \thmref{668P16}, both $\prM{1}^{\ast}(\rho)$ and $\prM{2}^{\ast}(\rho)$ are lifts of the canonical morphism $\mathbb{N}_{S''}^n \to \mathscr{Q}''$, where $\mathscr{Q}''$ is defined in Remark \thmref{668P8}. So there exists an element
\[u = (\Sequ{n}{u}) \in (\mathcal{O}_{S''}^{\ast})^n(S'')\]
such that $\varphi \circ \prM{1}^{\ast}(\rho) = u \cdot \prM{2}^{\ast}(\rho)$ and $\prM{23}^{\ast}(u_i) \circ \prM{12}^{\ast}(u_i) = \prM{13}(u_i)$ in $\mathcal{O}_{S'''}^{\ast}(S''')$ for all $i \in [1,n]$. We have $a_i \otimes 1 = u_i \cdot (1 \otimes a_i)$ in $A''$. Let $v_i$ denote the image of $u_i$ in $A_2^{\ast}$. Then $a_i \otimes 1 = v_i \cdot (1 \otimes a_i)$ in $A_2$ and $v_i$ defines an isomorphism $\psi_i \colon \mathcal{O}_{S_2} \xrightarrow{\sim} \mathcal{O}_{S_2}$ of \mDash{\mathcal{O}_{S_2}}{modules} such that $\prM{23}^{\ast}(\psi_i) \circ \prM{12}^{\ast}(\psi_i) = \prM{13}(\psi_i)$ on $S_3$. By flat descent of quasi-coherent sheaf, there exists an invertible \mDash{\mathcal{O}_S}{module} $\mathscr{L}_i$ and an isomorphism $\phi_i \colon q^{\ast}\mathscr{L}_i \xrightarrow{\sim} \mathcal{O}_{S_1}$ of \mDash{\mathcal{O}_{S_1}}{modules} such that
\[\xymatrix@C+2em{\prM{1}^{\ast}q^{\ast}\mathscr{L} \ar[r]^-{\prM{1}^{\ast}(\phi_i)} \ar@{=}[d] & \mathcal{O}_{S_2} \ar[d]^{\psi_i} \\ \prM{2}^{\ast}q^{\ast}\mathscr{L} \ar[r]_-{\prM{2}^{\ast}(\phi_i)} & \mathcal{O}_{S_2}}\]
is commutative, where $q \colon S_1 \to S_0$ is the canonical morphism. Since $A_0$ is a noetherian local ring, $\mathscr{L}_i \cong \mathcal{O}_{S_0}$. So $\phi_i$ defines an element $w_i \in A_1^{\ast}$ such that $v_i = w_i^{-1} \otimes w_i$. Put $b_i \defeq w_i a_i$. Then $b_i \otimes 1 = 1 \otimes b_i$ in $A_2$. By \cite[I, 2.18]{JSMil1}, $b_i \in A_0$. Since $A_1$ is flat over $\widehat{A}_0$, by Lemma \thmref{668P1}, $b_i = w'_ia_i$ for some $w'_i \in \widehat{A}_0^{\ast}$. Now replacing $T_1$ with $(w'_i)^{-1}T_1$ in \eqref{E:curve4}, we obtain
\begin{equation}\label{E:curve5}
\widehat{\mathcal{O}}_{X_0,x'_i} \cong \widehat{A}_0[[T_1,T_2]]/(T_1T_2 - b_i) \,.
\end{equation}
As $A_0 = \Hensel{(R_{\mathfrak{p}})}$, there exists a finitely generated \'{e}tale \mDash{R}{algebra} $B$, a prime ideal $\mathfrak{q}$ of $B$, elements $\Sequ{n}{c} \in \mathfrak{q}$, and an isomorphism $\nu \colon \Hensel{(B_{\mathfrak{q}})} \xrightarrow{\sim} A_0$ of \mDash{R}{algebras} such that $\nu(c_i) = b_i$ for all $i \in [1,n]$. Put $L \defeq \Spec B$ and $Z \defeq Y \times_{\Spec R} L$. Let $l \in L$ be the point defined by $\mathfrak{q}$. As $\kappa(l) = k'$, we may regard that $Z_l = (X_0)_{s_0}$. By \eqref{E:curve5} and Lemma \thmref{182P3}, we have
\[\widehat{\mathcal{O}}_{Z,x'_i} \cong \widehat{\mathcal{O}_{L,l}}[[T_1,T_2]]/(T_1T_2 - c_i) \,.\]
As $B$ is finitely generated over $\mathbb{Z}$, by Lemma \thmref{149T5} there exists a \lC of $Z \to L$ at $x'_i$. By base extension $S_0 \to L$, we obtains a \lC of $X_0 \to S_0$ at $x'_i$. So $X_0 \to S_0$ satisfies \nA. By Lemma \thmref{182P2} $X_0 \to S_0$ also satisfies \nC. By base extension $S \to S_0$, we know that $f \colon X \to S$ satisfies \nC and is a semistable curve over $S$.
\end{proof}

From above lemma and Corollary \thmref{672P6}, we obtains that

\begin{noteTh}
Any semistable curve over a locally noetherian scheme satisfies \nC, thus has a canonical semistable log structure.
\end{noteTh}

\begin{noteTh}
Let $S$ be a noetherian scheme and $f \colon X \to S$ be a proper and faithfully flat morphism. Then $X$ is a semistable curve over $S$ if and only if for every closed point $y \in S$, $X \times_S \Spec \kappa(y)_s \to \Spec \kappa(y)_s$ is a semistable curve.
\end{noteTh}

\begin{noteTh}
Let
\[\xymatrix{X' \ar[r]^{f'} \ar[d]_{p} \CartS & S' \ar[d]^{q} \\ X \ar[r]_{f} & S}\]
be a Cartesian square of schemes such that
\begin{enumerate}
\item $f'$ is a semistable curve,
\item $q$ is faithfully flat,
\item $f$ is proper and of finite presentation.
\end{enumerate}
Then $f$ is also a semistable curve.
\end{noteTh}

\begin{proof}
Obviously $f$ is also faithfully flat. So by Definition \thmref{182P8}, we may assume that $S = \Spec k$ where $k$ is a separably closed field, and $S'$ is affine. By Theorem \thmref{672P5}, we may further assume that $S'$ is of finite type over $S$. Then the theorem is by Theorem \thmref{181P6}.
\end{proof}

\clearP
\bibliographystyle{plain}
\bibliography{ref}

\end{document}